%     Ce ficher est du plain TeX.   
%     Partial Riemann Problem, Boundary Conditions and Gas dynamics 
%     Fran\c{c}ois Dubois, 29 juin 2000, 
%         edition  28, 29 decembre 2001,  08, 09 janvier, 
%                  28, 30 avril, 6, 7, 10  mai 2002.
%     version Linux du 10 mai  2002.
%     qques retouches le 11 sept 2005 et le 26 mars 2009.
%     version hal...   11 janvier 2011 
%%%%%%%%%%%%%%%%%%%%%%%%%%%%%%%%%%%%%%%%%%%%%%%%%%%%%%%%%%%%%%%%%% 
%                       thanks to Gabriel Turinici
\input epsfx.tex
%                       thanks to Gabriel Turinici
%%%%%%%%%%%%%%%%%%%%%%%%%%%%%%%%%%%%%%%%%%%%%%%%%%%%%%%%%%%%%%%%   

\overfullrule=0pt

% \nopagenumbers

%taille d'agrandissement  :
\magnification=1400
\hsize=12.2cm
\vsize=16.0cm
% marge gauche : \hoffset=-1cm
\hoffset=-0.4cm
\voffset=.8cm
\baselineskip 16 true pt

\font \smcaps=cmbx10 at 12 pt

% \font \pbf=cmb10   scaled 200

\font \caps=cmbx10 scaled 1200
\font \gcaps=cmbx10 scaled 1600
\font \smcaps=cmcsc10 
\font \pecaps=cmcsc10 at 9 pt

%  Pagination d'apres Raymond S\'eroul pages 231 et 68 
\newtoks \hautpagegauche  \hautpagegauche={\hfil}
\newtoks \hautpagedroite  \hautpagedroite={\hfil}
\newtoks \titregauche     \titregauche={\hfil}
\newtoks \titredroite     \titredroite={\hfil}
\newif \iftoppage         \toppagefalse   
\newif \ifbotpage         \botpagefalse    
\titregauche={\pecaps      Fran\c{c}ois Dubois } 
\titredroite={\pecaps     Partial Riemann problem, boundary  conditions, and gas dynamics }
\hautpagegauche = { \hfill \the \titregauche  \hfill  }
\hautpagedroite = { \hfill \the \titredroite  \hfill  }
\headline={ \vbox  { \line {  
\iftoppage    \ifodd  \pageno \the \hautpagedroite  \else \the
\hautpagegauche \fi \fi }     \bigskip  \bigskip  }}
\footline={ \vbox  {   \bigskip  \bigskip \line {  \ifbotpage  
\hfil {\oldstyle \folio} \hfil  \fi }}}

%petite boullette

%carre de fin de dmonstration
\def\sqr#1#2{{\vcenter{\vbox{\hrule height.#2pt \hbox{\vrule width .#2pt height#1pt 
\kern#1pt \vrule width.#2pt} \hrule height.#2pt}}}}
\def\square{\mathchoice\sqr64\sqr64\sqr{4.2}3\sqr33} 

% Les nombres mathmatiques

\def\R{{\rm I}\! {\rm R}}

% indice bas
\def\ib#1{_{_{_{\scriptstyle{#1}}}}}

% les valeurs absolues, les modules et les nornes
\def\abs#1{\mid \! #1 \! \mid }

\def\mod#1{\setbox1=\hbox{\kern 3pt{#1}\kern 3pt}%
\dimen1=\ht1 \advance\dimen1 by 0.1pt \dimen2=\dp1 \advance\dimen2 by 0.1pt
\setbox1=\hbox{\vrule height\dimen1 depth\dimen2\box1\vrule}%
\advance\dimen1 by .1pt \ht1=\dimen1
\advance \dimen2 by .01pt \dp1=\dimen2 \box1 \relax}

\def\nor#1{\setbox1=\hbox{\kern 3pt{#1}\kern 3pt}%
\dimen1=\ht1 \advance\dimen1 by 0.1pt \dimen2=\dp1 \advance\dimen2 by 0.1pt
\setbox1=\hbox{\kern 1pt  \vrule \kern 2pt \vrule height\dimen1 depth\dimen2\box1
\vrule
\kern 2pt \vrule \kern 1pt  }%
\advance\dimen1 by .1pt \ht1=\dimen1
\advance \dimen2 by .01pt \dp1=\dimen2 \box1 \relax}

\rm
\hyphenation {per-tur-ba-tion com-pres-sible com-pres-si-bles  me-cha-nics sprin-ger
con-ve-na-ble si-mu-la-tion si-tua-tion uni-di-men-sio-nal  com-pres-sible
com-pres-si-bles  me-cha-nics sprin-ger con-ve-na-ble si-mu-la-tion si-tua-tion 
com-pres-sible com-pres-si-bles  me-cha-nics sprin-ger con-ve-na-ble
si-mu-la-tion si-tua-tion uni-di-men-sio-nal thermo-dyna-miques }

\hyphenation {per-tur-ba-tion com-pres-sible com-pres-si-bles  me-cha-nics sprin-ger
con-ve-na-ble si-mu-la-tion si-tua-tion uni-di-men-sio-nal  com-pres-sible
com-pres-si-bles  me-cha-nics sprin-ger con-ve-na-ble si-mu-la-tion si-tua-tion 
com-pres-sible com-pres-si-bles  me-cha-nics sprin-ger con-ve-na-ble
si-mu-la-tion si-tua-tion uni-di-men-sio-nal thermo-dyna-miques ra-re-fac-tion 
pro-blem boun-da-ry}

% saut de ligne
\def\br {\break} 

%  page
  \def \page #1{\unskip\leaders\hbox to 1.3 mm {\hss.\hss}\hfill   {  $\,\,$ {\oldstyle #1}}}

%%%%%%%%%%%%%%%%%%%%%%%%%%%%%%%%%%%%%%%%%%%%%%%%%%%%%%%%%%%%%%%%%%%%%%%%%%%%%%%%%%%%

$\, $ 

\bigskip \bigskip \bigskip

\centerline {\gcaps   Partial Riemann problem,  }
\bigskip
\centerline {\gcaps    boundary  conditions, }
\bigskip
\centerline { \gcaps   and gas dynamics }
\null\vskip 1.1 cm 

\centerline{\caps Fran\c{c}ois Dubois }
\bigskip
\bigskip

\centerline { Conservatoire National des Arts et M\'etiers }

\centerline {  Department of Mathematics,   Paris, France. }

\bigskip
\centerline { Applications Scientifiques du Calcul Intensif }

\centerline {   b\^at. 506, BP 167,  F-91~403   Orsay  Cedex, France. } 

\centerline {   dubois@asci.fr  } 

\bigskip
\bigskip

\centerline { {\bf  June 29, 2000}\footnote {$ \, ^{^{\scriptstyle \bf \square}}$}
{Published in  {\it  Absorbing Boundaries and Layers, Domain Decomposition} 
\smallskip  \vskip -2pt  
{\it   Methods: Applications to Large Scale Computations},  pages 16-77,  
\smallskip  \vskip -2pt  
 Edited by   Lo\"{\i}c Tourrette and Laurence Halpern,  
\smallskip  \vskip -2pt 
Nova Science Publishers, Inc, New York,  2001,   ISBN: 1-56072-940-6.  
\smallskip  \vskip -2pt 
R\'evision 10 May  2002. Edition 11 January 2011.  }}

\bigskip \bigskip \bigskip \bigskip


\smallskip 


%%%%%%%%%%%%%%%%%%%%%%%%%%%%%%%%%%%%%%%%%%%%%%%%%%%%%%%%%%%%%%%%%%%%%%%%%%%%%%%%%%%%%%%%
\vfil \eject 
~

\bigskip
\bigskip 
\bigskip 
 \bigskip \noindent {\smcaps Contents}   
\bigskip 

\noindent  {\bf 0)  $ \,\, $    Introduction} \page {$\,\,$2}

\noindent {\bf 1)  $ \,\, $    Euler equations of gas dynamics} % \page {$\,\,$3}  

1.1 \quad  Thermodynamics \page {$\,\,$3}

1.2 \quad  Linear and nonlinear waves \page {$\,\,$6}

1.3 \quad  Riemann invariants and rarefaction waves \page {$\,\,$11}

1.4 \quad   Rankine-Hugoniot jump relations and shock waves  \page {$\,\,$14}

1.5 \quad    Contact discontinuities  \page {$\,\,$17}

1.6 \quad  Practical solution of the Riemann problem   \page {$\,\,$19}

\noindent  {\bf 2)  $ \,\, $    Partial Riemann problem for hyperbolic systems} 

2.1  \quad   Simple waves for an hyperbolic system of conservation laws   \page {$\,\,$24}

2.2  \quad   Classical Riemann problem between two states    \page {$\,\,$27}

2.3  \quad    Boundary manifold    \page {$\,\,$29}

2.4  \quad   Partial Riemann problem between a  state and a  manifold     \page {$\,\,$31}

2.5  \quad   Partial Riemann problem with an half-space   \page {$\,\,$35}

\noindent {\bf 3)  $ \,\, $    Nonlinear boundary conditions for gas dynamics }   

3.1  \quad   System of linearized Euler equations    \page {$\,\,$37}

3.2  \quad   Boundary  problem for linear hyperbolic systems  \page {$\,\,$38}

3.3  \quad    Weak Dirichlet nonlinear boundary condition  \page {$\,\,$40}

3.4  \quad    Some  fluid boundary conditions \page {$\,\,$42}

3.5  \quad    Rigid wall and moving solid boundary  \page {$\,\,$47}

\noindent {\bf 4)  $ \,\, $    Application to the finite volume method} 

4.1  \quad  Godunov finite-volume method  \page {$\,\,$51}

4.2  \quad   Boundary fluxes   \page {$\,\,$52}

4.3  \quad   Strong nonlinearity at the boundary  \page {$\,\,$53}

4.4  \quad   Extension to  second order accuracy and to two space
dimensions   {$\,\,$55} 

\noindent {\bf 5)  $ \,\, $    References. }  \page {$\,\,$61} 

%%%%%%%%%%%%%%%%%%%%%%%%%%%%%%%%%  fd le 15 aout 2005  %%%%%%%%%%%%%%%%%%%%
\toppagetrue  
\botpagetrue    
%%%%%%%%%%%%%%%%%%%%%%%%%%%%%%%%%  fd le 15 aout 2005  %%%%%%%%%%%%%%%%%%%%
 
%%%%%%%%%%%%%%%%%%%%%%%%%%%%%%%%%%%%%%%%%%%%%%%%%%%%%%%%%%%%%%%%%%%%%%%%%%%%%%%%%%%%%%%%
% \vfil \eject  
\bigskip  \bigskip  

\noindent  {\smcaps 0) $ \quad $  Introduction. }

\noindent   $\bullet \qquad \,\,\,\,\, $ 
We study in this article the problem of boundary conditions for the equations of gas
dynamics. We restrict ourselves to the model of perfect fluid and essentially to simple
unidimensional geometry. The problem can be considered from two points of view. First
the linearization of the Euler equations of gas dynamics allows the development of
the so-called method of characteristics and the linearized problem is mathematically 
well posed when boundary  data are associated with the characteristics that comes
inside the domain of study.  Second the  physical approach separates clearly fluid
boundary conditions where incomplete physical data are given, and rigid walls where the
 slip condition describes the interaction between the fluid and a given object
located inside or around the fluid. There is at our knowledge no simple way to
connect these  mathematical and physical approaches for strongly nonlinear
interactions.
 
\smallskip \noindent   $\bullet \qquad \,\,\,\,\, $ 
We introduce in this contribution the notion of partial Riemann problem.  Recall that
the Riemann problem describes a shock tube interaction  between two given states ; the
partial Riemann problem is a generalization of the previous concept and introduces the
notion of boundary manifold. In what follows, we first recall very classical notions
concerning gas dynamics and the associated Riemann problem. In a second part, we
introduce the partial Riemann problem for general systems of conservation laws and
proves that this problem admits a solution in some class of appropriate nonlinear
waves. In section 3, we recall the linearized analysis with the method of
characteristics, introduce the weak formulation of the Dirichlet boundary condition
for nonlinear situations in terms of the partial Riemann problem and show that lot of
physically relevant situations are described with this theoretical framework. In the
last paragraph, we   propose a practical implementation of the previous considerations
with the finite volume method. 

%%%%%%%%%%%%%%%%%%%%%%%%%%%%%%%%%%%%%%%%%%%%%%%%%%%%%%%%%%%%%%%%%%%%%%%%%%%%%%%%%%%%%%%%%%%%%%%%
\bigskip  \bigskip  
% \vfill \eject \noindent 
\noindent {\smcaps 1) $ \quad $  Euler equations of gas dynamics. }
\smallskip \noindent {\smcaps 1.1 } $ \,  $ { \bf  Thermodynamics. }

\noindent   $\bullet \qquad \,\,\,\,\, $ 
We study a perfect gas submitted to a motion with variable velocity in space and
time. We note first that  the primitive  unknowns of this problem are  the scalar
fields that  characterize the thermodynamics  of the gas, {\it i.e.}  density  $\,
\rho ,\, $ internal energy~$\, e ,\,$ temperature $\, T ,\,$  and pressure~$\,p \,$
(see {\it e.g.} the book of Callen [Ca85]).  In what follows, we suppose that the gas is a
polytropic perfect gas ; it has constant specific heats at constant volume $\, C_v \,$
and at constant pressure  $\, C_p .\,$ These two quantities do not depend on any
thermodynamic variable like temperature or pressure ; we denote by  $\, \gamma \,$
their ratio~: 

\noindent  (1.1.1) $\qquad \displaystyle 
\gamma \,\,= \,\, {{C_p}\over{C_v}} \,\, (= \,\, Cste) \,. \,$ 

\noindent
We suppose that  the gas satisfies the law of perfect gas that can be written
with the following form~: 

\noindent  (1.1.2) $\qquad \displaystyle 
p \,\,= \,\, (\gamma-1) \, \rho \, e \,.\, $

\noindent
As usual, internal energy and temperature are linked together  by the Joule-Thomson
relation~: 

\noindent  (1.1.3) $\qquad \displaystyle 
e \,\,= \,\, C_v \, T \,. \,$

\smallskip \noindent   $\bullet \qquad \,\,\, $ 
In the formalism proposed by Euler during the 18$^{0}$ century, the motion is
described with the help of an unknown vector field  $\, u \,$ which is a function of
space $ \,x \,$ and time ~$ \,t \,$~: 

\noindent  (1.1.4) $\qquad \displaystyle 
u \,\,= \,\, u(x,\,t) \,.\,$

\noindent
In the following, we will suppose that space $ \,x \,$ has only one dimemsion $\, (x
\in \R).$ We have four unknown functions (density, velocity, pressure and internal
energy) linked together by the state law (1.1.2) In consequence, we need three
complementary equations in order to define a unique solution of the problem. 

\smallskip \noindent   $\bullet \qquad \,\,\, $ 
The general laws of Physics assume that mass, momentum and total energy are conserved
quantities, at least in the context of classical physics associated to the paradigm
of  invariance for the  Galileo group of space-time transformations (see {\it e.g.} Landau
and Lifchitz [LL54]). When we write the conservation of mass, momentum and energy
inside an infinitesimal volume  $\, {\rm d}x \,$ advected with celerity   $\,
u(x,\,t) ,\,$ which is exactly the mean velocity of particules that compose the gas,
it is classical [LL54] to write the fundamental  conservation laws of Physics with the
help of divergence operators~: 

\smallskip \noindent  (1.1.5) $\qquad \displaystyle 
{{\partial \rho}\over{\partial t}} \,\,+\,\, {{\partial}\over{\partial x}}\,
\bigl( \rho \, u \bigr) \,\,= \,\, 0  \,$

\noindent  (1.1.6) $\qquad \displaystyle 
{{\partial  }\over{\partial t}} \bigl( \rho \, u \bigr)  \,\,+\,\,
{{\partial}\over{\partial x}}\, \bigl( \rho \, u^2 \,+\, p \bigr) \,\,= \,\, 0  \,$

\noindent  (1.1.7) $\qquad \displaystyle 
{{\partial}\over{\partial t}} \Bigl( {{1}\over{2}} \rho \, u^2 \,+\, \rho \,e \, 
\Bigr)  \,\,+\,\, {{\partial}\over{\partial x}}\,\Bigl( \, \bigl( {{1}\over{2}} \rho
\, u^2 \,+\,p\bigr) \, u \,+\, p \, u \,  \Bigr) \,\,= \,\, 0  \,. \,$

\smallskip  \noindent
We introduce the specific total energy  $\, E \,$ by unity of volume

\noindent  (1.1.8) $\qquad \displaystyle 
E \,\,= \,\, {{1}\over{2}}  u^2 \,+\, e \,$ 

\noindent
and the vector  $\, W \,$ composed by the ``conservative variables'' or more precisely
by the ``conserved variables''~: 

\noindent  (1.1.9) $\qquad \displaystyle 
W \,\, = \,\, \bigl( \, \rho \,,\, \rho \, u \,,\, \rho \, E \, \bigr)
^{\displaystyle \rm t} \,. \,$ 

\noindent 
The conservation laws (1.1.5)-(1.1.7) take the following general  form of a so-called
system of conservation laws~: 

\noindent  (1.1.10) $\qquad \displaystyle 
{{\partial W}\over{\partial t}} \,\,+\,\, {{\partial}\over{\partial x}}\, 
F(W)  \,\,= \,\, 0 \,$

\noindent
where the flux vector $\, W \longmapsto F(W) \,$ satisfies the following algebraic
expression~: 

\noindent  (1.1.11) $\qquad \displaystyle 
F(W)  \,\, = \,\, \bigl( \, \rho \,u \,,\, \rho \, u^2 \,+\,p  \,,\, \rho \, u \, E
\,+\, p \, u \, \bigr) ^{\displaystyle \rm t} \,. \,$

\smallskip \noindent   $\bullet \qquad \,\,\, $ 
The system of conservation laws (1.1.10) is an hyperbolic system of equations 
sometimes stated as a ``quasilinear system'' of conservation laws. From a mathematical
point of view, the actual state of knowledge (see {\it e.g.} Serre [Se96]) does not give any
answer to the question of the existence of a solution  $\,\, (x,\,t) \longmapsto
W(x,\,t)
\,\,$  when time $\, t \,$ is no longer  small enough, even if the initial condition 
$\, \, x \longmapsto  W(x,\,0) \, \,$ is a regular function of the space variable  $\,
x .\,$ Even if the problem is solved in the scalar case (Kruzkov [Kv70]), nonlinear
waves are present in system (1.1.10), they  can create discontinuities, and this
phenomenon makes   the general mathematical  study of system  (1.1.10) very tricky.
Moreover, the  unicity of irregular solutions of this kind of system is in general not
satisfied. The usual tentative to achieve uniqueness is associated with the
incorporation of the second principle of Thermodynamics. Carnot's  principle of
thermodynamics can be explicitly introduced with the specific entropy $\, s .\,$
Recall that for a perfect polytropic gas, the specific entropy is a function of
the state thermodynamic variables  $\, p \,$ and~$\, \rho \,$~: 

\noindent  (1.1.12) $\qquad \displaystyle 
s \,\,= \,\, C_v \,\, {\rm Log} \, \Bigl( 
{{\rho_0^{\gamma}}\over{p_0}} \, {{p}\over{\rho^{\gamma}}} \, \Bigr) \,$ 

\noindent 
and it is easy to see that the following so-called {\bf mathematical entropy} $\, \eta
({\scriptstyle \bullet}) \,$

\noindent  (1.1.13) $\qquad \displaystyle 
\eta(W) \,\,= \,\, - \rho \, s \,$ 

\noindent 
is a  convex function of the conservative variables if temperature is positive (we
refer {\it e.g.} to our proof  in [Du90]). It admits an entropy flux  $\, \xi(W) \,$ in the
sense of the mathematical theory proposed by Friedrichs and  Lax in 1971 [FL71] and we
have 

\noindent  (1.1.14) $\qquad \displaystyle 
\xi(W) \,\,= \,\, - \rho \,u\, s \,. \, $ 

\noindent
The second principle of increasing of physical entropy when time is increasing has
been formalized mathematically by Germain-Bader [GB53], Oleinik [Ol57], Godunov
[Go61] and  Friedrichs-Lax in 1971 [FL71] among others. It  takes inside the formalism
of mechanics of continuous media the following weak ``conservative form''~: 

\noindent  (1.1.15) $\qquad \displaystyle 
{{\partial }\over{\partial t}} \, \eta(W)  \,\,+\,\, {{\partial}\over{\partial
x}}\,\xi(W) \,\,\leq\,\, 0 \,. \,$

\noindent 
Inequality (1.1.15) is exactly an equality for a regular solution $\,\,W({\scriptstyle
\bullet} ,\, {\scriptstyle \bullet} ) \,\,$ of conservation law (1.1.10). It has to
be considered in the sense of distributions for a weak solution $\,\, W({\scriptstyle
\bullet} ,\, {\scriptstyle \bullet} )\,\,$ of conservation law (1.1.10). 

\smallskip \noindent   $\bullet \qquad \,\,\, $ 
The mathematical entropy (1.1.13) is a strictly convex function of the conservative
variables if temperature is positive. We  introduce the Fr\'echet derivative $\,\, {\rm
d}\eta(W) \,\,$  of entropy $\,\, \eta ({\scriptstyle \bullet} )  \,\,$  with respect 
to the conservative variables  $\, W \,$~: 

\noindent  (1.1.16) $\qquad \displaystyle 
{\rm d}\eta(W) \,\,= \,\, \phi \,\, {\rm d}W \,$ 

\noindent 
and define the so-called {\bf  entropy variables}  $\, \phi .\,$ The precise calculus
of entropy variables is elementary from the traditional expression of the second
principle of thermodynamics. If a container has  volume  $\, V ,\,$ mass $\, M ,\,$
internal energy $\, {\cal E} ,\,$ then entropy  $\, S \,$  is a {\bf function}
$\,\Sigma({\scriptstyle \bullet}) \,$ of the above variables~: 

\noindent  (1.1.17) $\qquad \displaystyle 
S \,\,= \,\, \Sigma \, (M,\,V,\, {\cal E}) \,,\, $

\noindent 
and moreover, function $\,\Sigma({\scriptstyle \bullet}) \,$ is an extensive  function
(homogeneous of degree 1) of the preceeding variables, {\it id est } 
$\,\,  \Sigma \, (\lambda \, M,\,\lambda \, V,\, \lambda \, {\cal E}) \,=\,
\lambda \, \Sigma \, (M,\,V,\, {\cal E}) \,,$\br
$ \forall \,\,\, \lambda \, > 0 \,. \,
$  The classical differential relation between these quantities can be stated as

\noindent  (1.1.18) $\qquad \displaystyle 
{\rm d}{\cal E} \,\, =\,\, T \, {\rm d}S \,\,-\,\,p \, {\rm d}V \,\,+\,\, \mu \, 
{\rm d}M \,.\,$

\noindent 
Notice that the third intensive thermodynamic variable $\, \mu \,$ is just the {\bf 
specific chemical  potential}, {\it i.e.} the chemical potential by unit of mass. 
Then the global values are related for corresponding specific quantities with the help
of the following relations~: 

\noindent  (1.1.19) $\qquad \displaystyle 
M \,\, =\,\, \rho \, V \,,\,\quad {\cal E} \,\,= \,\, e \,M \,,\, \quad S \,\,=\,\, s
\,M \,\,,\, \quad  {\cal E} \,+\, p \, V \,-\, T \,  S \,\,= \,\, \mu \, M \,.\, $

\noindent
Due to the homogeneity of function  $\, \Sigma({\scriptstyle \bullet}) ,\,$  the
mathematical entropy  $\,\eta({\scriptstyle \bullet}) \,$ can be simply expressed with
the help of this thermostatic primitive  function~: 

\noindent  (1.1.20) $\qquad \displaystyle 
\eta(W) \,\, = \,\, - \Sigma \, (\rho \,,\, 1 \,,\, \rho \, e ) \,.\,$ 

\noindent
Taking into account the relations (1.1.8) and (1.1.19), we easily differentiate the
relation (1.1.20) and it comes 

\noindent  (1.1.21) $\qquad \displaystyle 
{\rm d}\eta(W) \,\,=\,\, {{1}\over{T}}\, \Bigl( \mu - {{u^2}\over{2}} \Bigr) \,{\rm
d}\rho \,\,+ \,\, {{u}\over{T}} \,{\rm d}(\rho \,u) \,\, -\,\, {{1}\over{T}} \, {\rm
d}(\rho \,E) \,. \,$   

\noindent
By comparison between (1.1.16) and (1.1.21) we have 

\noindent  (1.1.22) $\qquad \displaystyle 
\phi \,\,=\,\, {{1}\over{T}} \, \Bigl( \, \mu - {{u^2}\over{2}}  \,,\, u \,,\, -1
\, \Bigr) \,. $ 

\bigskip 
% \vfill \eject 
\noindent {\smcaps 1.2 } $ \,  $ { \bf 	Linear and nonlinear waves. }

\noindent   $\bullet \qquad \,\,\, $ 
In this section, we construct particular solutions of the so-called Riemann problem.
First recall that the Riemann problem consists in searching an entropy solution ({\it
i.e.} a solution satisfying (1.1.15) in the sense of distributions, see {\it e.g.}
Godlewski-Raviart [GR96] or Serre [Se96]) $\,\,  W(x,\,t)\,\,\,  (x\in\R,\,\, t > 0)
\,\, $ of the following Cauchy  problem~: 

\noindent  (1.2.1) $\qquad \displaystyle 
{{\partial W}\over{\partial t}} \,\,+\,\, {{\partial}\over{\partial x}}\, F(W)  \,\,=
\,\, 0 \,$

\setbox21=\hbox {$\displaystyle  W_l \,\,, \qquad x \,< \, 0 \,  $}
\setbox22=\hbox {$\displaystyle  W_r \,\,, \qquad x \,> \, 0 \,.\,   $}
\setbox30= \vbox {\halign{#&# \cr \box21 \cr \box22    \cr   }}
\setbox31= \hbox{ $\vcenter {\box30} $}
\setbox44=\hbox{\noindent  (1.2.2) $\displaystyle  \qquad   W(x \,,\,0) \,\,= \,\,
\left\{ \box31 \right. $}  

\noindent $ \box44 $

\noindent
We first remark that conservation law (1.2.1) is {\bf invariant under a change of
the space-time scale}.  Space-time transform $\,T_{\lambda} \,$ parameterized by  $\,
\lambda > 0, \,$ and  defined by the relation  $\,\,  T_{\lambda}(x,\,t) = (\lambda
\, x,\, \lambda \, t ) \,\, $ can be applied on any   (weak) solution $\,\, 
W({\scriptstyle \bullet}\,,\, {\scriptstyle \bullet}) \,\,$  of system 
(1.2.1)-(1.2.2) and generates a new weak solution  $\,\, T_{\lambda}\, W(x,\,t) \,\,$
defined by the condition  $\,\, (T_{\lambda}\, W)(x,\,t) \equiv W(
T_{\lambda}(x,\,t))\,.\,$ We remark also that the initial condition (1.2.2) is
invariant  under space dilatation, {\it i.e.} 

\noindent  (1.2.3) $\qquad \displaystyle 
W(\lambda \,x \,,\, 0) \,\,= \,\, W(x\,,\,0) \,\,, \qquad \forall \, \lambda \,> \, 0
\, . \,$

\noindent 
Then an hypothesis of {\bf unicity} for weak solution of problem (1.2.1)-(1.2.2) shows
that the solution  $\,\, W({\scriptstyle \bullet}\,,\, {\scriptstyle \bullet})
\,\,$ must be self-similar. This property can be  expressed by the relation 

\noindent  (1.2.4) $\qquad \displaystyle 
W(\lambda \,x \,,\, \lambda \,t) \,\,= \,\, W(x\,,\,t) \,\,, \qquad \forall \, \lambda
\,> \, 0 \,. \,$

\noindent
The relation (1.2.4) claims that solution  $\, W(x,\,t) \,$ must be searched under
a selfsimilar   form, {\it i.e.} under  the form of a function of the  variable
$\,{{x}\over{t}}\,$ (see more details {\it e.g.} in the book of  Landau and Lifchitz
[LL54])~:

\noindent  (1.2.5) $\qquad \displaystyle 
W(x \,,\, t)  \,\,=\,\,   U(\xi)  	\,\,, \qquad  \xi \,=\, {{x}\over{t}} \,. \,$

\smallskip \noindent   $\bullet \qquad \,\,\, $ 
As a   first step, we consider regular solutions  $\,  \xi \longmapsto  U(\xi)  \,$
of the Riemann problem (1.2.1) (1.2.2).  We introduce the representation (1.2.5) inside
conservation law (1.2.1) and we obtain by this way~: 

\noindent  (1.2.6) $\qquad \displaystyle 
{\rm d}F(U(\xi)) \, {\scriptstyle \bullet} \, {{{\rm d}U}\over{{\rm d}\xi}} \,\, =\,\,
\xi
\,\,   {{{\rm d}U}\over{{\rm d}\xi}} \, .\,$

\noindent
We deduce from (1.2.6) that we are necessarily  in one of the two following cases~:
either vector   $\, {{{\rm d}U}\over{{\rm d}\xi}} \,$ is equal to zero or this vector
is not equal to zero. In the first case, the solution is  a constant state and in the
second opportunity, the vector  $\, {{{\rm d}U}\over{{\rm d}\xi}} \,$ is necessarily
equal to some eigenvector   $\, R(W)\,$  of the jacobian matrix  $\, {\rm d}F(W) .\,$
In this second case, we have the classical relation between the jacobian matrix,
eigenvector and eigenvalue $\,\,  \lambda (W) \,\,$~: 

\noindent  (1.2.7) $\qquad \displaystyle 
{\rm d}F(W) \, {\scriptstyle \bullet} \,  R(W) \,\,=\,\, \lambda (W) \,   R(W) \, $

\noindent 
in the particular case where  $\,\, W = U(\xi) .\,$  By identification between the
two relations (1.2.6) and (1.2.7), we deduce that the vectors $\,
{{{\rm d}U}\over{{\rm d}\xi}} \,$  and   $\, R(U)\,$ are proportional and we deduce
also~:  

\noindent  (1.2.8) $\qquad \displaystyle 
\lambda\bigl(   U(\xi) \bigr) \,\,= \,\, \xi \,.$

\noindent 
We have derived the conditions that characterize a so-called {\bf rarefaction wave.} 

\smallskip \noindent   $\bullet \qquad \,\,\, $ 
It is also possible to suppose that function  $\,  \xi \longmapsto U(\xi) \,$  
admits some point of discontinuity at  $\,  \xi  =  \sigma .\, $ It satisfies
necessarily the Rankine-Hugoniot relations that link  the jump of state $\, [ \, U\,
], \,$ the jump of the flux $\, [ \,  F(U)\, ] \,$ and the celerity
$\, \sigma  \,$ of the discontinuity profile~: 

\noindent  (1.2.9) $\qquad \displaystyle 
\bigl[ \,  F(U)\, \bigr] \,\,= \,\, \sigma  \, \bigl[ \, U\, \bigr] \,. \,$

\noindent 
Recall that the Rankine-Hugoniot relations express that the discontinuous function 
$\,  \xi  \longmapsto U(\xi) \,$ is a weak solution of conservation (1.2.1) (see {\it e.g.}
the book of Godlewski and Raviart [GR96]). By this way, we are deriving a so-called
{\bf shock wave}. 

\smallskip \noindent   $\bullet \qquad \,\,\, $ 
We have observed in the previous section that the construction of a rarefaction wave
is associated with some eigenvalue of the jacobian matrix   $\, {\rm d}F(W). \,$ In
the case of the Euler equations of gas dynamics, these eigenvalues are simply
computed with the help of the so-called nonconservative form of the equations that 
are obtained with the introduction  of specific entropy  $\, s \,$  (see [LL54] for the
proof)~:  

\noindent  (1.2.10) $\qquad \displaystyle 
{{\partial \rho }\over{\partial t}} \,\,+\,\, u \, {{\partial \rho }\over{\partial
x}} \,\,+\,\, \rho \, {{\partial u }\over{\partial x}}  \,\,\,\,= \,\, 0 \,$

\noindent  (1.2.11) $\qquad \displaystyle 
{{\partial u }\over{\partial t}} \,\,+\,\, u \, {{\partial u }\over{\partial
x}} \,\,+\,\, {{1}\over{\rho}} \, {{\partial p }\over{\partial x}}  \,\,= \,\, 0 \,$

\noindent  (1.2.12) $\qquad \displaystyle 
{{\partial s }\over{\partial t}} \,\,+\,\, u \, {{\partial s }\over{\partial
x}} \qquad  \qquad \, \,\,\, = \,\, 0 \,. \,$

\noindent 
Then we set~: 

\noindent  (1.2.13) $\qquad \displaystyle 
Z (W) \,\,\equiv\,\,  \bigl(  \, \rho \,,\, u \,,\, s \,\bigr )^{\displaystyle \rm t}
\,$

\noindent  (1.2.14) $\qquad \displaystyle 
B(Z)  \,\,=\,\, \pmatrix{  u & \rho & 0 \cr  \displaystyle  {{c^2}\over{ \rho}} & u&
\displaystyle {{1}\over{\rho}} \,   {{\partial p}\over{\partial s}}(\rho,\,s) \cr 0 & 0
& u \cr} \,$

\noindent
and the Euler equations that took the  form (1.2.10)-(1.2.12) can be also written with
the above matrix~:  

\smallskip \noindent  (1.2.15) $\qquad \displaystyle 
{{\partial Z}\over{\partial t}} \,\,+\,\, B(Z) \, {\scriptstyle \bullet}\, 
{{\partial Z}\over{\partial x}}  \,\,= \,\, 0 \, . \, $

\noindent 
Celerity $\, c \,$ for  sound waves  used in formula (1.2.14) is defined
by 

\noindent  (1.2.16) $\qquad \displaystyle  
c^2  \,\,= \,\,   {{\partial p}\over{\partial \rho }}(\rho,\,s) \,$ 

\noindent 
and for a perfect polytropic gas with ratio  $\, \gamma \,$ of specific heats, we
have~: 

\noindent  (1.2.17) $\qquad \displaystyle 
c  \,\,= \,\,  \sqrt{ {{\gamma \, p }\over{\rho}} } \,.\,$
 
\smallskip \noindent   $\bullet \qquad \,\,\, $ 
If we wish to diagonalize the matrix  $\, {\rm d}F(W)  ,\,$ it is sufficient to make
this job for matrix  $\, B(Z)\,$ because they are {\bf conjugate} as we show in the
following. We introduce the vector variable  $\,Z\,$ inside equation (1.2.1)~: 
 
\noindent  (1.2.18) $\qquad \displaystyle
{\rm d}W(Z)  \, {\scriptstyle \bullet}\, {{\partial Z}\over{\partial t}} \,\,+\,\,
{\rm d}F(W) \, {\scriptstyle \bullet}\,  {\rm d}W(Z) \, {\scriptstyle \bullet}\,
{{\partial Z}\over{\partial x}}  \,\,= \,\, 0 \,$

\noindent 
and by comparison with relation (1.2.15) we have necessarily~: 
 
\noindent  (1.2.19) $\qquad \displaystyle
B(Z)  \,\,=\,\,   ({\rm d}W(Z))^{-1} \, {\scriptstyle \bullet}\, {\rm d}F(W) \,
{\scriptstyle \bullet}\, {\rm d}W(Z) \,.\,$

\noindent 
This conjugation relation shows that  $\, \widetilde{R}(Z) \equiv  ({\rm d}W(Z))^{-1}
\, {\scriptstyle \bullet}\, R(W(Z)) \,$ is an eigenvector for the matrix  $\,
B(Z(W))\,$  with the eigenvalue  $\, \lambda(W) \,$    if   $\, R(W) \,$ is some
eigenvector of the jacobian matrix  $\, {\rm d}F(W) \,$ that satisfies the relation
(1.2.7). The diagonalization of matrix  $\, B(Z)\, $ is straightforward. We find 

\noindent  (1.2.20) $\qquad \displaystyle
\lambda_1 \,\,= \,\, u-c \quad < \quad \lambda_2 \,\,= \,\, u \quad < \quad \lambda_3
\,\,= \,\, u+c \,$

\noindent  
with associated eigenvectors given by the following formulae~: 

\smallskip \noindent  (1.2.21) $\qquad \displaystyle
\widetilde{R}_1(Z) \,\,=\,\, \pmatrix {\rho \cr -c \cr 0 \cr } \,\,,\quad 
\widetilde{R}_2(Z) \,\,=\,\, \pmatrix { {{\partial p}\over{\partial s}}  \cr 0 \cr
-c^2 \cr }   \,\,,\quad \widetilde{R}_3(Z) \,\,=\,\, \pmatrix {\rho \cr c \cr 0 \cr }
\,. \,$

\noindent 
We remark that the derivation of scalar field  $\, \lambda_1 \,$  (respectively 
$\, \lambda_3 $) inside direction  $\, \widetilde{R}_1 \,$  (respectively $\, 
\widetilde{R}_3 $) is never null 

\noindent  (1.2.22) $\qquad \displaystyle
{\rm d}  \lambda_1 (W(Z)) \, {\scriptstyle \bullet}\, \widetilde{R}_1(W) \,\,\neq \,\,0
\,\,,\qquad  {\rm d}  \lambda_3 (W(Z)) \, {\scriptstyle \bullet}\, \widetilde{R}_3(W)
\,\,\neq \,\,0 \,\,,\qquad \forall \, W \,$

\noindent
whereas we have the opposite situation for the second eigenvalue $\, \lambda_2 \,$~:  

\noindent  (1.2.23) $\qquad \displaystyle
{\rm d}  \lambda_2 (W(Z)) \, {\scriptstyle \bullet}\, \widetilde{R}_2(W) \,\,= \,\,0
\,\,,\qquad \forall \, W \,. \,$

\noindent  
For this reason, we will say that the eigenvalues  $\,  \lambda_1  \,$ and $ \,
\lambda_3 \,$ define {\bf genuinely nonlinear} fields whereas the eigenvalue $\,\, 
 \lambda_2 = u \, \,$  defines a {\bf linearly degenerate} field. Rarefaction and
shock  waves are always associated with genuinely nonlinear fields and in what follows,
we distinguish between   1-rarefaction wave and 3-rarefaction wave. 

% \smallskip \smallskip \vskip 3,5cm  \smallskip  \qquad  \qquad   
% \special{illustration  fig.1.1.epsf scaled  600}  \smallskip  \smallskip 
\bigskip 
\centerline { \epsfysize=4,5cm    \epsfbox  {fig.1.1.epsf} }
\smallskip  \smallskip

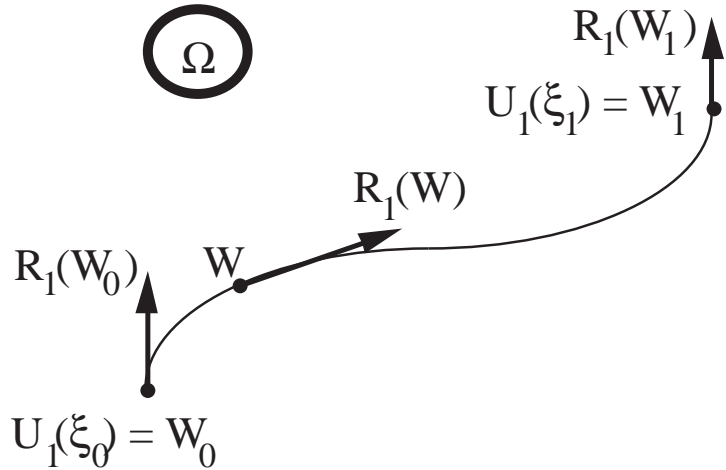
\noindent {\bf Figure 1.1}	\quad {\it Rarefaction wave associated with eigenvalue  
$\, \lambda_1 \,=\, u - c \,$   in the space of states ; the curve is everywhere
tangent to eigenvector $\, R_1(W) .\,$ }
\smallskip  \smallskip 

% \bigskip \smallskip \smallskip \vskip 3,5cm  \smallskip  \qquad  \qquad   
% \special{illustration  fig.1.2.epsf scaled  600}  \smallskip  \smallskip 
\bigskip 
\centerline { \epsfysize=4,5cm    \epsfbox  {fig.1.2.epsf} }
\smallskip  \smallskip

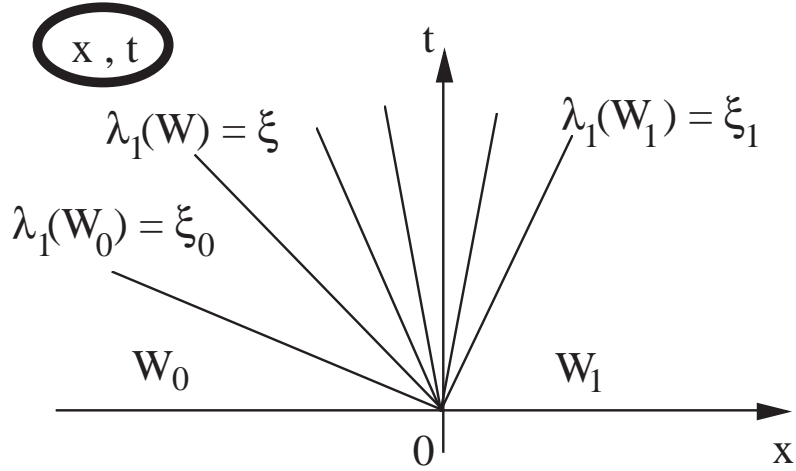
\noindent {\bf Figure 1.2} \quad  {\it 	Rarefaction wave $\, U_1(\xi) \,$   associated
with eigenvalue   $\, \lambda_1 \,=\, u - c \,$   inside the space-time plane. }
\smallskip  \smallskip 

\smallskip \noindent   $\bullet \qquad \,\,\, $ 
A 1-rarefaction wave is a function  $\, \xi \longmapsto U_1(\xi) \,$ that satisfies 

\noindent  (1.2.24) $\qquad \displaystyle
{{{\rm d}}\over{{\rm d}\xi}} \bigl( U_1(\xi) \bigr) \,\,$ proportional to $ 
\displaystyle R_1 \bigl( U_1(\xi) \bigr)  \,\,  $

\noindent 
and we have an analogous definition for a 3-wave   $\, \xi \longmapsto U_3(\xi) \,$~: 
$\,\, {{{\rm d}}\over{{\rm d}\xi}} \bigl( U_3(\xi) \bigr) \,\,$ is proportional to $  
R_3 \bigl( U_3(\xi) \bigr)  .\,  $ Then if we integrate the vector field
$\, \widetilde{R}_1(Z) \, $ {\it i.e.} if we solve the ordinary differential equation
$\,\,\,   {{{\rm d}Z}\over{{\rm d}\xi}} = \widetilde{R}_1(Z(\xi)) \, ,\, \,$ we find a
function  $\, \, \xi \longmapsto  Z(U_1(\xi)) \, \,$ that defines a 1-rarefaction wave,
say for variable  $\, \xi \,$  contained between two limiting values  $\, \xi_0 \,$ 
and  $\, \xi_1 . \,$ We can draw the solution of differential equation (1.2.24) inside
the space of states   $\,W\,$  ; we find a curve $\,\, \xi \longmapsto U_1(\xi) \,\,$ 
that satisfies equation (1.2.24) and the initial condition 

\noindent  (1.2.25) $\qquad \displaystyle
U_1(\xi_0)  \,\,=\,\,   W_0 \,.\, $ 

\noindent 
Moreover this curve   is defined up to  $\, \xi = \xi_1 \,$  where the final state 
$\,  W_1 \,$ is achieved (see Figure 1.1) 

\noindent  (1.2.26) $\qquad \displaystyle
U_1(\xi_1)  \,\,=\,\,   W_1 \,. \, $ 

\noindent 
In space-time  $\, (x,\, t) ,\,$  we observe that we have to consider three rates of
flow for celerities  $\, \xi = {{x}\over{t}} \,$~: 

\noindent  (1.2.27) $\qquad \displaystyle
U_1(\xi) \,\,=\,\,   W_0		\quad  \quad {\rm for }  \quad \quad \xi \, \leq \, \xi_0
\qquad
\qquad  $ (constant state)

\noindent  (1.2.28) $\qquad \displaystyle
U_1(\xi)  \,\, $ variable  \quad  for $  \quad  \,\,\, \xi_0 \,\leq \, \xi \, \leq \,
\xi_1 \quad \,\, \,$ (rarefaction wave)

\noindent  (1.2.29) $\qquad \displaystyle
U_1(\xi) \,\,=\,\,   W_1		\quad  \quad {\rm for } \quad  \quad  \xi \, \geq \, \xi_1
\qquad \qquad  $ (constant state)

\noindent  
as proposed at Figure 1.2. Notice that relation (1.2.8) imposes a particular value 
$\, \xi_0 \,$  for celerity~: 

\noindent  (1.2.30) $\qquad \displaystyle
\xi_0 \,\,=\,\,   u(W_0)  -  c(W_0) \,$ 

\noindent
and an analogous value for end-point   $\, \xi_1 .\,$  It is simple to show (see {\it e.g.}
Smoller [Sm83] or Godlewski-Raviart [GR96]) that inequality 

\noindent  (1.2.31) $\qquad \displaystyle
 \xi_0  \,\, \leq \,\,  \xi_1 \,$

\noindent 
is necessary if $\, W_1 \,$ is linked to  $\, W_0 \,$  through a -1- or a
3-rarefaction wave. In the particular case of a 3-rarefaction wave, relations (1.2.27)
to (1.2.29) are still correct, but relation (1.2.30) must be replaced by 

\noindent  (1.2.32) $\qquad \displaystyle
\xi_0 \,\,=\,\,   u(W_0)  +  c(W_0) \,$ 

\noindent  
and we have also an analogous relation for   $\, \xi_1 .\,$

\bigskip \bigskip
% \vfill \eject 
\noindent {\smcaps 1.3 } $ \,  $ {\bf Riemann invariants and rarefaction waves.  }

\noindent   $\bullet \qquad \,\,\, $ 
The practical computation of the curve $\,\, \xi \longmapsto U_1(\xi) \,\,$ satisfying
the relations (1.2.24) and (1.2.25) uses the notion of Riemann invariant.  By
definition, a 1-Riemann invariant  (respectively a 3-Riemann invariant)  is a
function  $\, \, W \longmapsto \beta^1(W) \,\, $   (respectively  $\, \, W
\longmapsto  \beta^3(W) $)  that is constant along the curves of ~1-rarefactions 
(respectively along the curves of ~3-rarefactions) and satisfies by definition~: 

\noindent  (1.3.1) $\qquad \displaystyle
{\rm d}  \beta^1(W) \, { \scriptstyle \bullet } \, R_1(W) \,\,=\,\, 0 \qquad \forall
\, W \,$ 

\noindent 
(respectively   $\, {\rm d}  \beta^3(W) \, { \scriptstyle \bullet } \, R_3(W) \,=\, 
0 \,$   for each state   $W$). If we express this relation in terms of
nonconservative variables $Z$, we set  $\,\, \widetilde{\beta}(Z) \equiv
\beta\bigl(W(Z)\bigr) \,$ and we have~: 

\noindent  $\displaystyle
{\rm d}\widetilde{\beta}(Z) \, {\scriptstyle \bullet}\,  {\widetilde R}(Z)
\,\,=\,\, {\rm d} \beta (W(Z)) \, {\scriptstyle \bullet}\, {\rm d}W(Z) \,
{\scriptstyle \bullet}\, \bigl( {\rm d}W(Z) \bigr)^{-1} \, {\scriptstyle \bullet}\,
R(W(Z)) \,\,=\,\, 0 \,\, \,$ if $\, \beta( {\scriptstyle \bullet}) \,$ is a Riemann
invariant for the field associated with eigenvector $\,\, R ({\scriptstyle \bullet}) .
\,$ We deduce from the previous calculus the   relation~: 
 
\noindent  (1.3.2) $\qquad \displaystyle
{\rm d}  \beta^1(W(Z)) \, { \scriptstyle \bullet } \, {\widetilde R}_1(Z) \,\,=\,\, 0
\qquad \forall \, Z \, . \,$

\noindent
Taking into account the particular form (1.2.21) of the vector $\, 
{\widetilde R}_1(Z),\,$ the two following functions 

\noindent  (1.3.3) $\qquad \displaystyle
\beta^1_1(W) \,\,=\,\, s \,$

\noindent  (1.3.4) $\qquad \displaystyle
\beta^1_2(W) \,\,=\,\, u \,+\, \int_{ \displaystyle \rho_0}^{ \displaystyle \rho} \,
{{c(\theta,\,s)}\over{\theta}} \, {\rm d}\theta \,\,\,\,=\,\,\,\, u \,+\, {{2 \,
c}\over{\gamma \!-\! 1}}  \,  $

\noindent 
are particular ~1-Riemann invariants ; they satisfy together the relation (1.3.2).
The expression of the ~3-Riemann invariants is obtained by an analogous way~:

\smallskip \noindent  (1.3.5) $\qquad \displaystyle
\beta^3_1(W) \,\,=\,\, s \,$
\smallskip \noindent  (1.3.6) $\qquad \displaystyle
\beta^3_2(W) \,\,=\,\, u \,-\, \int_{ \displaystyle \rho_0}^{ \displaystyle \rho} \,
{{c(\theta,\,s)}\over{\theta}} \, {\rm d}\theta \,\,\,\,=\,\,\,\, u \,-\, {{2 \,
c}\over{\gamma \!-\! 1}}  \, . \,  $

\noindent 
The states  $\,W \,$  on a ~1-rarefaction wave issued from state  $\,W_0 \,$ are
explicited  by using the relation (1.3.1) for the Riemann invariants
proposed at relations (1.3.3) and (1.3.4).  This fact express that on
a 1-rarefaction wave, the two associated Riemann invariants are constant ;  we have 

\noindent  (1.3.7) $\qquad \displaystyle
s   \,\,=\,\, s_0 \,$

\noindent  (1.3.8) $\qquad \displaystyle
u \,+\, {{2 \, c}\over{\gamma \!-\! 1}} \,\,=\,\,  u_0 \,+\, {{2 \, c_0}\over{\gamma
\!-\! 1}}  \,$

\noindent  (1.3.9) $\qquad \displaystyle
\xi \,\,=\,\, u-c \,\,, \qquad \xi \, \geq \, u_0 - c_0 \,. $

\smallskip \noindent   $\bullet \qquad \,\,\, $  
We detail now the particular algebraic form that takes  the description of the link
beween a state   $\,W \,$  and its initial datum   $\,W_0 \,$  through a
~1-rarefaction wave inside a perfect polytropic gas. We first set (we refer for the
details to the book of Courant and Friedrichs [CF48]) that inside a ~1-rarefaction
wave, the pressure $\,p \,$ is a decreasing function of velocity~: 

\noindent  (1.3.10) $\qquad \displaystyle
p  \,\, \leq \,\, p_0 \,\,, \qquad u \,\, \geq \,\, u_0 \,\,,\qquad W \,\, $   issued
from $\,\, W_0 \,\,$     {\it via} a 1-rarefaction

\noindent
and due to the expression of the entropy $\,\,s \,\,$ for a polytropic perfect gas, 

\noindent  (1.3.11) $\qquad \displaystyle
s \,\,$ is a function of variable  $ \displaystyle \,\, {{p}\over{\rho^{\gamma}}} \,$

\noindent
it comes, taking into account the relations (1.2.17), (1.3.6) and (1.3.7),

\noindent  (1.3.12) $\qquad \displaystyle
u  \,-\,   u_0 \,+\, {{\sqrt{1 \!-\! \mu^4}}\over{\mu^2}} \,
{{p_0^{{\scriptstyle 1}\over{\scriptstyle 2   \gamma}}}\over{\sqrt{\rho_0}}} \, 
\biggl( p^{{\scriptstyle \gamma \!-\! 1 }\over{\scriptstyle 2   \gamma}} \,- \,
p_0^{{\scriptstyle \gamma \!-\! 1 }\over{\scriptstyle 2   \gamma}}  \biggr) \,\,
=\,\, 0 \qquad $ (1-rarefaction) 

\noindent  
parameterized by the non-dimensional coefficient   $\,\mu > 0 \,$ [note that this
parameter  $\, \mu \,$ has nothing to do with the chemical potential, even if the same
letter is used !] defined by 

\noindent  (1.3.13) $\qquad \displaystyle
\mu^2 \,\,=\,\, {{ \gamma \!-\! 1} \over{ \gamma \!+\! 1}} \,. \,$

\noindent 
In the plane  $\, (u, p) \,$ composed by velocity and pressure, the graphic
representation of the ~1-rarefaction (relation (1.3.12) under condition (1.3.10)) is
proposed on Figure~1.3.

% \bigskip \smallskip \smallskip \vskip 4,0cm  \smallskip  \qquad  \qquad   
% \special{illustration  fig.1.3.epsf scaled  600}  \smallskip  \smallskip 
% \bigskip 
\centerline { \epsfysize=5,0cm    \epsfbox  {fig.1.3.epsf} }
\smallskip  \smallskip

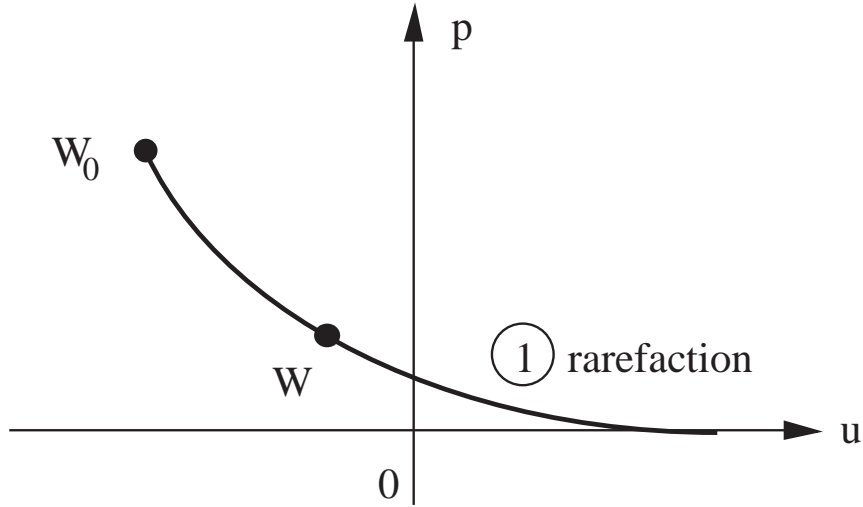
\noindent   {\bf Figure 1.3}	\quad  {\it  1-rarefaction wave linking state $\, W_0 \,$ 
with state  $\, W \,$   in the plane of velocity and pressure. }
\smallskip  \smallskip 

\smallskip \noindent   $\bullet \qquad \,\,\, $  
For a ~3-rarefaction wave, we look for (due to reasons that will be explained in
what follows) an upstream state  $\, W\,$ linked to a {\bf downstream} state  $\,
W_0 .\,$ We obtain 

\noindent  (1.3.14) $\qquad \displaystyle
s   \,\,=\,\, s_0 \,$

\noindent  (1.3.15) $\qquad \displaystyle
u \,-\, {{2 \, c}\over{\gamma \!-\! 1}} \,\,=\,\,  u_0 \,-\, {{2 \, c_0}\over{\gamma
\!-\! 1}}  \,$

\noindent  (1.3.16) $\qquad \displaystyle
\xi \,\,=\,\, u + c \,\,, \qquad \xi \, \leq \, u_0  +  c_0 \,. \,  $

\noindent 
It is also easy to prove that in a ~3-rarefaction,  velocity is a nondecreasing
function of pressure, {\it i.e.}~: 

\noindent  (1.3.17) $\qquad \displaystyle
p  \,\, \leq \,\, p_0 \,\,, \qquad u \,\, \leq \,\, u_0 \,\,,\qquad W_0 \,\, $   issued
from $\,\, W \,\,$     {\it via}  a ~3-rarefaction. 

\noindent 
A computation analogous to the one presented for ~1-rarefactions shows 

\noindent  (1.3.18) $\qquad \displaystyle
u  \,-\,   u_0 \,-\, {{\sqrt{1 \!-\! \mu^4}}\over{\mu^2}} \,
{{p_0^{{\scriptstyle 1}\over{\scriptstyle 2   \gamma}}}\over{\sqrt{\rho_0}}} \, 
\biggl( p^{{\scriptstyle \gamma \!-\! 1 }\over{\scriptstyle 2   \gamma}} \,- \,
p_0^{{\scriptstyle \gamma \!-\! 1 }\over{\scriptstyle 2   \gamma}}  \biggr) \,\,
=\,\, 0 \qquad $ (3-rarefaction).  

\noindent
The comparison between relations (1.3.12) and (1.3.18) induces us to set 

\noindent  (1.3.19) $\qquad \displaystyle
\psi( p \,;\, \rho_0 \,,\, p_0 \,;\, \gamma ) \,\,\equiv \,\,  {{\sqrt{1 \!-\!
\mu^4}}\over{\mu^2}} \, {{p_0^{{\scriptstyle 1}\over{\scriptstyle 2  
\gamma}}}\over{\sqrt{\rho_0}}} \,  \biggl( p^{{\scriptstyle \gamma \!-\! 1
}\over{\scriptstyle 2   \gamma}} \,- \, p_0^{{\scriptstyle \gamma \!-\! 1 }\over
{\scriptstyle 2   \gamma}}  \biggr) \,. \, $

\noindent
The graph  in the plane  $\, (u,\, p) \,$  of the curve of ~3-rarefaction (equation
(1.3.18) under the constraint (1.3.17)) is presented on Figure 1.4.

% \bigskip \smallskip \smallskip \vskip 4cm  \smallskip  \qquad  \qquad   
% \special{illustration  fig.1.4.epsf scaled  600}  \smallskip  \smallskip 
\bigskip 
\centerline { \epsfysize=5,0cm    \epsfbox  {fig.1.4.epsf} }
\smallskip  \smallskip

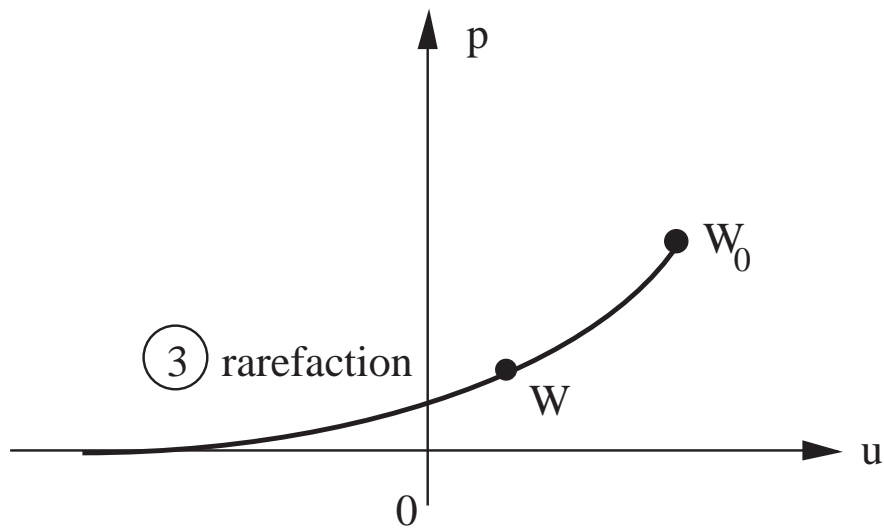
\noindent   {\bf Figure 1.4}	\quad  {\it  Curve in the plane of velocity and pressure
showing the set   of states $\, W \,$ such that the 3-rarefaction wave that begins at
state   $\, W \,$     ends at the particular state $\, W_0 .\,$  }
\smallskip  \smallskip 

\smallskip \noindent   $\bullet \qquad \,\,\, $  
It can be usefull to precise the (varying) state inside a ~1-rarefaction wave as a
function of celerity  $\, \xi = {{x}\over{t}} .\,$ Taking into account the relation
(1.2.8), we have on one hand~: 

\noindent  (1.3.20) $\qquad \displaystyle
u \,-\,c \,\,= \,\, \xi \,$

\noindent
and on the other hand taking into account the relation (1.3.8), it comes easily, for 
$\, u_l - c_l \,< \, \xi \,< \, u_1 - c_1 \,$~: 

\noindent  (1.3.21) $\qquad  \displaystyle
u \,\, = \,\, {{\gamma \!-\!1}\over{\gamma \!+\!1}} \, u_l \,\,+\,\, {{2}\over{\gamma
\!+\!1}} \, c_l \,+\,  {{2}\over{\gamma \!+\!1}} \, \xi \,$ 

\noindent  (1.3.22) $\qquad \displaystyle
c \,\, = \,\, {{\gamma \!-\!1}\over{\gamma \!+\!1}} \, u_l \,\,+\,\, {{2}\over{\gamma
\!+\!1}} \, c_l \,-\,  {{\gamma \!-\!1}\over{\gamma \!+\!1}} \, \xi \,.  \,$ 

\noindent
The elimination of pressure among the relations (1.2.17) and (1.3.7) allows the
evaluation of density. It is the same set of operations for a ~3-rarefaction wave.
The relation (1.2.8) gives the expression of the celerity of the 3-wave as a function
of the data~: 

\noindent  (1.3.23) $\qquad \displaystyle
u \,+\,c \,\,= \,\, \xi \,$

\noindent 
and the Riemann invariant (1.3.15) allows the evaluation of velocity  $\,u\,$  and
sound celerity $\,c\,$  as a function of the data~: 

\noindent  (1.3.24) $\qquad  \displaystyle
u \,\, = \,\, {{\gamma \!-\!1}\over{\gamma \!+\!1}} \, u_r \,\,-\,\, {{2}\over{\gamma
\!+\!1}} \, c_r \,+\,  {{2}\over{\gamma \!+\!1}} \, \xi \,$ 

\noindent  (1.3.25) $\qquad \displaystyle
c \,\, = \,\,- {{\gamma \!-\!1}\over{\gamma \!+\!1}} \, u_r \,\,+\,\, {{2}\over{\gamma
\!+\!1}} \, c_r \,+\,   {{\gamma \!-\!1}\over{\gamma \!+\!1}} \, \xi \,$ 

\noindent 
and under the conditions  $\, u_2 + c_2 \,< \, \xi \,< \,  u_r + c_r  . \,$

\bigskip \noindent {\smcaps 1.4 } $ \,  $ { \bf 	 Rankine-Hugoniot jump relations and
shock waves.  }

\noindent   $\bullet \qquad \,\,\, $ 
We consider now two particular states   $\,W_0\,$   and  $\,W\,$  linked together
{\it via} the Rankine-Hugoniot jump relations (1.2.9) and the associated entropy
inequality. We first remark that the Galilean invariance of the equations of gas
dynamics allows us to think the Physics inside the reference frame with a velocity
exactly equal to the celerity   $\, \sigma \,$  of the discontinuity. Then the jump 
equations (1.2.9) are written again under the more detailed form 

\noindent  (1.4.1) $\qquad \displaystyle
\bigl[ \, \rho \,(u-\sigma) \, \bigr] \,\,= \,\, 0 \,$

\noindent  (1.4.2) $\qquad \displaystyle
\bigl[ \, \rho \,(u-\sigma)^2 \,+ \, p  \bigr] \,\,= \,\, 0 \,$

\noindent  (1.4.3) $\qquad \displaystyle
\biggl[ \, \rho \,(u-\sigma)\,\Bigl(\, e \,+\, {{(u-\sigma)^2}\over{2}} \, \Bigr)
 \,+ \, p\,(u-\sigma)  \biggr] \,\,= \,\, 0 \,$

\noindent  
as it is also easy to derive directly. We consider the classical expression for mass
flux that cross  the shock wave~: 

\noindent  (1.4.4) $\qquad \displaystyle
m  \,\,=\,\,  \rho \,(u-\sigma)  \, . \,$

\noindent 
The associated conditions of increasing  physical entropy through a shock wave
(see {\it e.g.} [CF48] or [LL54])  claim that we have  

\noindent  (1.4.5) $\qquad \displaystyle
m \,\,> \,\,  0 \qquad $ 	through a 1-shock

\noindent  (1.4.6) $\qquad \displaystyle
m \,\,< \,\,  0 \qquad $ 		through a 3-shock.

\noindent 
The numbering of shock waves can be explained by an argument of continuity as
follows. If the jumps inside relations   (1.4.1) to (1.4.3)  are weak enough, it is
possible to show (see details in [Sm83] or [GR96] for example) that the celerity  
$\, \sigma \,$  of the discontinuity relative to both states $\, W_0 \,$  and  $\, W
\,$  has a limit value equal to the common value $\, u-c  \,$ for a ~1-shock and  
 $\, u+c \,$ for a ~3-shock wave. The particular case where the mass flux $\, m \,$ is
equal to zero will be considered afterwards and corresponds to a   slip surface or
 {\bf contact discontinuity}. 
 
\smallskip  \noindent   $\bullet \qquad \,\,\, $  
We detail the algebra that is necessary in order to express that the state  $\,W\,$ is
obtained from an upstream state  $\,W_0 \,$ through a ~1-shock wave with celerity 
$\, \sigma .\,$ We first remark that the entropy condition implies the following
family  of inequalities~: 

\noindent  (1.4.7) $\qquad \displaystyle
\rho \,\,> \,\,  \rho_0 \,\,,\qquad p \,\,> \,\, p_0  \,\,,\qquad  u-c \,\, < \,\,
\sigma \,\, < \,\, u_0 - c_0  \,\,,\qquad s \,\,> \,\, s_0 \,$

\noindent
when  $\,W\,$ is issued from $\,W_0 \,$ through a ~1-shock wave. We denote also by 
$\,h\,$   and   $\, \tau \,$  the specific enthalpy and the specific volume 

\noindent  (1.4.8) $\qquad \displaystyle 
h  \,\,=\,\,   e  \,+\, {{p}\over{\rho}} \,$

\noindent  (1.4.9) $\qquad \displaystyle 
\tau \,\,=\,\, {{1}\over{\rho}} \,.$

\noindent  
From relations (1.4.1) to (1.4.3) we deduce, taking into account (1.4.4), (1.4.8) and
(1.4.9)~: 

\noindent  (1.4.10) $\qquad \displaystyle 
\bigl[ \, u \, \bigr] \,\,= \,\, m \, \bigl[ \, \tau \, \bigr] \,$

\noindent  (1.4.11) $\qquad \displaystyle 
m^2 \,  \bigl[ \, \tau \, \bigr] \,+ \,  \bigl[ \, p \, \bigr] \,\,=\,\, 0 \,$

\noindent  (1.4.12) $\qquad \displaystyle 
\bigl[ \, h\, \bigr] \,+\, {{m^2}\over{2}} \,  \bigl[ \, \tau^2 \, \bigr] \,\,=\,\,0
\,.\,$

\noindent
In the particular case of a polytropic perfect gas, we eliminate the mass flow
$\,m\,$  from relations (1.4.11) and (1.4.12) and we express the enthalpy as a
function of pressure and specific volume, {\it i.e.} 

\noindent  (1.4.13) $\qquad \displaystyle 
h \,\,=\,\, {{\gamma}\over{\gamma \!-\!1}} \, p \, \tau \,. \,$

\noindent  
After some lines of elementary algebra we  express the specific volume $\, \tau \,$
downstream to the shock as a function of upstream data, downstream pressure and
variable   $\, \mu \,$  introduced at relation (1.3.13). We obtain~: 

\noindent  (1.4.14) $\qquad \displaystyle 
\tau \,\,=\,\,{{p_0 \,+\, \mu^2 \,p}\over{p \,+\, \mu^2 \,p_0}} \, \tau_0 \,. \,$

\noindent
We report this particular expression inside relation (1.4.11)   and obtain by this way
the square of the mass flux across the shock~: 

\noindent  (1.4.15) $\qquad \displaystyle 
m^2 \,\,=\,\, {{p \,+\, \mu^2 \,p_0}\over{(1\!-\!\mu^2)\, \tau_0}} \,\,, \qquad W \,$
issued from  $\, W_0 \,$   by a ~1-shock.

%%%%%%%%%%%%%%%%%%%%%%%%%%%%%%%%%%%%%%%%%%%%%%%%%%%%%%%%%%%%%%%%%%%%%%%%%%%%%%%%%%% 
% \bigskip  
\centerline { \epsfysize=6,0cm    \epsfbox  {fig.1.5.epsf} }
\smallskip  \smallskip

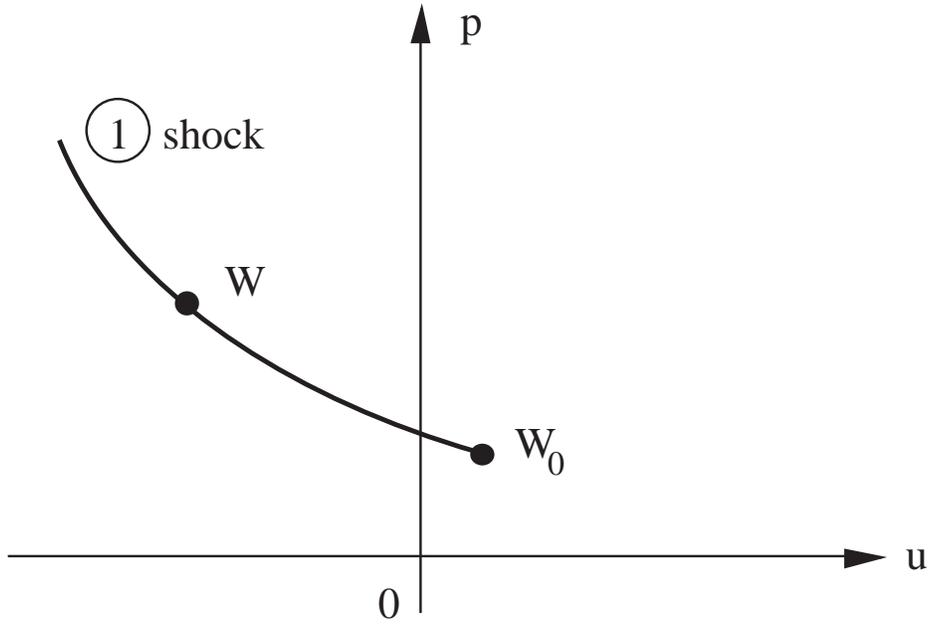
\noindent   {\bf Figure 1.5}	\quad  {\it   Curve showing (in the velocity-pressure 
plane) the set of states  $ \, W \, $  issued from the particular state  $ \, W_0 \,$  
through a  1-shock wave. }
\smallskip  \smallskip 
%%%%%%%%%%%%%%%%%%%%%%%%%%%%%%%%%%%%%%%%%%%%%%%%%%%%%%%%%%%%%%%%%%%%%%%%%%%%%%%%%%%
%%%%%%%%%%%%%%%%%%%%%%%%%%%%%%%%%%%%%%%%%%%%%%%%%%%%%%%%%%%%%%%%%%%%%%%%%%%%%%%%%%% 
% \bigskip 
\centerline { \epsfysize=5,0cm    \epsfbox  {fig.1.6.epsf} }
\smallskip  \smallskip

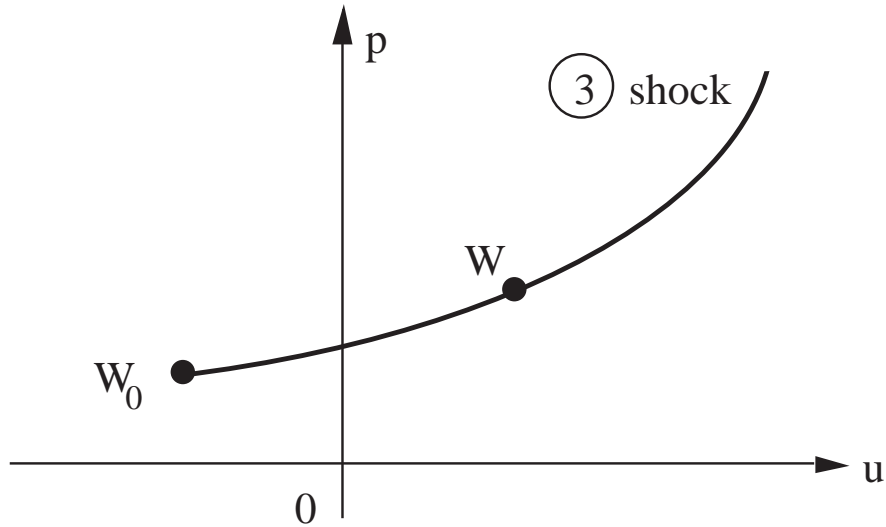
\noindent   {\bf Figure 1.6}	\quad  {\it   Set of states  $ \, W\, $    such that 
the 3-shock wave links  state  $ \, W \, $    to the particular state  $ \, W_0 .\,$  }
 \smallskip  \smallskip 
%%%%%%%%%%%%%%%%%%%%%%%%%%%%%%%%%%%%%%%%%%%%%%%%%%%%%%%%%%%%%%%%%%%%%%%%%%%%%%%%%%%

\smallskip  \noindent   $\bullet \qquad \,\,\, $  
For a ~1-shock wave, the relation (1.4.2) joined with (1.4.10) allows us to precise
the jump of velocity as a function of pressure and upstream state~: 

\noindent  (1.4.16) $\qquad \displaystyle
u \,-\, u_0 \,+\, \sqrt{{{1\!-\!\mu^2}\over{\rho_0 \,(p +\mu^2 \,p_0)}}} \, \bigl(
p-p_0 \bigr) \,\,=\,\, 0 \,$

\noindent 
when  $\,W\,$ is issued from $\,W_0 \,$ through a ~1-shock wave. In the plane of
velocity-pressure variables, we can propose (see Figure 1.5) the curve 
characterized by the equation (1.4.16) and the inequalities (1.4.7).  

\smallskip  \noindent   $\bullet \qquad \,\,\, $  
For a ~3-shock wave, the relations (1.4.10) to (1.4.13) are still valid. We do not
consider in what follows  a state  $\,W\,$ issued from  $\,W_0\,$ {\it via} a
~3-shock wave but we reverse the roles and consider the set of states  $\,W\,$ that
are {\bf  upstream}  to the particular state   $\,W_0\,$   through a ~3-shock wave. We
 reverse the zero index in previous formulae. With these particular hypotheses
concerning our new notations, a simple calculus shows that relation (1.4.15) remains 
still correct and we have 

\noindent  (1.4.17) $\qquad \displaystyle 
m^2 \,\,=\,\, {{p \,+\, \mu^2 \,p_0}\over{(1\!-\!\mu^2)\, \tau_0}} \,\,, \qquad  W_0
 \,$ issued from  $ \, W \,$   by a ~3-shock.
 
\noindent 
The algebraic relation that express the jump of velocity as a function of downstream
state  $\, W_0 \,$ and upstream pressure takes now the form 

\noindent  (1.4.18) $\qquad \displaystyle
u \,-\, u_0 \,-\, \sqrt{{{1\!-\!\mu^2}\over{\rho_0 \,(p +\mu^2 \,p_0)}}} \, \bigl(
p-p_0 \bigr) \,\,=\,\, 0 \,$

\noindent
when $\,W_0 \,$  is issued from $\,W \,$ through a ~3-shock wave. Taking into account
these new notations, the inequalities for entropy condition can be written as 

\smallskip \noindent  (1.4.19) $\qquad \displaystyle
\rho_0 \,\,<\,\, \rho \,\,,\qquad p_0 \,\,< \,\, p \,\,,\qquad u_0 + c_0 \,\,< \,\,
\sigma \,\,< \,\,u + c \,\,, \qquad s_0 \,\,< \,\,s \,$

\noindent
if $\,W_0 \,$  is issued from $\,W \,$ through a ~3-shock wave. A ~3-shock wave is
graphically represented represented at Figure 1.6. By comparison between the relations
(1.4.16) and (1.4.18), it is natural to define 

\noindent  (1.4.20) $\qquad \displaystyle
\varphi \, (   p \,;\, \rho_0 \,,\, p_0 \,;\, \gamma ) \,\,\equiv \,\, 
\sqrt{{{1\!-\!\mu^2} \over{\rho_0 \,(p +\mu^2 \,p_0)}}} \, \bigl( p-p_0 \bigr) \,\,$ 

\noindent
where parameters  $\,\mu\,$   and   $\, \gamma \,$   are linked by relation  
(1.3.13).

\bigskip \noindent {\smcaps 1.5} $ \,  $ { \bf Contact discontinuities. }

\noindent   $\bullet \qquad \,\,\, $  
When we have studied the rarefaction waves, we have introduced the notion of linearly
degenerated field at relation (1.2.23) and we have observed that  the second
characteristic field of the Euler equations of gas dynamics is effectively linearly
degenerated. This fact means that the second eigenvalue $\, \lambda_2 \,$ is also a
Riemann invariant~:

\noindent  (1.5.1) $\qquad \displaystyle
\beta^2_1(W) \,\,=\,\, u \,. \,$

\noindent 
Taking into account the particular expression (1.2.21) of eigenvectors, it is easy to
see that the pressure is also a ~2-Riemann invariant independent from the first
one~:  

\noindent  (1.5.2) $\qquad \displaystyle
\beta^2_2(W) \,\,=\,\, p \,.  \,$

\noindent
We can search a selfsimilar regular  wave $\, \, \xi \longmapsto  U(\xi) \, \,$
solution of the differential equation that expresses proportionality between 
$\, {{{\rm d}U} \over{{\rm d}\xi}} \,$    and   $\, R_2(U), \,$ but in that case
the necessary condition (1.2.8) implies that variable  $\, \xi \,$   is not allowed to
get any variation~! We have in fact by derivation of relation (1.2.8) 

\noindent  (1.5.3) $\qquad \displaystyle
{\rm d}  \lambda_2 (U) \, {\scriptstyle \bullet}\,  {{{\rm d}U} \over{{\rm d}\xi}}
\,\,= \,\, 1 \,$

\noindent 
whereas this last quantity is identically equal to zero (relation (1.2.23)) if we
consider the hypothesis of linearly degeneration of second characteristic  field $\, 
\lambda_2 \equiv u .\,$  

% \bigskip \smallskip \smallskip \vskip 3,5cm  \smallskip  \qquad  \qquad   
% \special{illustration  fig.1.7.epsf scaled  600}  \smallskip  \smallskip 
\bigskip 
\centerline { \epsfysize=4,5cm    \epsfbox  {fig.1.7.epsf} }
\smallskip  \smallskip

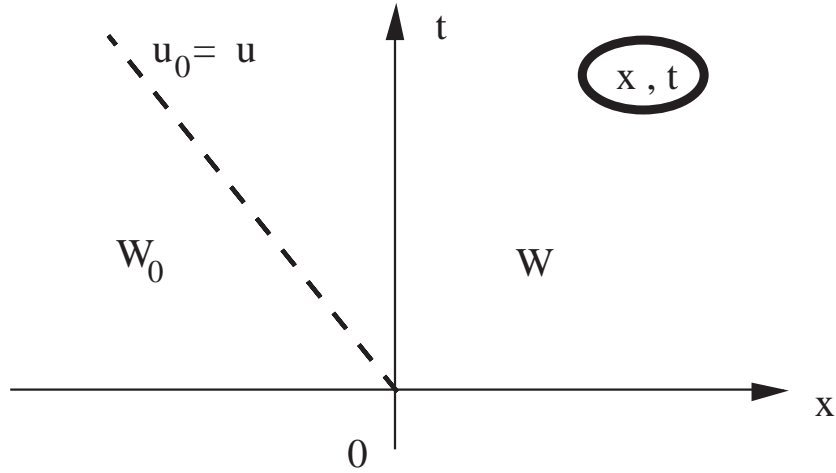
\noindent   {\bf Figure 1.7}	\quad  {\it  Contact discontinuity between the two 
states   $ \, W_0 \,$   and $ \, W .\, $   }
\smallskip  \smallskip 

\smallskip \noindent   $\bullet \qquad \,\,\, $  
We remark  also that  if the function  $\, \xi \longmapsto  U(\xi) \,$   is equal to
some integral curve of vector field $\, R_2 ,\,$  {\it i.e.} 

\noindent  (1.5.4) $\qquad \displaystyle
{{{\rm d}U} \over{{\rm d}\xi}} \,\,= \,\,  R_2(U) \,$ 

\noindent 
then we have between states  $\, W_0 \,$ and  $\, W \,$ the following calculus~: 

\noindent  $ \displaystyle
F(W)  \,-\,  F(W_0)	\,\,\,=\,\, \int_{\displaystyle \xi_0}^{\displaystyle \xi} \, {\rm
d}F(U(\eta)) \, {\scriptstyle \bullet} \,  {{{\rm d}U} \over{{\rm d}\eta}} \, {\rm
d}\eta \,\, = \,\, \int_{\displaystyle \xi_0}^{\displaystyle \xi} \, {\rm
d}F(U(\eta)) \, {\scriptstyle \bullet} \,  R_2(U(\eta))  \, {\rm d}\eta \,$

\noindent  $ \displaystyle \qquad \qquad \qquad \quad \,\,\, \,\,
= \,\, \int_{\displaystyle \xi_0}^{\displaystyle \xi} \,\lambda_2(U(\eta))  \, 
R_2(U(\eta)) \, {\rm d}\eta \,\,$

\noindent  $ \displaystyle \qquad \qquad \qquad \quad \,\,\, \,\,
= \,\, \lambda_2   \,  \int_{\displaystyle \xi_0}^{\displaystyle \xi} \,  
R_2(U(\eta))  \, {\rm d}\eta \,\,\,=\,\,\, \lambda_2 \,(W - W_0) \, . \,$

\noindent 
The fourth step in the previous relations  is a consequence of the linear
degenerescence of velocity  $\, u \,$   which is a constant along the ~2-curve defined
at relation (1.5.4). We have just shown that we have a Rankine-Hugoniot jump 
relation  between states $\, W_0 \,$ and $\, W \,$~:
 
\noindent  (1.5.5) $\qquad \displaystyle 
\bigl[ \,  F(W)\, \bigr] \,\,= \,\,  \lambda_2   \, \bigl[ \, W\, \bigr] \,. \,$

\noindent
In the space of states, the rarefaction waves defined by the differential relation
(1.5.4) and the shock curves that we founded at relation (1.5.5) are identical ! The
second characteristic field is linearly degenerate and defines a curve in the space
of states that joins the two states  $\, W_0 \,$  and   $\, W \,$ according to the
relations 

\noindent  (1.5.6) $\qquad \displaystyle 
u \,\,= \,\, u_0 \,$

\noindent  (1.5.7) $\qquad \displaystyle 
p \,\,= \,\, p_0  \,. \,$

\noindent 
In space-time variables, these two states are linked through a contact discontinuity
with celerity  $\, \sigma \,$ equal to  $\, \lambda_2 \,$ due to relation (1.5.5)~: 

\noindent  (1.5.8) $\qquad \displaystyle 
\sigma \,\,=\,\, u \,\,= \,\, u_0 \,. \,$

\noindent 
Such a wave is called {\bf contact discontinuity} or {\bf slip line}. We have 

\setbox21=\hbox {$\displaystyle  W_0 \,\,, \qquad \xi \,< \,  u_0 \,\,=\,\, \sigma 
\,  $}
\setbox22=\hbox {$\displaystyle  W \,\,, \,\, \qquad \xi \,> \, u_0 \,\,=\,\, \sigma 
\, . \,  $}
\setbox30= \vbox {\halign{#&# \cr \box21 \cr \box22    \cr   }}
\setbox31= \hbox{ $\vcenter {\box30} $}
\setbox44=\hbox{\noindent  (1.5.9) $\displaystyle  \qquad   U_2(\xi)  \,\,= \,\,
\left\{ \box31 \right. $}  

\noindent $ \box44 $

\noindent
A representation in space-time of relation (1.5.9) is proposed at Figure 1.7. In the
velocity-pressure plane, the (nonlinear !) projections of states  $\,W_0\,$  and  
$\,W\,$ coincide, as previously established in relations (1.5.6) and (1.5.7). 

\bigskip \noindent {\smcaps 1.6} $ \,  $ { \bf  Practical solution of the Riemann
problem. }

\noindent   $\bullet \qquad \,\,\, $  
We propose to solve the Riemann problem (1.2.1)-(1.2.2) between two states  $\,W_l \,$
and   $\, W_r ,\,$ and  we remark that the general theory proposed by Lax (see {\it e.g.} 
[Lax73]) can be applied for gas dynamics. We search two intermediate states $\,W_1 \,$ 
and  $\,W_2 \,$ such that 

\noindent  (1.6.1) $\qquad \displaystyle 
W_1 \,\, $    is issued from state    $\,W_l \,$  by a ~1-wave 

\noindent  (1.6.2) $\qquad \displaystyle 
W_2 \,\, $   is issued from state    $\,W_1 \,$  by a ~2-wave 

\noindent  (1.6.3) $\qquad \displaystyle 
W_r \,\, $   is issued from state   $\,W_2 \,$  by a ~3-wave.

\noindent 
As a first step, we restrict ourselves to the search of velocity and pressure which is
common to both states   $\,W_1 \,$ and  $\,W_2 \,\,$ due to the fact that the ~2-wave
is a contact discontinuity~: 

\noindent  (1.6.4) $\qquad \displaystyle 
u_1 \,\, = \,\, u_2 \,\, = \,\, u^{*} \, \, $

\noindent  (1.6.5) $\qquad \displaystyle 
p_1 \,\, = \,\, p_2 \,\, = \,\, p^{*} \,. \, $

\smallskip \noindent   $\bullet \qquad \,\,\, $  
The relation (1.6.1) means that state   $\,W_1 \,$  is issued from state  $\,W_l \,$ 
through a ~1-rarefaction wave (relation (1.3.10)) or a ~1-shock wave (inequalities
(1.4.7)). We have in consequence~: 

\setbox21=\hbox {$\displaystyle  u_1 \,- \, u_l \,+\,  \psi( p_1 \,;\, \rho_l \,,\,
p_l \,;\, \gamma ) \,\,=\,\, 0 \,\,,\qquad p_1 \,< \, p_l \,   $}
\setbox22=\hbox {$\displaystyle  u_1 \,- \, u_l \,+\,  \varphi( p_1 \,;\, \rho_l \,,\,
p_l \,;\, \gamma ) \,\,=\,\, 0 \,\,,\qquad p_1 \,> \, p_l \,  . \,  $}
\setbox30= \vbox {\halign{#&# \cr \box21 \cr \box22    \cr   }}
\setbox31= \hbox{ $\vcenter {\box30} $}
\setbox44=\hbox{\noindent  (1.6.6) $\displaystyle  \qquad \left\{ \box31 \right. $} 

\noindent $ \box44 $

\noindent 
In a similar way, we remark that relation (1.6.3) expresses that state $\,W_2 \,$ is
the upstream state if a ~3-wave whose downstream state is exactly  the datum $\,W_r
.\,\,$ We can use a ~3-rarefaction wave between   $\,W_2 \,$   and  $\,W_r \,$  
(relation (1.3.17)) or a ~3-shock wave (relation (1.4.19)). It comes~: 

\setbox21=\hbox {$\displaystyle  u_2 \,- \, u_r \,-\,  \psi( p_2 \,;\, \rho_r \,,\,
p_r \,;\, \gamma ) \,\,=\,\, 0 \,\,,\qquad p_2 \,< \, p_r \,   $}
\setbox22=\hbox {$\displaystyle  u_2 \,- \, u_r \,-\,  \varphi( p_2 \,;\, \rho_r \,,\,
p_r \,;\, \gamma ) \,\,=\,\, 0 \,\,,\qquad p_2 \,> \, p_r \,  . \,  $}
\setbox30= \vbox {\halign{#&# \cr \box21 \cr \box22    \cr   }}
\setbox31= \hbox{ $\vcenter {\box30} $}
\setbox44=\hbox{\noindent  (1.6.7) $\displaystyle  \qquad \left\{ \box31 \right. $}

\noindent $ \box44 $

\noindent
It is sufficient to write the equation (1.6.4) that links velocities  $\,u_1 \,$  and 
$\,u_2 \,$  ({\it i.e.} $\, u_1 = u_2 = u^* $) under the condition (1.6.5) related to
the common pressures   $\,p_1 \,$  and   $\,p_2 \,\,$  ($p_1 = p_2 = p^*$) to set
completely the problem (see also Figure 1.8). 

%%%%%%%%%%%%%%%%%%%%%%%%%%%%%%%%%%%%%%%%%%%%%%%%%%%%%%%%%%%%%%%%%%%%%%%%%%%%%%%%%%%% 
% \bigskip  
% \vfill \eject 
\centerline { \epsfysize=5,5cm    \epsfbox  {fig.1.8.epsf} }
\smallskip  \smallskip

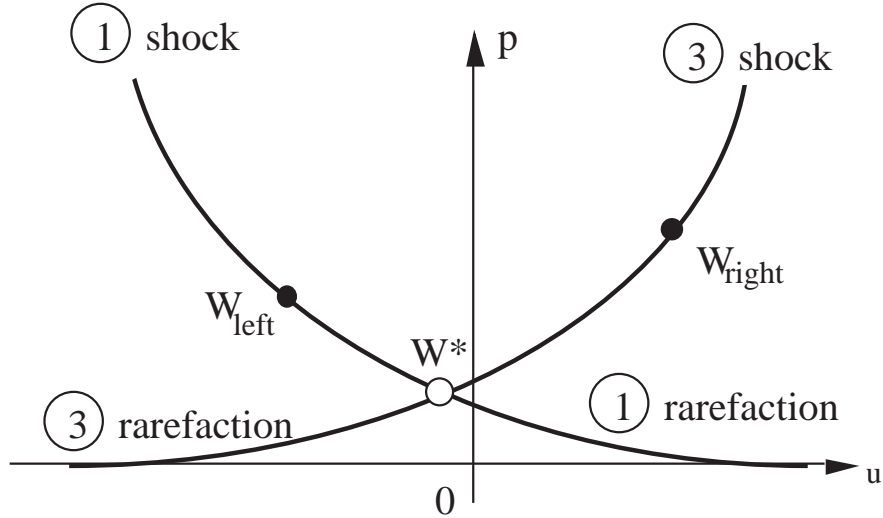
\noindent   {\bf Figure 1.8}	\quad  {\it Solution of the Riemann problem  
$\, \, R(W_{\rm left} \,,\,  W_{\rm right}) \,$   in the velocity-pressure plane. }
% \smallskip  \smallskip 
%%%%%%%%%%%%%%%%%%%%%%%%%%%%%%%%%%%%%%%%%%%%%%%%%%%%%%%%%%%%%%%%%%%%%%%%%%%%%%%%%%%% 

\smallskip \noindent   $\bullet \qquad \,\,\, $  
The numerical resolution of problem (1.6.4)-(1.6.7) can be done by Newton iterations
indexed by some integer $\, k \,$ and presented on Figure 1.9. Starting from a given
pressure  $\,p_k \,,\,$   we easily compute with relations (1.6.6) and (1.6.7)
velocities $\, u_{1,\,k} \,\,$ and  $\, u_{3,\,k} \, \,$ respectively associated with
the ~1-wave  and the ~3-wave. We evaluate the locus of intersection of the two
tangent lines issued from  the two corresponding velocities  in order to define a new 
value  $\,p_{k\!+\!1} \,$ for pressure at iteration   $\,\,k\!+\!1 .\,$ The
initialization of the algorithm can be obtained by the intersection of the two
rarefaction waves, {\it i.e.} by solving the following system with unknowns   $\, (
u_{O} ,\, p_{O} )
\,$~: 

\setbox21=\hbox {$\displaystyle  u_{O} \,- \, u_l \,+\,  \psi( p_{O} \,;\, \rho_l
\,,\, p_l \,;\, \gamma ) \,\,=\,\, 0    $}
\setbox22=\hbox {$\displaystyle   u_{O} \,- \, u_r \,+\,  \psi( p_{O} \,;\, \rho_r
\,,\, p_r \,;\, \gamma ) \,\,=\,\, 0   \,. \,  $}
\setbox30= \vbox {\halign{#&# \cr \box21 \cr \box22    \cr   }}
\setbox31= \hbox{ $\vcenter {\box30} $}
\setbox44=\hbox{\noindent  (1.6.8) $\displaystyle  \qquad \left\{ \box31 \right. $} 
\noindent $ \box44 $

%%%%%%%%%%%%%%%%%%%%%%%%%%%%%%%%%%%%%%%%%%%%%%%%%%%%%%%%%%%%%%%%%%%%%%%%%%%%%%%%%%%% 
% \bigskip  \vskip -,8 cm
\centerline { \epsfysize=6,5cm    \epsfbox  {fig.1.9.epsf} }
\smallskip  \smallskip

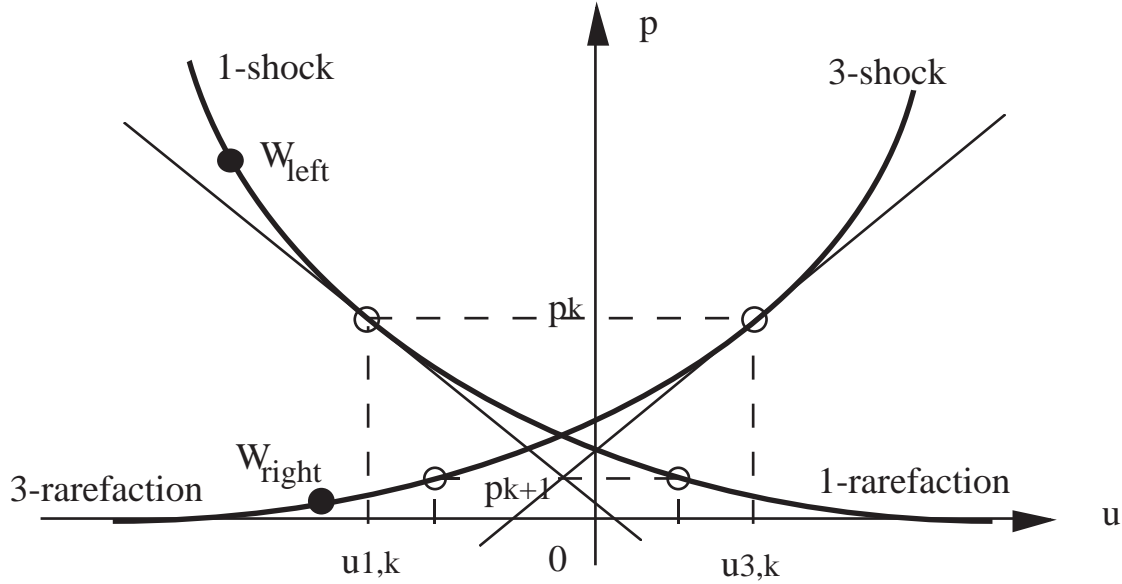
\noindent   {\bf Figure 1.9}	\quad  {\it Newton iterations for the resolution of 
the Riemann problem in the  plane of velocity and pressure. }
% \smallskip  \smallskip 
%%%%%%%%%%%%%%%%%%%%%%%%%%%%%%%%%%%%%%%%%%%%%%%%%%%%%%%%%%%%%%%%%%%%%%%%%%%%%%%%%%%%

\noindent
Note that in relations (1.6.8), the index ``O'' stands for ``Osher'' because solving
equations  (1.6.8) is essentially what is to be done  for the computation of the Osher
[Os81] flux decomposition (see {\it e.g.} [Du87]).  This system of equations, even if it is a
nonlinear one, can be solved exactly with some explicit algebra and pressure  $\,
p_{O} \,$ is finally evaluated thanks to the relation~: 

\noindent  (1.6.9) $\qquad \displaystyle
p_{O}^{{\scriptstyle \gamma \!-\! 1 }\over{\scriptstyle 2   \gamma}} \,\,= \,\,
{{{{(\gamma \!-\!1)}\over{2}}\,(u_l-u_r)  \,+\, c_l+c_r}\over
{c_l\,\bigl({{1}\over{p_l}}\bigr)^{{\scriptstyle \gamma \!-\! 1 }\over{\scriptstyle 2  
\gamma}} \,+\, c_r\,\bigl({{1}\over{p_r}}\bigr)^{{\scriptstyle \gamma \!-\! 1
}\over{\scriptstyle 2   \gamma}} }}  \,\,, \,$

\noindent
where sound celerities  $\,c_l\,$   and  $\,c_r\,$ are evaluated  with relation
(1.2.17). We remark that relation (1.6.9) defines effectively a {\bf positive}
pressure if the following  relation of {\bf non vacuum} appearance is satisfied~: 

\noindent  (1.6.10) $\qquad \displaystyle
u_r \,-\, u_l \,\,\leq\,\, {{2}\over{\gamma \!-\!1 }} \, (c_l + c_r) \,.\,$

%%%%%%%%%%%%%%%%%%%%%%%%%%%%%%%%%%%%%%%%%%%%%%%%%%%%%%%%%%%%%%%%%%%%%%%%%%%%%%%%%%%
\bigskip 
\centerline { \epsfysize=4.3cm    \epsfbox  {fig.1.10.epsf} }
\smallskip  \smallskip

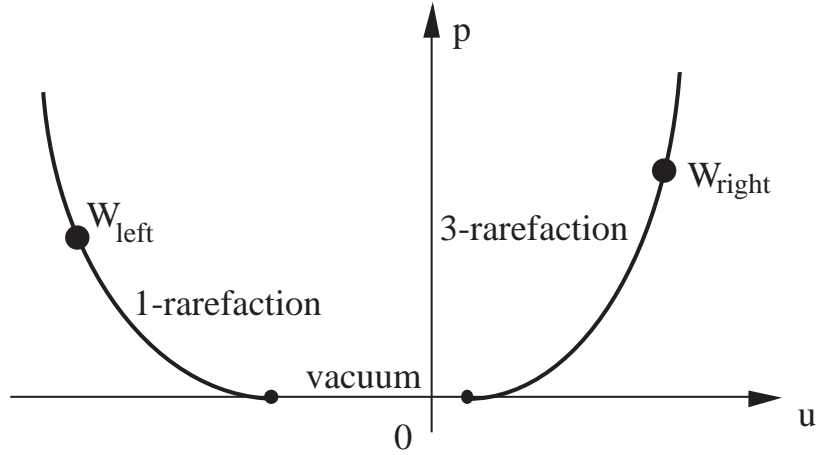
\noindent   {\bf Figure 1.10}	\quad  {\it Apparition of cavitation phenomenon 
in the velocity-pressure plane. }
\smallskip  \smallskip 
%%%%%%%%%%%%%%%%%%%%%%%%%%%%%%%%%%%%%%%%%%%%%%%%%%%%%%%%%%%%%%%%%%%%%%%%%%%%%%%%%%%

%%%%%%%%%%%%%%%%%%%%%%%%%%%%%%%%%%%%%%%%%%%%%%%%%%%%%%%%%%%%%%%%%%%%%%%%%%%%%%%%%%% 
\bigskip 
\centerline { \epsfysize=5,0cm    \epsfbox  {fig.1.11.epsf} }
\smallskip  \smallskip

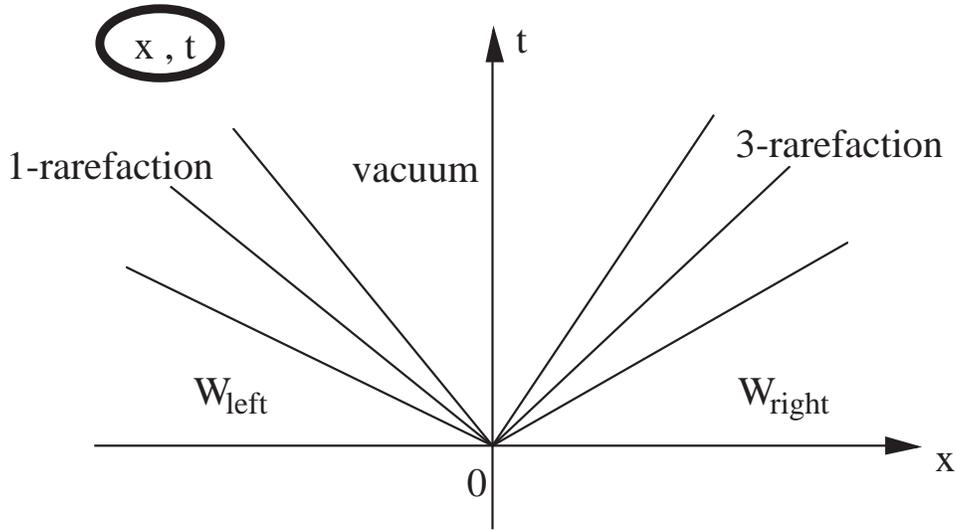
\noindent   {\bf Figure 1.11}	\quad  {\it Apparition of cavitation phenomenon. Solution of the 
Riemann  problem in the space-time plane. }
\smallskip  \smallskip 
%%%%%%%%%%%%%%%%%%%%%%%%%%%%%%%%%%%%%%%%%%%%%%%%%%%%%%%%%%%%%%%%%%%%%%%%%%%%%%%%%%%

\smallskip \noindent   $\bullet \qquad \,\,\, $  
If in contrary relation (1.6.10) does not hold,   the two curves associated with both
rarefaction waves do not intersect  in the velocity-pressure plane. A vacuum appears
(see Figures 1.10 and 1.11). The solution of the Riemann problem is no longer
mathematically well defined in the sense of Lax and contains a zone without any matter
defined by celerities $\, \xi \,$  such that 

\noindent  (1.6.11) $\qquad \displaystyle
u_l \,+\, {{2 \, c_l}\over{\gamma \!-\!1 }} \,\, \leq \,\, \xi \,\, \leq \,\,
u_r \,-\, {{2 \, c_r}\over{\gamma \!-\!1 }} \,. \,$

\noindent 
For such celerities, the pressure and the density are  null whereas velocity is not
defined. This cavitation process remains an exception but can be still completely
solved as we  have just seen.  

\smallskip \noindent   $\bullet \qquad \,\,\, $  
When relation (1.6.10) is satisfied, the Newton algorithm illustrated at Figure 1.9
is convergent towards a pair  $\, (u^*,\,p^*) \,$  composed by the common velocity  $\,
u^* \in \R \,$  and pressure  $\, p^* > 0 \,$ of the two intermediate states  $\,
W_1 \,$ and  $ \, W_2 .\,$ The calculus of the density of these two
intermediate states depends on the choice of the wave effectively used for the
resolution of the Riemann problem. If pressure  $\, p^* \,$ is less or equal to the
left pressure  $ \,p_l ,\,$  the ~1-wave is a rarefaction wave ;  then the entropy
remains constant and we have 

\noindent  (1.6.12) $\qquad \displaystyle
\rho_1 \,\,= \,\, \Bigl( {{p^*}\over{p_l}} \Bigr)^{\!{{\scriptstyle  1 }\over
{\scriptstyle  \gamma}}} \, \rho_l \,\,\,, \qquad \qquad  \qquad p^* \,<\, p_l \,. \,$

\noindent
On the opposite case,  we use a ~1-shock wave and taking into account the relation
(1.4.9), we get finally 

\noindent  (1.6.13) $\qquad \displaystyle
\rho_1 \,\,= \,\, {{p^* \,+\, \mu^2 \,p_l}\over{p_l \,+\, \mu^2 \,p^*}} \,\, \rho_l 
 \,\,\,, \qquad  \qquad p^* \,>\, p_l \,\,. \,$

\noindent 
The celerity $\, \sigma_1 \,$ of the shock wave can be explicited~: 

\noindent  (1.6.14) $\qquad \displaystyle
\sigma_1 \,\,= \,\, u_l \,-\, \sqrt{ {{p^* \,+\, \mu^2 \,p_l}\over{(1\!-\!\mu^2) \,
\rho_l}} }  \,\,\, \qquad p^* \,>\, p_l \,\,. \,$

\smallskip \noindent   $\bullet \qquad \,\,\, $  
For the ~3-wave, there is an analogous discussion that conducts finally to the
following  relations~:

\noindent  (1.6.15) $\qquad \displaystyle
\rho_3 \,\,= \,\, \Bigl( {{p^*}\over{p_r}} \Bigr)^{\!{{\scriptstyle  1 }\over
{\scriptstyle  \gamma}}} \, \rho_r\,\,\,, \qquad \qquad  \qquad p^* \,<\, p_r \,$

\noindent  (1.6.16) $\qquad \displaystyle
\rho_3 \,\,= \,\, {{p^* \,+\, \mu^2 \,p_r}\over{p_r \,+\, \mu^2 \,p^*}} \,\, \rho_r 
 \,\,\,, \qquad  \qquad p^* \,>\, p_r \,\, \,$
 
\noindent  (1.6.17) $\qquad \displaystyle
\sigma_3 \,\,= \,\, u_r \,+\, \sqrt{ {{p^* \,+\, \mu^2 \,p_r}\over{(1\!-\!\mu^2) \,
\rho_r}} }  \,\,\, \qquad p^* \,>\, p_r \,\,. \,$

\noindent
The qualitative comportment in space-time of the solution of a Riemann problem for
 gas dynamics is illustrated on Figure 1.12. 

%%%%%%%%%%%%%%%%%%%%%%%%%%%%%%%%%%%%%%%%%%%%%%%%%%%%%%%%%%%%%%%%%%%%%%%%%%%%
% \smallskip \smallskip \vskip 4,8cm  \smallskip  \qquad  \qquad   
% \special{illustration  fig.1.12.epsf scaled  600}  \smallskip  \smallskip 
\bigskip 
\centerline { \epsfysize=5,8cm    \epsfbox  {fig.1.12.epsf} }
\smallskip  \smallskip

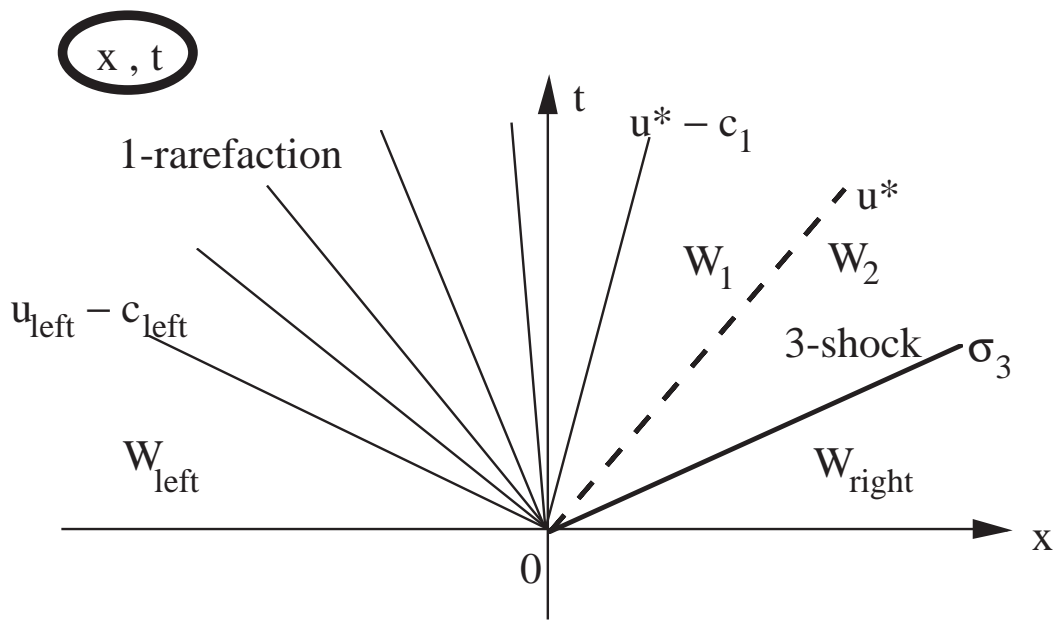
\noindent   {\bf Figure 1.12}	\quad  {\it Solution of the Riemann problem  
$\, \, R(W_{\rm left} \,,\,  W_{\rm right}) \,$  in the space-time plane.}
%%%%%%%%%%%%%%%%%%%%%%%%%%%%%%%%%%%%%%%%%%%%%%%%%%%%%%%%%%%%%%%%%%%%%%%%%%%%

% \vfill \eject  
%%%%%%%%%%%%%%%%%%%%%%%%%%%%%%%%%%%%%%%%%%%%%%%%%%%%%%%%%%%%%%%%%%%%%%%%%%%%%%%%%%%%%%%%%%
\bigskip    \bigskip    
\noindent  {\smcaps 2) $ \,\,\,\,\,\,\, $   Partial Riemann problem for hyperbolic systems.} 

\noindent   $\bullet \qquad \,\,\,\,\, $ 
In this section we generalize in two directions  the notion of Riemann problem
presented in the previous section for the particular case of gas dynamics. First we
consider a general nonlinear system of conservation laws and second we introduce the
notion of {\bf partial Riemann problem} between a state and a manifold.

\smallskip \noindent {\smcaps 2.1 } $ \,\, ${ \bf Simple waves for an hyperbolic
system  of conservation laws.}

\noindent   $\bullet \qquad \,\,\, $ 
We study a system of conservation laws in one space dimension. 
The unknown function $\,\, \, \R  \times  [0,\,+\infty[ \,  \, \ni (x,\,t)
\longmapsto  W(x,\,t) \in \Omega \,\, \,$ takes its values inside a convex open cone
$\, \Omega\,$ included in $\, \R^{m}\,$ ($m$ is a fixed positive integer and $\, m=3 \,
$ in the particular case of the Euler equations of gas dynamics studied previously)~:

\noindent  (2.1.1) $\qquad \displaystyle
\forall \, W \in \Omega\,, \quad \forall \, \lambda \, {\rm >}\,0\, , \quad \lambda 
\, W \in \Omega\, , \qquad \Omega \subset \R^{m} \,.\,  $

\noindent
The flux function $\,\, \Omega \ni W  \longmapsto F(W) \in \R^{m} \,\,$ is 
supposed to be sufficiently regular (of $\, {\cal C}^{2}\,$  class typically). The
conservation law takes the following classical form (see {\it {\it e.g.}} Godlewski-Raviart
[GR96])~: 

\noindent  (2.1.2) $\qquad \displaystyle
{{\partial }\over{\partial t}} W(x,\,t)\, + \, {{\partial}\over{\partial x}}\,
F(W(x,\,t)) \,\,=\,\,0\,. $

\noindent
We suppose that the system of conservation laws (2.1.2) is a strictly hyperbolic  
system, that is for each $\,W \in \Omega \,$, there exists $\,m\,$ real distinct 
eigenvalues $\, \lambda_{j}(W)\,$ satisfying conventionnaly  the ordering
condition 

\noindent  (2.1.3) {$\qquad \displaystyle \lambda_{1}(W) \,\, < \,\,
\lambda_{2}(W)  \,\, < \,\, \cdots \,\, < \,\, \lambda_{m}(W) \,, \qquad
W \in \Omega $}

\noindent 
and associated with eigenvectors $\,\, R_{j}(W) \, \in \R^{m}\,\,$:

\noindent  (2.1.4) {$\qquad \displaystyle 
{\rm d}F(W) \, {\scriptstyle \bullet} \,  R_{j}(W)
\,\,=\,\, \lambda_{j}(W)\,\,  R_{j}(W)  \,, \qquad W \in \Omega \,,\quad
j\,=\,1\,,\cdots \,, m\,. $}

\noindent 
Each vector $\, r \in \R^{m} \,$ can be decomposed in the basis of eigenvectors $\,  
R_{j}(W) \,$ and the coordinates in this basis define the family of
 left-eigenvectors $\,\, \bigl(l_{j}(W)\bigr) _{j=1,\cdots \, m} \, \,$ which is
exactly the  dual basis of system $\, \smash{\bigl(R_{j}(W)\bigr) _{j=1,\cdots \, m}}
\,$~:

\noindent  (2.1.5) {$\qquad \displaystyle 
r\,\, \equiv \,\, \sum_{k=1}^{k=m} \, \bigl( l_{j}(W) \,,\, r   \bigr) \,
R_{j}(W)\,  $}

\noindent 
with the classical property~:

\noindent  (2.1.6) {$\qquad \displaystyle 
 \bigl( l_{j}(W) \,,\,  R_{k}(W)   \bigr) \,\, = \,\, \delta_{jk} \,\, . \, $
 
\smallskip \noindent  {\bf Hypothesis 1. $\quad$  Genuinely nonlinear or linearly
degenerate fields. }

 \noindent 
We restrict ourselves to systems of conservation laws such that for each integer  $\,
j \,$ with $\,\, 1 \leq \, j \, \leq \, m , \, $ the $\,j$-th  field is
supposed to be 

\noindent ${\scriptstyle \bullet} \quad$  either genuinely nonlinear, that is 

\noindent  (2.1.7) $\quad \displaystyle 
{\rm d}\lambda_{j}(W) \, {\scriptstyle \bullet} \,  R_{j}(W) \, \equiv
1\,,\quad \forall \,  W \in \Omega \, \quad   (j$-th field  genuinely nonlinear)

\noindent ${\scriptstyle \bullet} \quad$  either linearly degenerate, {\it id est}

\noindent  (2.1.8) $\quad \displaystyle 
{\rm d}\lambda_{j}(W) \, {\scriptstyle \bullet} \,  R_{j}(W) \, \equiv
0\,,\quad  \forall \,  W \in \Omega \, \quad   (j$-th field  linearly degenerate).

\smallskip \noindent   $\bullet \qquad \,\,\, $ 
When the $\, j$-th field is genuinely nonlinear, it is possible to construct
the so-called $\,j-$wave issued from a particular state $\, W_0 \in \Omega.\,$ This
(genuinely) nonlinear wave can be considered from two points of view.  On the first
hand, it is a curve $\, \epsilon \longmapsto \smash{\chi\ib{j}(\epsilon \,;\, W_0)
}\,$ inside the space $\, \Omega \,$ of all the states and on the other hand for real
variable  $\, \epsilon \,$ fixed sufficiently small, it is possible to construct a
self-similar weak solution of   the conservation law (2.1.2) between state $\,W_0\,$
and state $\, \smash{\chi\ib{j}(\epsilon\,;\, W_0) }.\,$ 

\smallskip \noindent   $\bullet \qquad \,\,\, $ 
When $ \, \epsilon > 0 \,$ this particular wave is a $j-${\bf rarefaction} wave and is
defined as an integral curve of the vector field $\,\, \Omega \ni W \longmapsto 
R_{j}(W) \in \R^m \,\,$   in the space of states~:

\noindent  (2.1.9) $\qquad \displaystyle 
\,\, \, {{\partial}\over{\partial \epsilon}}  \bigl( \smash{\chi\ib{j}(\epsilon\,;\,
W_0) \bigr) \,\,= \,\,  R_{j} \bigl( \chi\ib{j}(\epsilon\,;\, W_0) \bigr)
\,\,, \quad \epsilon}  > 0
\,$

\noindent  (2.1.10) $\qquad \displaystyle 
 \smash{\chi\ib{j}(0,\, W_0) \,\quad \,= \,\, W_0 }\, . $

\noindent 
In space-time plane, this $j-$curve allows to construct a continuous rarefaction,
which is a particular solution of the conservation law (2.1.2)~: 

\noindent  (2.1.11) $\quad \displaystyle 
W(x,\,t) \,\,= \,\, W_0 \,\, \qquad  \qquad \qquad \qquad  \,\,\,\,\,  {\rm
if} \quad {{x}\over{t}} \,< \, \lambda_{j} (W_0) \, $

\noindent  (2.1.12) $\quad \displaystyle 
W(x,\,t) \,\,= \,\,  \chi\ib{j}\Bigl({{x}\over{t}}-\lambda_{j}(W_0) \,;\, W_0
\Bigr)\quad \, {\rm if} \,\,  \lambda_{j} (W_0)  \,\leq  \,{{x}\over{t}} \,\leq 
\,  \lambda_{j} \bigl(  \chi\ib{j}(\epsilon\,;\, W_0) \bigr) \, $

\noindent  (2.1.13) $\quad \displaystyle 
W(x,\,t) \,\,= \,\, \smash{\chi\ib{j}(\epsilon\,;\, W_0)  \,\, \qquad  \qquad \quad
\,\, \,\,  {\rm if} \quad  {{x}\over{t}} \,  > \,  \lambda_{j} \bigl( 
\chi\ib{j}(\epsilon\,;\, W_0) \bigr) } \,\,. $

\noindent
Moreover we have 

\noindent  (2.1.14) $\qquad \displaystyle 
 \lambda_{j} \bigl(W(x,\,t) \bigr) \,\, \equiv  \,\, {{x}\over{t}}  \qquad {\rm if}
\quad  \lambda_{j} (W_0)  \,\leq  \,{{x}\over{t}} \,\leq  \,  \lambda_{j}
\bigl(  \smash{\chi\ib{j}(\epsilon\,;\, W_0) \bigr) }\,. $

\smallskip \noindent   $\bullet \qquad \,\,\, $ 
When $ \, \epsilon < 0 \,$ the  $\, j$-th  nonlinear wave    is a  $j-${\bf
shock}  wave  of celerity $\,\, \sigma_{j}(\epsilon\,;\, W_0) \,$  satisfying
the entropy condition ({\it e.g.} in the sence of Lax [Lax73]). The states $\,W_0\,$ and $\,
\smash{\chi\ib{j}(\epsilon\,;\, W_0) }\,$ are linked with   celerity  $\,\,
\sigma_{j}(\epsilon\,;\, W_0) \,$ in state space according to the   Rankine-Hugoniot
jump conditions~:  

\noindent  (2.1.15) $\quad \displaystyle 
F\bigl( \chi\ib{j}(\epsilon\,;\, W_0) \bigr) \,-\, F(W_0) \,\,\equiv \,\,
\sigma_{j}(\epsilon\,;\, W_0) \,\bigl( \chi\ib{j}(\epsilon\,;\, W_0) \,-\,
W_0 \bigr) \,\,, \quad \epsilon < 0 \,. $

\noindent 
In space-time space, this discontinuous  self-similar  $j-$shock wave   $\,\, \, \R  
\times  [0,\,+\infty[ \,  \, $ $ \, \ni (x,\,t) \longmapsto  W(x,\,t) \in \Omega \,\,
\,$ of strength $\, \abs{ \epsilon} \, $ is characterized by the following two 
conditions 

\noindent  (2.1.16) $\qquad \displaystyle 
W(x,\,t)\,\,= \,\, W_0 \,\, \qquad  \qquad \qquad {\rm
if} \,\, {{x}\over{t}} \,< \, \sigma_{j}(\epsilon\,;\, W_0) \,$ 

\smallskip \noindent  (2.1.17) $\qquad \displaystyle 
W(x,\,t) \,\,= \,\, \smash{\chi\ib{j}(\epsilon\,;\, W_0)  \,\, \qquad \quad  {\rm
if} \,\,  {{x}\over{t}} }\,  > \,  \sigma_{j}(\epsilon\,;\, W_0)  \,. \, $ 

\smallskip \noindent   $\bullet \qquad \,\,\, $ 
When the  $\, j$-th field is linearly degenerate, the construction of the 
$\,j-$wave issued from the particular state $\, W_0 \in \Omega\,$ is still possible.
Due to the condition (2.1.8), the $j^{\rm o}\, $ eigenvalue $\, \lambda_{j}(W) \,$
is constant along the integral curve defined in relations (2.1.9) (2.1.10) and we have
also all along this curve the following jump relation~: 

\noindent  (2.1.18) $\qquad \displaystyle 
F\bigl(\chi\ib{j}(\epsilon\,;\, W_0) \bigr) \,-\, F(W_0) \,\,\equiv \,\,
\lambda_{j}(\epsilon\,;\, W_0) \,\bigl( \chi\ib{j}(\epsilon\,;\, W_0) \,-\,
W_0 \bigr) \,\,,\,\,\forall \, \epsilon \,. \,$

\noindent 
In space-time plane, a self-similar  $j-${\bf contact discontinuity} can be constructed
between states $\, W_0\,$ and $\,  \chi\ib{j}(\epsilon\,;\, W_0) \,$ and we have~: 

\noindent  (2.1.19) $\qquad \displaystyle 
W(x,\,t) \,\,= \,\, W_0 \,\, \qquad  \qquad \qquad {\rm
if} \quad {{x}\over{t}} \,< \, \lambda_{j}(\epsilon\,;\, W_0) \,$ 

\noindent  (2.1.20) $\qquad \displaystyle 
W(x,\,t) \,\,= \,\, \chi\ib{j}(\epsilon\,;\, W_0)  \,\, \quad  \qquad {\rm
if} \quad {{x}\over{t}} \,  > \,  \lambda_{j}(\epsilon\,;\, W_0)  \,. \, $ 

\noindent 
All this material is summarized in the following result (see {\it e.g.} [GR96]). 

\smallskip 
% \vfill \eject 
\noindent  {\bf Proposition 1. $\quad$  }

 \noindent 
We suppose that the hyperbolic system of conservation laws (2.1.2) satisfies
Hypothesis 1. Then for each state $\, W_0 \in \Omega $ and for each integer $\,j
\,\,(0\leq j \leq m) ,\,$ there exists some vicinity  $\, \Theta_j  \, $ of $\,0\,$ in
$\,\R\,$ and  there exists a $j-$wave $\,\, \smash{ \chi\ib{j} \,: \,  \Theta_j \times
\Omega \ni (\epsilon\,, } $~\br
 $ \smash{ \,W_0) \longmapsto  \chi\ib{j}(\epsilon\,;\, \,W_0)  }   \in \Omega  \,\,$ 
which is a regular  (of $\, {\cal C}^{2}\,$  class) curve. This corresponds to 
a $j-$rarefaction when the $\, j^{\rm o}\, $ field is genuinely nonlinear and $\,
\epsilon > 0  ,\,$  to an admissible discontinuity (weak shock  satisfying an entropy
condition) of conservation law (2.1.2)   when the $\, j^{\rm o}\,
$ field is genuinely nonlinear and $\,  \epsilon < 0 \,$ or to a  $j-$contact
discontinuity  when  the $\, j^{\rm o}\, $ field is  linearly degenerate. 
Moreover, the  $\, j^{\rm o}\, $ wave $\,  \chi\ib{j}({\scriptstyle \bullet} \,;\, W_0) \,$ 
is continuously derivable at starting point $\,W_0\,$ and we have~: 

\noindent  (2.1.21) $\qquad \displaystyle 
 \chi\ib{j}(\epsilon\,;\, W_0) \,\,= \,\, W_0 \,+\, \epsilon \,R_{j}(W_0) \,+ \,
O\bigl(\epsilon^2 \bigr) \,, \quad \epsilon \in  \Theta_j  \,. \, $

\smallskip \noindent   $\bullet \qquad \,\,\, $ 
We remark that in space-time, a $j$-wave defines an entropy solution of hyperbolic
system (2.1.1). 

\bigskip  \noindent {\smcaps 2.2 } $ \, ${ \bf Classical Riemann problem between two
states.}

\noindent   $\bullet \qquad \,\,\, $ 
We consider the solution of the particular Cauchy problem associated with the
conservation law (2.1.1) and the particular initial datum

\setbox21=\hbox {$\displaystyle  W_l \,\,, \qquad x \,< \, 0 \,  $}
\setbox22=\hbox {$\displaystyle  W_{\! r} \,\,, \qquad x \,> \, 0 \,.\,   $}
\setbox30= \vbox {\halign{#&# \cr \box21 \cr \box22    \cr   }}
\setbox31= \hbox{ $\vcenter {\box30} $}
\setbox44=\hbox{\noindent  (2.2.1) $\displaystyle  \qquad   W(x \,,\,0) \,\,= \,\,
\left\{ \box31 \right. $}  

\noindent $ \box44 $

\noindent
composed by two constant states on each side of a discontinuity located at $\, x =
0. \,$  In what follows, this problem is refered as the {\bf classical Riemann
problem} and is denoted by $\, \, R(W_l \,,\,  W_{\! r}) .\,$ We have the
following Theorem due to Lax (see {\it e.g.}  [Lax~73]).

\bigskip 
%%%%%%%%%%%%%%%%%%%%%%%%%%%%%%%%%%%%%%%%%%%%%%%%%%%%%%%%%%%%%%%%%%%%%%%%%%%% 
\bigskip 
\centerline { \epsfysize=4,5cm    \epsfbox  {fig.2.1.epsf} }
\smallskip  \smallskip
 
\centerline { \epsfysize=4,0cm    \epsfbox  {fig.2.1.bis.epsf} }
\smallskip  \smallskip

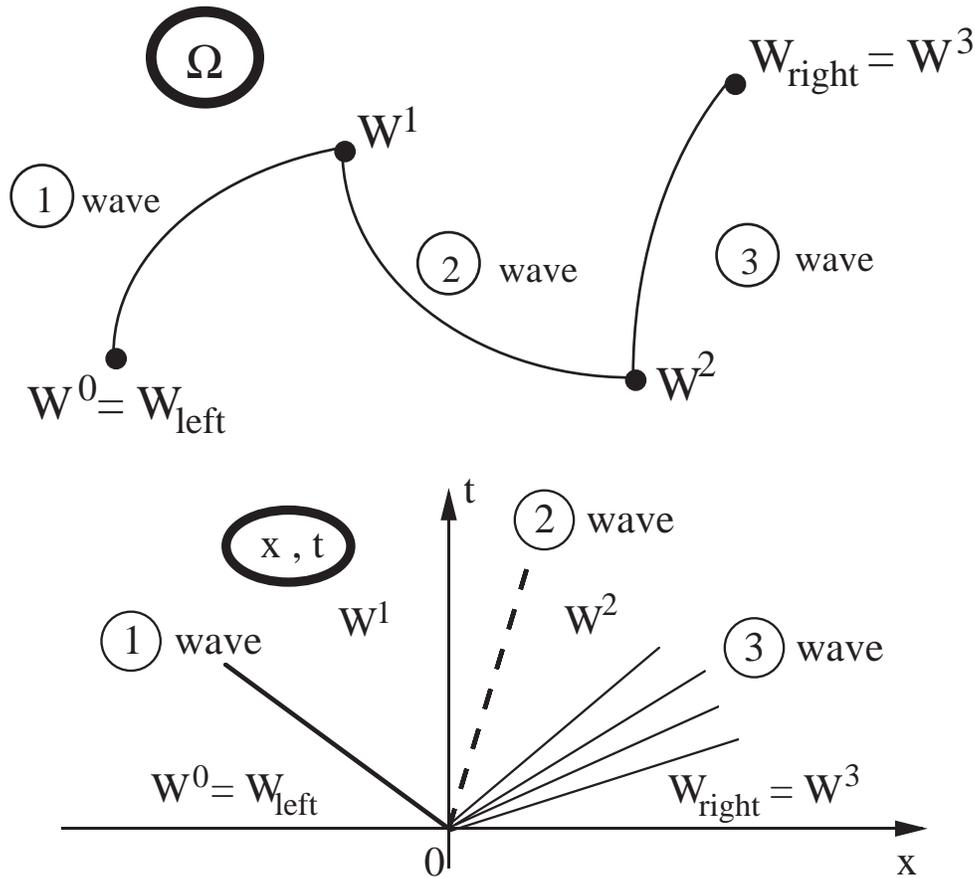
\noindent   {\bf Figure 2.1}	\quad   {\it Solution of the Riemann problem 
in the state space  $ \, \Omega \,$   and in the space-time plane.}
\smallskip  \smallskip 
%%%%%%%%%%%%%%%%%%%%%%%%%%%%%%%%%%%%%%%%%%%%%%%%%%%%%%%%%%%%%%%%%%%%%%%%%%%%

% \smallskip 
 \vfill \eject 
\noindent  {\bf Theorem 1. $\quad$  Solution of the classical Riemann
problem.}

\noindent 
We suppose that the hyperbolic system of conservation laws (2.1.2) satisfies
Hypothesis 1~: for $ \,\, j \in \{1,\, \cdots ,\, m\} ,\,$ the $\, j^{\rm o}\, $ field
is either genuinely nonlinear or  linearly degenerate. Then for each $\,\, W_l \in
\Omega \,\,$ there exists a vicinity  $\, {\cal Y} \, $ of state $\,\, W_l \,\,$ 
$\, ({\cal Y} \in  {\cal V}(W_l) ,\,$  set of all  vicinities of the
particular state $\, W_l \,$) such that for each state $\, W_r \, $ lying in $\, 
{\cal Y} \, $ $\,\, ( W_r \in {\cal Y} ) ,\, \,$  the Riemann problem $\, R(W_l \,,\, 
W_{\! r}) \,$ has a unique  entropy solution composed by at most $\, (m+1) \,$ states
separated by (at most) $\, m \,$ elementary  waves. 

\smallskip \noindent   $\bullet \qquad \,\,\, $ 
Figure 2.1 shows what has to be done for the resolution of the classical Riemann
problem in the particular case $\,\,  m = 3 \, \,$ that corresponds to the Euler
equations of gas dynamics in one space dimension. The two intermediate states $\, W^1
\,$ and $\,W^2 \,$ are linked with data $\, W_l\,$ and $\, W_r \,$ acording to the wave
relations~: 

\noindent  (2.2.2) $\qquad \displaystyle 
W^1 \,\,= \,\, \chi\ib{1} (\epsilon\ib{1} \,;\, W_l) \,$

\noindent  (2.2.3) $\qquad \displaystyle 
W^2 \,\,= \,\, \chi\ib{2} (\epsilon\ib{2} \,;\, W^1) \,$

\noindent  (2.2.4) $\qquad \displaystyle 
W_r \,\,= \,\, \chi\ib{3} (\epsilon\ib{3} \,;\, W^2) \, \, $

\noindent 
and  by elimination of the two intermediate states $\, W^1 \,$ and $\, W^2 \,$ we
obtain

\noindent  (2.2.5) $\qquad \displaystyle 
\chi\ib{3} \Bigl( \epsilon\ib{3} \,;\, \chi\ib{2} \bigl( \epsilon\ib{2} \,;\,
\chi\ib{1} ( \epsilon\ib{1} \,;\, W_l ) \bigr) \Bigr) \,\,= \,\, W_r\,$

\noindent 
which is a set of three equations     with three scalar unknowns $\, \epsilon\ib{1} ,\,
\epsilon\ib{2} \,$ and $\, \epsilon\ib{3} .\,$ When parameters $\,  \epsilon\ib{j} \,$
are determined (in practice for gas dynamics  with the method presented in the
previous section) each wave is acting in space-time space as presented on Figure 2.1. 

\smallskip \noindent   $\bullet \qquad \,\,\, $ 
The proof of Theorem 1  in the general case consists in studying the chaining of $\, m
\,$ elementary waves, {\it i.e.} the mapping $\,\, \chi \,\,$ defined by the relations 

\setbox21=\hbox {$\displaystyle  \R^m \,\, \supset \,\,  \Theta_1 \times \, \cdots \,
\, \times \Theta_m \,\, \ni \, (\epsilon\ib{1} \,,\, \cdots \,,\, \epsilon\ib{m})
\, \equiv\, \epsilon \, \longmapsto \,  \chi(\epsilon) \, \in \Omega  $}
\setbox22=\hbox {$\displaystyle   \chi(\epsilon) \equiv\,  \chi\ib{m} \bigl(
\epsilon\ib{m} \,;\, \chi\ib{m-1} \bigl( \epsilon\ib{m-1} \,;\, \cdots \,;\, \chi
_{1}(\epsilon\ib{1} \,;\, W_l) \cdots \bigr) \bigr) \,.\,    $}
\setbox30= \vbox {\halign{#&# \cr \box21 \cr \box22    \cr   }}
\setbox31= \hbox{ $\vcenter {\box30} $}
\setbox44=\hbox{\noindent  (2.2.6) $\displaystyle  \qquad  
\left\{ \box31 \right. $}  

\noindent $ \box44 $

\noindent 
Two local properties have to be derived concerning on one hand the mapping  $\, \chi
\,$ itself at the origin~: 

\noindent  (2.2.7) $\qquad \displaystyle 
\chi(0) \,\,= \,\, W_l  \, $

\noindent 
and on the other hand the tangent vector field $\, {\rm d}\chi \,$ considered at the
same point has a very simple expression~: 

\noindent  (2.2.8) $\qquad \displaystyle 
{\rm d}\chi(0) \, {\scriptstyle \bullet} \, \eta \,\,= \,\, \sum_{j=1}^{m} \,
\eta_{j} \,  R_{j}(W_l)  \,\,, \qquad \eta \,=\, \bigl (\eta_{1},\, \cdots \,,
\eta_{m} \bigr) \in \R^m \,.\,  \, $

\noindent
Then the local inversion theorem proves that the equation 

\noindent  (2.2.9) $\qquad \displaystyle 
\chi(\epsilon) \,\,= \,\, W_r \,$ 

\noindent 
has a unique solution.   For the details of this proof, we refer to Godlewski and
Raviart [GR96]. $ \hfill \square \kern0.1mm $

\smallskip \noindent   $\bullet \qquad \,\,\, $ 
In space-time domain, the entropic selfsimilar solution $\,\, \R \times $ $ ]0,\,
+\infty [  \, \ni $ $(x,\,t)  \mapsto U({{x}\over{t}} \,; \, W_l ,\, W_r ) \,\,$ of
the Riemann problem $\,\,  R( W_l ,\, W_r ) \,\,$ is constructed as follows. Let 

\setbox21=\hbox {$\displaystyle  
W^0 \,=\, W_l \,, \cdots ,\,\, W^j= \chi\ib{j} ( \epsilon\ib{j} \,;\, W^{j-1}
) \,, \cdots ,\, $}
\setbox22=\hbox {$\displaystyle \qquad  \qquad \qquad \qquad \qquad \qquad 
\cdots ,\, W^m \,=\, \chi\ib{m} ( \epsilon\ib{m} \,;\, W^{m-1} ) \,=\, W_r \, $}
\setbox30= \vbox {\halign{#&# \cr \box21 \cr \box22    \cr   }}
\setbox31= \hbox{ $\vcenter {\box30} $}
\setbox44=\hbox{\noindent  (2.2.10) $\displaystyle  \qquad  
\left\{ \box31 \right. $}  

\noindent $ \box44 $

\noindent 
be the intermediate states, $\,\, \mu_{j}^-(W_l ,\, W_r ) \,\,$ the smallest wave
celerity of the $j^o$ wave and $\,\, \mu_{j}^+(W_l ,\, W_r ) \,\,$ the
corresponding maximal wave celerity. If the  $j^o$ wave is a rarefaction, we have
from (2.1.11) and (2.1.13)~: 

\noindent  (2.2.11) $\quad \displaystyle 
\mu_{j}^-(W_l ,\, W_r ) \,\,= \,\, \lambda_{j}(W^{j-1}) \,,\quad 
\mu_{j}^+(W_l ,\, W_r ) \,\,= \,\, \lambda_{j}(W^{j}) \,$ 

\noindent 
whereas in case of a $j$-shock wave or $j$-contact discontinuity, we have due to
(2.1.16) and (2.1.17)~: 

\noindent  (2.2.12) $\quad \displaystyle 
\mu_{j}^-(W_l ,\, W_r ) \,\,= \,\,\mu_{j}^+(W_l ,\, W_r ) \,\,= \,\,
\sigma_{j}(\epsilon\ib{j} \,;\,  W^{j-1}) \, .\,$ 

\noindent 
Then the selfsimilar solution $\,\, U(\xi \,; \, W_l ,\, W_r ) \,\,$ satisfies 

\setbox21=\hbox {$\displaystyle  W^0 \,=\, W_l  \,,\,   $}
\setbox31=\hbox {$\displaystyle  \!\!\!\!\!\!\!\!\!\!\!\!\!\!
\xi < \, \mu_{1}^-(W_l ,\, W_r ) $}
\setbox22=\hbox {$\displaystyle  \vdots $}
\setbox23=\hbox {$\displaystyle \chi\ib{j}\Bigl( {{x}\over{t}}-\mu_{j}^-(W_l ,\,
W_r ) \Bigr) ,\,\,    $}
\setbox33=\hbox {$\displaystyle  \!\!\!\!\!\!\!\!\!\!\!\!\!\!
\mu_{j}^-(W_l ,\, W_r ) < \xi < \mu_{j}^+(W_l ,\, W_r )  $}
\setbox24=\hbox {$\displaystyle  W^{j}   \,,\,   $}
\setbox34=\hbox {$\displaystyle  \!\!\!\!\!\!\!\!\!\!\!\!\!\!
\mu_{j}^+(W_l ,\, W_r ) < \xi < \mu_{j+1}^-(W_l ,\, W_r )  $}
\setbox25=\hbox {$\displaystyle  \vdots $}
\setbox26=\hbox {$\displaystyle  W^m \,=\, W_r  \,,\,   $}
\setbox36=\hbox {$\displaystyle  \!\!\!\!\!\!\!\!\!\!\!\!\!\!
\xi > \, \mu_{m}^+(W_l ,\, W_r ) \,. \, $}
\setbox40= \vbox {\halign{#&# \cr \box21 &   \box31  \cr \box22 \cr  \box23 \cr  
& \box33  \cr  
\box24 &   \box34  \cr \box25 \cr  \box26 &  \box36 \cr   }}
\setbox41= \hbox{ $\vcenter {\box40} $}
\setbox44=\hbox{\noindent  (2.2.13) $\displaystyle  \,    U(\xi \,; \, W_l ,\, W_r
) \,= \! \left\{ \box41 \right. $}  

\noindent $ \box44 $

\bigskip  
% \vfill \eject 
\noindent {\smcaps 2.3 } $ \, ${ \bf Boundary  manifold.}

\noindent   $\bullet \qquad \,\,\, $ 
The Riemann problem is a usefull tool to prescribe weakly a boundary condition for an
hyperbolic system of conservation laws when a right state $\, W_r \,$ is
supposed to be given (see Section 3.3). For physically relevant conditions in gas
dynamics, the data at the boundary are of the type ``the pressure is known'' or
``total pressure and total temperature are given'', and in consequence  do not define
explicitely a single state $\,\, W_r \,. \,$  Nevertheless, a {\bf manifold} of states
is associated with these sets of incomplete {\bf boundary} data, as we will develop in
Section 3.4. This physical situation motivates the following definition of a boundary
manifold. 

% \bigskip \smallskip \smallskip \vskip 7,5cm  \smallskip  \qquad  \qquad   
% \special{illustration  fig.2.2.epsf scaled  600}  \smallskip  \smallskip 
\bigskip 
\centerline { \epsfysize=8.5cm    \epsfbox  {fig.2.2.epsf} }
\smallskip  \smallskip

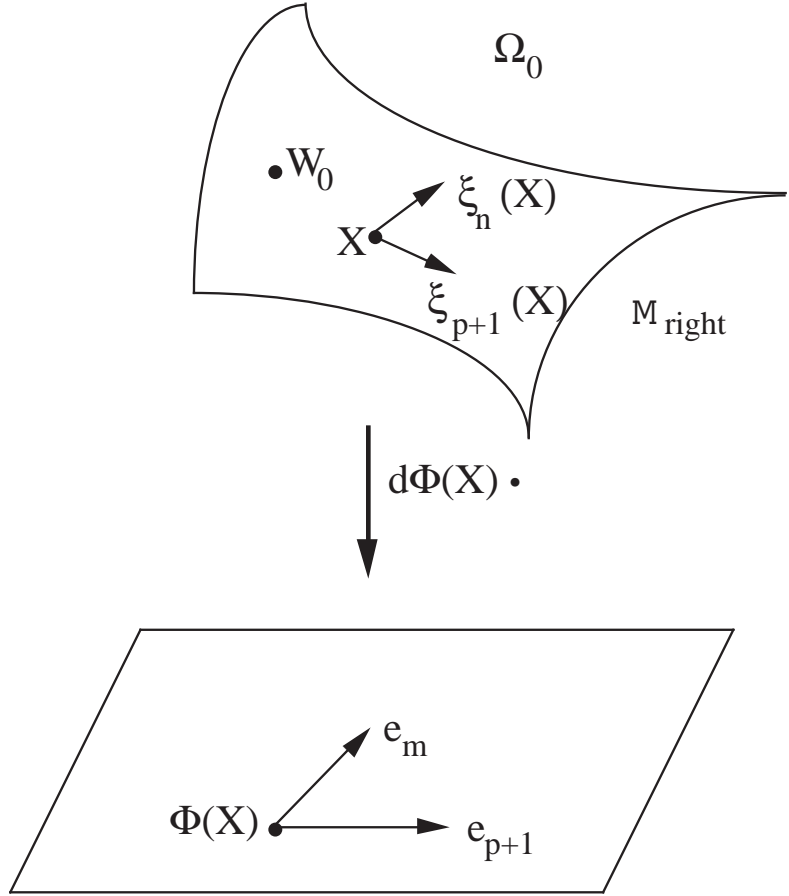
\noindent   {\bf Figure 2.2}	\quad  {\it Tangent vectors at the particular state
$\, X \, $  of  the boundary manifold $\, {\cal M}_{\rm right}  .\,$ }
\smallskip \smallskip

\smallskip \noindent   $\bullet \qquad \,\,\, $ 
We restrict ourselves to a local manifold $\,\, {\cal M}_r  \,\,$ around a given
state $\, W_0  .\,$ Consider a state $\, W_0 \in \Omega  ,\,$  a vicinity $\,
\Omega_0 \,$ of this state, {\it i.e.} 

\noindent  (2.3.1) $\qquad \displaystyle 
W_0 \in  \Omega \,\,,\quad \Omega_0 \in {\cal V}(W_0) 
\,\,,\quad \Omega_0 \subset \, \Omega \subset \, \R^m \,$ 

\noindent 
and an invertible regular local chart $\, \Phi \,$ defined on the vicinity $\, \Omega_0
\,$ and taking its values in some vicinity  $\, \Theta \subset \, \R^m \,$ of $\,0\,$
in $\, \R^m \,$~:

\noindent  (2.3.2) $\qquad \displaystyle 
\Omega_0 \ni  X \longmapsto \Phi(X) \in \Theta \, \subset \, \R^m \,\,, \quad \Theta
\in {\cal V}(0) \,.  \,$ 

\noindent 
Mapping  $\, \Phi \,$ is  one to one,  of $\, {\cal C}^1 \,$ class ($ \Phi \in 
{\cal C}^1 (\Omega_0, \, \Theta)) \,$  and its inverse mapping  $\, \Phi^{-1} \,$ is
also of $\, {\cal C}^1 \,$ class  ($ \Phi^{-1} \in  {\cal C}^1 (\Theta , \,
\Omega_0)) .\,$ The boundary  manifold $\, {\cal M}_{right}  \equiv {\cal M}_r  \,$  
is here defined {\bf  locally}  as the set of states $\, X \,$ in the
vicinity of  $\, \Omega_0 \,\,(X \in \Omega_0 ) \, $ satisfying the equations  

\noindent  (2.3.3) $\qquad \displaystyle 
\Phi_{1}(X) \,=\, \cdots \, = \, \Phi_{p}(X) \,=\, 0 \,,\quad X  \in \Omega_0
\,.\,$ 

\noindent
The index $\, p \,$ is a fixed integer such that  $\, 0 \leq p \leq m ~\,$ and is the 
 co-dimension of the boundary manifold $\,\, {\cal M}_{r} . \,\,$ Recall that for
each state $\, X \,$ (or point $\, X \in {\cal M}_{r}) ,\,$  the vector $\, \Phi(X)
\,$  in $\, \R^m \,$ has its  $\,p\,$ first coordinates   equal to zero~: 

\noindent  (2.3.4) $\qquad \displaystyle 
\Phi({\cal M}_r) \,\, \subset \,\, \bigl\{ (0,\, \cdots \,,\,0,\, y ) , \,\,\, y \in
\R^{m-p} \, \bigr\} \,. \,$ 

\smallskip \noindent   $\bullet \qquad \,\,\, $ 
A system of tangent vector fields $\,\, \bigl( \xi_{p+1}(X) \,,\, \cdots \,,\,
\xi_{m}(X) \bigr) \,$ at point $\, X \in {\cal M}_r \,\,$ is obtained by lifting
the tangent mapping $\, {\rm d}\Phi(X) \,$ at point $\, X \in  \Omega_0 .\,$ Let
$\,\, e_{k} \,\,$ be the $\,k^{\rm o}\,$ vector of the canonical basis in linear
space $\, \R^m .\,$ We have by definition~: 

\noindent  (2.3.5) $\qquad \displaystyle 
{\rm d}\Phi(X) \, {\scriptstyle \bullet} \, \xi_{k}(X) \,\,= \,\, e_{k} \,\,,
\quad X \in {\cal M}_r \,\,,\quad k \geq p+1 \,$ 

\noindent 
as illustrated on Figure 2.2. We make a new hypothesis.

\smallskip \noindent  {\bf Hypothesis 2. $\quad$ Transversality. }

 \noindent 
We suppose that the following family $\, \Sigma(W_0) \,$   of $\,m \,$ vectors~: 

\noindent  (2.3.6) $\qquad \displaystyle 
\Sigma(W_0) \,\,\equiv  \,\, \bigl( R_{1}(W_0) \,,\, \cdots \,,\, R_{p}(W_0)
\,,\, \xi_{p+1}(W_0) \,,\, \cdots \,,\, \xi_{m}(W_0) \bigr) \,$

\noindent
is a basis of linear space $\, \R^m \,.\,$ 

\bigskip  \noindent {\smcaps 2.4 } $ \, $ {\bf Partial Riemann problem between a 
 state and a  manifold. }

\noindent   $\bullet \qquad \,\,\, $ 
The partial Riemann problem $\,\,P(W_l,\, {\cal M}_r) \,\,$  between a state $\, W_l
\,$ and a manifold $\, {\cal M}_r \,$  is by definition the Cauchy problem for the
hyperbolic system of conservation laws (2.1.2) associated with the following
constraints relative to the  initial condition~:

\setbox21=\hbox {$\displaystyle  = \,\, W_l \,\,, \quad \qquad x \,< \, 0 \,  $}
\setbox22=\hbox {$\displaystyle  \in  {\cal M}_r \,\,,  \quad \qquad x \,> \, 0 \,.\,  
$}
\setbox30= \vbox {\halign{#&# \cr \box21 \cr \box22    \cr   }}
\setbox31= \hbox{ $\vcenter {\box30} $}
\setbox44=\hbox{\noindent  (2.4.1) $\displaystyle  \qquad   W(x \,,\,0) \, \,\,
\left\{ \box31 \right. $}  

\noindent $ \box44 $

\noindent
The above definition of the partial Riemann problem $\,\,P(W_l,\, {\cal M}_r) \,\,$ 
has been first proposed in a particular case in [Du87]. We have also used it in
[Du88], [DLF89],  [DLL91] and [CDV92].  We first remark that if $\, p = m \,$ and $\,\,
{\cal M}_r \,=\, \{W_r\} \,,\,$ then the partial Riemann problem $\, \, P(W_l,\,{\cal
M}_r) \,$  reduces to the classical Riemann problem $\, \, R(W_l,\,W_r) .\,$  The
following result has been first presented in [Du98]. 

\smallskip  \noindent  {\bf Theorem 2.}
 {\bf Existence of a solution for the partial Riemann problem.}

%  $\hfill$   {\bf Existence of a solution for the partial Riemann problem.}

 \noindent 
We suppose that the hyperbolic system of conservation laws (2.1.2) satisfies
Hypothesis 1 and that the boundary manifold $\, {\cal M}_r \,$ is  defined as above
with a  given state $\,\, W_0 \in \Omega \,$ that satisfies Hypothesis 2,  a vicinity
$\, \Omega_0 \,$ of $\, W_0 \,\,\, (\Omega_0 \in {\cal V}(W_0)) \,$  and a  local chart
$ \Omega_0 \ni  X \longmapsto \Phi(X) \in \Theta \, \subset \, \R^m \,, $ $\,\,
\Theta \in {\cal V}(0) \,\,$ and satisfying the relations  (2.3.3)-(2.3.4).  Then there
exists a vicinity $\, \Omega_1 \, \subset \Omega_0 \, $ of state $\, W_0 \,$ such that
for every state $\, W_l \in \Omega_1 ,\,$ the partial Riemann problem $\, P(W_l,\,
{\cal M}_r) \,$ defined by (2.1.2) and (2.4.1)  has an entropy  solution in space $\,
\Omega \,$  composed by at most $\, (p+1)\,$ states 

\setbox21=\hbox {$\displaystyle  W_l \,=\, W^0 \,,\,\, W^1 \,=\, 
\chi\ib{1}(\epsilon\ib{1},\, W^0 ) \,,\,\, \cdots \,, \,\, W^{p-1} \,=\,
\chi\ib{p-1}(\epsilon\ib{p-1},\, W^{p-2}) ,  $}
\setbox22=\hbox {$\displaystyle  W^p \,\,= \,\, \chi\ib{p}(\epsilon\ib{p},\,
W^{p-1}) \,, \qquad  W^p \in  {\cal M}_r , $}
\setbox30= \vbox {\halign{#&# \cr \box21 \cr \box22    \cr   }}
\setbox31= \hbox{ $\vcenter {\box30} $}
\setbox44=\hbox{\noindent  (2.4.2) $\displaystyle  \, 
\left\{ \box31 \right. $}  

\noindent $ \box44 $

\noindent
separated by (at most) $\, p \,$ simple waves. In space-time, the solution   $\, \R
\times ]0,\, +\infty [  $ $\, \ni $~$(x,\,t)  \longmapsto U({{x}\over{t}} \,; \, W_l
,\, {\cal M}_r ) \,\,$ is obtained by superposition of the $p$ simple waves considered
previously. With the notations proposed in (2.2.11) and (2.2.12), we have

\setbox21=\hbox {$\displaystyle  W^0 \,=\, W_l  \,,\,   $}
\setbox31=\hbox {$\displaystyle  \xi < \, \mu_{1}^- $}
\setbox22=\hbox {$\displaystyle  \vdots $}
\setbox23=\hbox {$\displaystyle \chi\ib{j}\Bigl( {{x}\over{t}}-\mu_{j}^- \Bigr)
,\,\,    $}
\setbox33=\hbox {$\displaystyle  \mu_{j}^- < \xi < \mu_{j}^+   $}
\setbox24=\hbox {$\displaystyle  W^{j}   \,,\,   $}
\setbox34=\hbox {$\displaystyle  \mu_{j}^+ < \xi < \mu_{j+1}^-  $}
\setbox25=\hbox {$\displaystyle  \vdots $}
\setbox26=\hbox {$\displaystyle  W^p \,\in\, {\cal M}_r  \,\,, \qquad    $}
\setbox36=\hbox {$\displaystyle  \xi > \, \mu_{p}^+ \,. \, $}
\setbox40= \vbox {\halign{#&# \cr \box21 &   \box31  \cr \box22 \cr  \box23 &  
\box33  \cr  \box24 &   \box34  \cr \box25 \cr  \box26 &  \box36 \cr   }}
\setbox41= \hbox{ $\vcenter {\box40} $}
\setbox44=\hbox{\noindent  (2.4.3) $\displaystyle  \qquad    U(\xi \,; \, W_l ,\, W_r
) \,\,=\,\, \left\{ \box41 \right. $}  

\noindent $ \box44 $

% \bigskip \smallskip \smallskip \vskip 3,5cm  \smallskip  \qquad  \qquad   
% \special{illustration  fig.2.3.epsf scaled  600}  \smallskip  \smallskip 
\bigskip 
\centerline { \epsfysize=5,0cm    \epsfbox  {fig.2.3.epsf} }
\smallskip  \smallskip

% \bigskip \smallskip \smallskip \vskip 3,5cm  \smallskip  \qquad  \qquad   
% \special{illustration  fig.2.3.bis.epsf scaled  600}  \smallskip  \smallskip 
% \bigskip 
\centerline { \epsfysize=4,5cm    \epsfbox  {fig.2.3.bis.epsf} }
\smallskip  \smallskip

\noindent   {\bf Figure 2.3}	\quad  {\it 	Resolution of the partial Riemann problem 
$\,\,P(W_{\rm left},\, {\cal M}_{\rm right}) \,\,$  in the state space  $ \, \Omega
\,$  and in the space-time plane  for a manifold $\, {\cal M}_{\rm right}  \,$ of
codimension $\, p=2. \,$ }
\smallskip  \smallskip 

% \vfill \eject ~ $\,$  \smallskip
% \vskip 4,5cm  \smallskip  \qquad  \qquad   
% \special{illustration  fig.2.4.epsf scaled  600}  \smallskip  \smallskip 
\bigskip 
\centerline { \epsfysize=6,0cm    \epsfbox  {fig.2.4.epsf} }
\smallskip  \smallskip

% \bigskip \smallskip \smallskip \vskip 3,0cm  \smallskip  \qquad  \qquad   
% \special{illustration  fig.2.4.bis.epsf scaled  600}  \smallskip  \smallskip 
\bigskip 
\centerline { \epsfysize=4,5cm    \epsfbox  {fig.2.4.bis.epsf} }
\smallskip  \smallskip

\noindent   {\bf Figure 2.4}	\quad  {\it  	Resolution of the partial Riemann problem 
$\,\,P(W_{\rm left},\, {\cal M}_{\rm right}) \,\,$   in the state space $ \, \Omega
\,$  and in the space-time plane  for a manifold  $\, {\cal M}_{\rm right}  \,$ of
codimension  $\, p=1. \,$ }
\smallskip  \smallskip 

\smallskip \noindent   $\bullet \qquad \,\,\, $ 
The proof of Theorem 2  follows the construction illustrated on  Figures 2.3 and 2.4
for two particular cases. The idea is to apply the  implicit function theorem  to the
mapping defined as follows. We first introduce a vicinity $\,\, \Delta \,=\, \Theta_1
\times \, \cdots \, \, \times \Theta_p \,\,$ of $\, 0 \, $ in $ \, \R^p \,$~: 

\noindent  (2.4.4) $\qquad \displaystyle 
\R^p \,\, \supset \,\, \Theta_1 \times \, \cdots \, \, \times \Theta_p \,\, \equiv \,\,
\Delta \,\, \ni \, (\epsilon\ib{1} \,,\, \cdots \,,\, \epsilon\ib{p}) \, \equiv
\, \epsilon \,,\qquad \Delta  \in  {\cal V}(0) \,, \, $

\noindent 
and a sub-vicinity $\, \, \Omega_1 \, \subset \Omega_0 \,\,$ of given state $\, W_0 \,$
in such a way that the mapping  $\,\, \Psi_{p}\,\,$ defined by   the chaining of
the {\bf first $\, p \,$ simple waves} is well defined~: 

\setbox21=\hbox {$\displaystyle  \Delta \, \times \, \Omega_1 \,\,\,\ni \,\,\, 
(\epsilon\,,\,W) \, \longmapsto \,   \Psi_{p}(\epsilon,\, W) \,\, \in \, \Omega_0
\,  $}
\setbox22=\hbox {$\displaystyle   \Psi_{p}(\epsilon,\, W) \,\,=\,\,  
\chi\ib{p} \bigl( \epsilon\ib{p} \,;\, \chi\ib{p-1} \bigl( \epsilon\ib{p-1}
\,;\, \cdots \,;\, \chi _{1}(\epsilon\ib{1} \,;\, W) \cdots \bigr) \bigr)
\,.\,    $}
\setbox30= \vbox {\halign{#&# \cr \box21 \cr \box22    \cr   }}
\setbox31= \hbox{ $\vcenter {\box30} $}
\setbox44=\hbox{\noindent  (2.4.5) $\displaystyle  \qquad  
\left\{ \box31 \right. $}  

\noindent $ \box44 $

\smallskip \noindent   $\bullet \qquad \,\,\, $ 
We iterate this mapping with the $\,p \,$ first components $\,\Phi_{1},\, \cdots
\,,\, \Phi_{p} \,$ of the  local chart $\, \Phi \,$ defined at relations (2.3.2)
to (2.3.4)~: 

\setbox21=\hbox {$\displaystyle  \Delta \, \times \, \Omega_1 \,\,\,\ni \,\,\, 
(\epsilon\,,\,W) \, \longmapsto \,   \varphi (\epsilon,\, W) \,\, \in \, \R^p  \,  $}
\setbox22=\hbox {$\displaystyle   \varphi (\epsilon,\, W) \,\,=\,\,  
\bigl( \, \Phi_{1} \bigl( \Psi_{p}(\epsilon,\, W) \bigr) \,,\,  \cdots
\,,\,\Phi_{p} \bigl( \Psi_{p}(\epsilon,\, W) \bigr) \, \bigr) \,. \,  $}
\setbox30= \vbox {\halign{#&# \cr \box21 \cr \box22    \cr   }}
\setbox31= \hbox{ $\vcenter {\box30} $}
\setbox44=\hbox{\noindent  (2.4.6) $\displaystyle  \qquad  
\left\{ \box31 \right. $}  

\noindent $ \box44 $

\noindent 
We construct a  solution of the partial Riemann problem composed by the conservation
law  (2.1.2) and the constraints (2.4.1) for the initial conditions  with the $p$ first
waves issued from the left state $\, W_l\,. \, $ In other terms, we search   a
right state $\, W_r \in {\cal M}_r \,$ as in the  Riemann problem $\,\, R(W_l,\,W_r)
\,\,$  {\it id est}~ under the form 

\noindent  (2.4.7) $\qquad \displaystyle 
W_r \,\,= \,\, \Psi_{p}(\epsilon,\, W_l) \,\,, \qquad W_r \in {\cal M}_r \,. \,$ 

\noindent 
The determination of the state $\,\, W_r \,\,$ satisfying the conditions (2.4.7) is
equivalent to the research of parameter  $\, \epsilon \in \Delta \subset \R^p \,\,$
that satisfy the equation 

\noindent  (2.4.8) $\qquad \displaystyle 
\varphi(\epsilon,\, W_l) \,=\, 0 \,. \,$ 

\smallskip \noindent   $\bullet \qquad \,\,\, $ 
We have the  following natural property~:  

\noindent  (2.4.9) $\qquad \displaystyle 
\Psi_{p}(0,\, W_0) \,\,= \,\, W_0 \, \in \, {\cal M}_r \,\,, \qquad \varphi(0,\,
W_0) \,\, = \,\, 0 \,\,\, \in \, \R^p \,$ 

\noindent
and moreover from Hypothesis 2, the $\,\,  m \times m \,\, $ matrix composed by the
family $\,\,  \Sigma(W_0) \,=\, \bigl( R_{1}(W_0) \,,\, \cdots \,,\, R_{p}(W_0)
\,,\, \xi_{p+1}(W_0) \,,\, \cdots \,,\, \xi_{m}(W_0) \bigr) \,\, \,$ has a rank
equal to $\,m .\,$ Then the same property holds after applying the linear one to one
mapping $\, {\rm d}\Phi(W_0) \,$~: 

\noindent  (2.4.10) $\quad \displaystyle 
{\rm rank } \, \bigl( \, {\rm d}\Phi(W_0)\, {\scriptstyle \bullet} \,  R_{1}(W_0)
\,,\, \cdots \,,\,   {\rm d}\Phi(W_0)\, {\scriptstyle \bullet} \,  R_{p}(W_0) \,,\,
e_{p+1}  \,,\, \cdots \,,\,  e_{m} \, \bigr) \,=\, m \,$ 

\noindent 
due to the definition (2.3.5) of tangent vectors $\,\, \xi_{p+1}(W_0) \,,\, \cdots
\,,\, \xi_{m}(W_0) .\,$  In consequence, when we look to the matrix  defined at 
the relation (2.4.10), we observe that  the block composed by the $p$ first lines and
the $p$ first columns at the top and  the left of  this matrix has a rank exactly 
equal to $\,p .\,$ Moreover, this $\, p \times p \,$  matrix is exactly equal to $\, 
  {{\partial \varphi}\over{\partial \epsilon}} (0,\,W_0) \, $ and in
consequence this  jacobian matrix is  invertible~: 

\noindent  (2.4.11) $\qquad \displaystyle 
{{\partial \varphi}\over{\partial \epsilon}} (0,\,W_0) \, \,$  is an  invertible $\,\,
p  \times  p \,\,$ matrix.

\smallskip \noindent   $\bullet \qquad \,\,\, $ 
Now the theorem of implicit functions proves,  with the eventual  constraint that 
vicinity $\,\, \Omega_1 \,\,$ may have to be reduced,  that the equation 

\noindent  (2.4.12) $\quad \displaystyle 
\varphi(\epsilon,\, W) \,=\, 0 \,\,$ 

\noindent
admits a unique  solution $\,\, (\epsilon,\, W) \,\,$ in the vicinity of $\,\,
(0,\,W_0) \,$  and it  takes the form $\,\, \epsilon \,=\, \pi(W) \,$~: 

\setbox21=\hbox {$\displaystyle  \exists \,\,\, \Omega_1  \ni  W \, \longmapsto \, 
\epsilon =  \pi(W) \,\, \in \Delta \subset  \R^p \,\,,\quad \pi \in {\cal
C}(\Omega_1, \, \Delta) \,\,\,  $  such that  }
\setbox22=\hbox {$\displaystyle  \bigl( \, \varphi (\epsilon,\, W) \,=\, 0
\,,\,\, W \in \Omega_1 \, \bigr)   \,\, \Longrightarrow \,\, \bigl( \, \epsilon 
\,=\, \pi(W) \, \bigr) \,. \,$ }
\setbox30= \vbox {\halign{#&# \cr \box21 \cr \box22    \cr   }}
\setbox31= \hbox{ $\vcenter {\box30} $}
\setbox44=\hbox{\noindent  (2.4.13) $\displaystyle  \quad  
\left\{ \box31 \right. $}  

\noindent $ \box44 $

\noindent 
The general structure of the solution of equation  (2.4.8) (or of the equivalent
equation (2.4.12)) is a consequence of the implicit function theorem   presented in
(2.4.13)~: it  allows the determination of the strengh $\, \epsilon \,$ of the waves
as a function of left state $\, W_l .\,$  Then Theorem 2  is established.   $
\hfill \square \kern0.1mm $

\smallskip \noindent   $\bullet \qquad \,\,\, $ 
The extension of the previous notion   to a partial Riemann problem  $\, \, P({\cal
M}_l ,\, $ $W_r) \, \,$ composed by a boundary manifold $\,\,  {\cal M}_{\rm left} 
\equiv  {\cal M}_l \,\,$  and a state $\,\, W_{\rm right} \equiv  W_r \,\,$   is
defined by the conservation law (2.1.2) and the initial conditions 

\setbox21=\hbox {$\displaystyle  \in  {\cal M}_l \,\,,\quad \qquad x \,< \, 0 \,\, $}
\setbox22=\hbox {$\displaystyle  = \,\, W_r \,\,, \,\, \,\, \qquad x \,> \, 0 \,.\, $}
\setbox30= \vbox {\halign{#&# \cr \box21 \cr \box22    \cr   }}
\setbox31= \hbox{ $\vcenter {\box30} $}
\setbox44=\hbox{\noindent  (2.4.3) $\displaystyle  \qquad   W(x \,,\,0) \, \,\,
\left\{ \box31 \right. $}  

\noindent $ \box44 $

\noindent
If $\, {\rm codim }({\cal M}_l) = p ,\,$ the construction proposed at Theorem 2 can
be extended without difficulty. For each state $\, W_r \,$ sufficiently close to the
manifold $\, {\cal M}_l ,\,$ the partial Riemann problem $\, P({\cal M}_l ,\, W_r) \,$
admits an entropy solution composed by the {\bf last}  $\, p \, $ waves of the Riemann
problem~: 

\setbox21=\hbox {$\displaystyle  W^0  \in  {\cal M}_l \,,\,\, W^{1} \,=\,
\chi\ib{m-p+1}(\epsilon\ib{m-p+1},\, W^0 ) \,,\,\, \cdots \,, \,\,  $}
\setbox22=\hbox {$\displaystyle    W^{p-1} \,\,= \,\, \chi\ib{m-1}(\epsilon\ib{m-1}
,\, W^{p-2}) \,,\,\,\,  W^p \,\,= \,\, \chi\ib{m}(\epsilon\ib{m} ,\,
W^{p-1})  \,=\, W_r \,.\,  $}
\setbox30= \vbox {\halign{#&# \cr \box21 \cr \box22    \cr   }}
\setbox31= \hbox{ $\vcenter {\box30} $}
\setbox44=\hbox{\noindent  (2.4.15) $\displaystyle  \quad   
\left\{ \box31 \right. $}  

\noindent $ \box44 $

\bigskip  \noindent {\smcaps 2.5 } $ \, $ {\bf Partial Riemann problem with an
half-space. }

\noindent   $\bullet \qquad \,\,\, $ 
The notion of partial Riemann problem can be extended to the particular situation of
a  ``half space'' defined as follows. We first introduce the set $\, \, {\cal
M}_{right} \,\equiv \, {\cal M}_{r} \,$ defined by 

\noindent  (2.5.1) $\quad \displaystyle 
 {\cal M}_{r} \,\,= \,\, \{ \, W \in \Omega \,,\,\, \lambda_{1}(W) \, \geq \, 0 \,
\} \,. \,$ 

\noindent 
The  partial Riemann problem $\, P(W_l \,,\,  {\cal M}_{r} ) \,$  between the
particular state $\, W_l \,$ and the half space $\, {\cal M}_{r} \,$ is still defined
by the partial differential equation (2.1.2) and the initial constraints (2.4.1). We
have the 

\smallskip \noindent  {\bf Proposition 2.  }

 $\hfill$   {\bf  Partial Riemann problem with a particular  half space.}

 \noindent 
We suppose that the ~1-field is genuinely nonlinear. A particular solution of  partial
Riemann problem (2.2.1)-(2.4.1) between left state $\, W_l \,$ and  right  half space
$\,\,  {\cal M}_{r} \,\,$ defined in  (2.5.1) can be constructed as follows~:

\noindent ${\scriptstyle \bullet} \quad$  if $\,\, \lambda_{1} (W_l) \, \geq \, 0
\,,\,\, $ then $\,\, W_l \in  {\cal M}_{r} \,\,$ and $\,\, \R \times [0,\, +\infty [
\, \ni   (x,\,t) \longmapsto W(x,\,t) \in \Omega \,$ is equal to a constant state~: 

\noindent  (2.5.2) $\quad \displaystyle 
W(x,\,t) \,\, \equiv \,\, W_l \,\,,\qquad \lambda_{1} (W_l)  \, \geq \, 0 \,. \,$ 

\noindent ${\scriptstyle \bullet} \quad$  if $\,\, \lambda_{1}  (W_l) \, < \, 0
\,,\,\, $ the following rarefaction wave~: 

\setbox21=\hbox {$\displaystyle  W_l \,\, \qquad  \qquad \qquad \qquad   \,\,\,\,\, 
{\rm if} \quad {{x}\over{t}} \,< \, \lambda_{1} (W_l) \, $  }
\setbox22=\hbox {$\displaystyle   \chi\ib{1}\Bigl({{x}\over{t}}-
\lambda_{1}(W_l) \,,\, W_l \Bigr)\quad \, {\rm if} \,\,  \lambda_{1} (W_l) 
\,\leq  \,{{x}\over{t}} \,\leq  \, 0 \, \, $ }
\setbox23=\hbox {$\displaystyle  \chi\ib{1} \bigl(  -\lambda_{1} (W_l) \,,\, W_l
\bigr)   \qquad \,\,\,\,    {\rm if} \quad  {{x}\over{t}} \,  > \, 0 \,\,$ }
\setbox30= \vbox {\halign{#&# \cr \box21 \cr \box22 \cr \box23   \cr   }}
\setbox31= \hbox{ $\vcenter {\box30} $}
\setbox44=\hbox{\noindent  (2.5.3) $\displaystyle  \qquad  W(x,\,t) \,\,= \,\, 
\left\{ \box31 \right. $}  

\noindent $ \box44 $

\noindent
is an entropy solution of the partial Riemann problem  (2.1.2) (2.4.1) (2.5.1).

\smallskip \noindent   $\bullet \qquad \,\,\, $ 
The proof of Proposition 2 is a direct consequence of the notion of rarefaction wave
explicited at relations (2.1.11) to (2.1.13). When $\, {\cal M}_{r} \,$ is the half
space of ``supersonic outflow'' $\, (\lambda_{1} \equiv u - c \, \geq \, 0 $ for
the Euler equations of gas dynamics), the solution of the partial Riemann problem is
trivial when state $\, W_l\,$ satisfies  this condition. On the contrary, the
solution of the  partial Riemann problem is constructed with the help of a ~1-wave for
linking the  ``subsonic left state'' $\, W_l \,$ ($u_l \,<\, c_l$) to the half space
$\,  {\cal M}_{r} .\,$ The ``end-point'' of the rarefaction is the {\bf sonic} state
$\, \, W_r =  \chi\ib{1} \bigl(  -\lambda_{1} (W_l) \,,\, W_l \bigr)   \in  {\cal
M}_{r} \,$ for particular celerity $\, {{x}\over{t}} \,=\, 0 \,\,$ (see
Figure~2.5).~$ \hfill \square \kern0.1mm $ 

% \smallskip 
%%%%%%%%%%%%%%%%%%%%%%%%%%%%%%%%%%%%%%%%%%%%%%%%%%%%%%%%%%%%%%%%%%%%%%%%%%%%%%%%%%%%%% 
% \bigskip 
\centerline { \epsfysize=6,0cm    \epsfbox  {fig.2.5.epsf} }
\smallskip  \smallskip

% \bigskip  
\centerline { \epsfysize=4,5cm    \epsfbox  {fig.2.5.bis.epsf} }
\smallskip  \smallskip

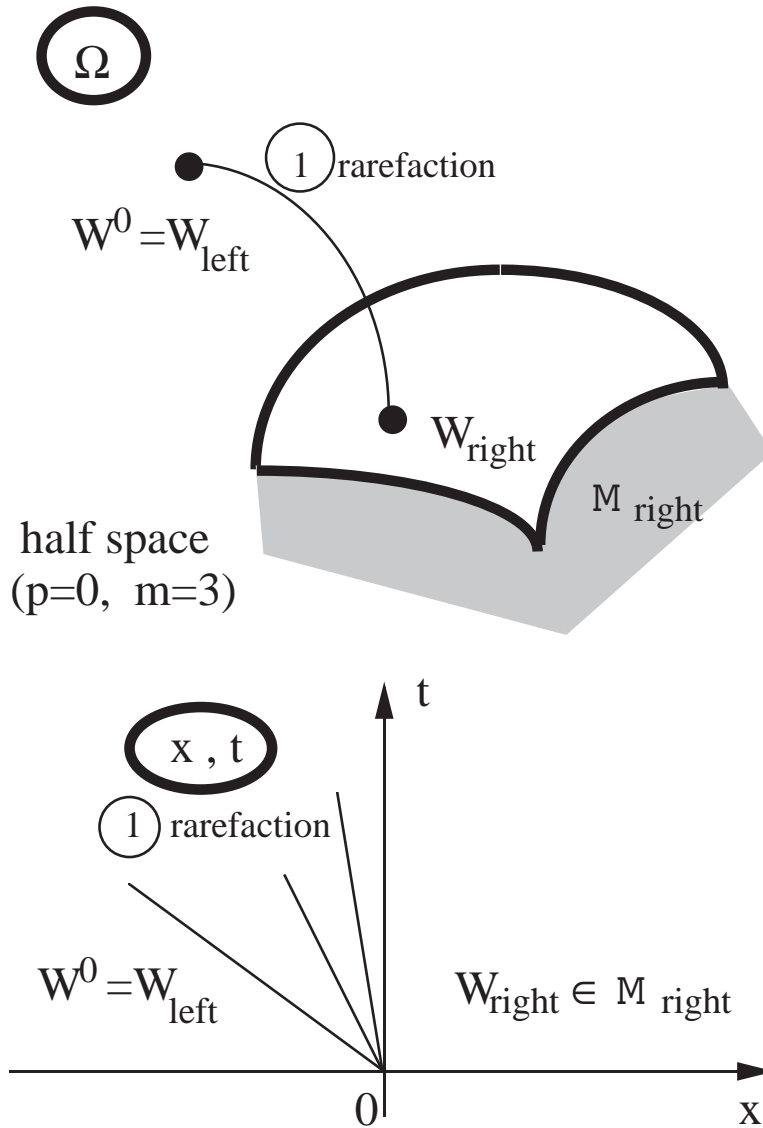
\noindent   {\bf Figure 2.5}	\quad  {\it   	Resolution of the partial Riemann problem 
$\,\,P(W_{\rm left},\, {\cal M}_{\rm right}) \,\,$    in the state space  $ \, \Omega
\,$ and in the space-time plane  for an half space  $\, {\cal M}_{\rm right}  \,$
describing a supersonic outflow. }

% \vfill \eject
%%%%%%%%%%%%%%%%%%%%%%%%%%%%%%%%%%%%%%%%%%%%%%%%%%%%%%%%%%%%%%%%%%%%%%%%%%%%%%%%%%%%%%%
 \bigskip  \bigskip  
\noindent  {\smcaps 3) $ \quad  $Nonlinear boundary conditions for gas dynamics.} 

\noindent   $\bullet \qquad \,\,\, $ 
In this section we use the notion of partial Riemann problem  to consider boundary
conditions for gas dynamics.  We first make the link between this notion and the way
that  linear hyperbolic systems  are well posed in the sense of least squares. Then we
detail some cases where  physically relevant boundary conditions can be interpreted with
particular  partial Riemann problems. 

\smallskip \noindent {\smcaps 3.1 } $ \, $ {\bf  System of linearized Euler
equations.}

\noindent   $\bullet \qquad \,\,\, $ 
Most of the known mathematical results concern hyperbolic {\bf linear} equations. We
linearize the Euler equations of gas dynamics around a {\bf constant state }  $\,
W_{0} \,$ with density $\, \rho_0 ,\,$ velocity $\, u_0 ,\,$ pressure $\, p_0 \,$
and sound celerity $\, c_0 .\,$ We set~: 

\noindent  (3.1.1) $\qquad \displaystyle 
W \,\,= \,\, W_{0} \,+\, W' \, \,$ 

\noindent 
and we neglect second order terms relatively to the variable $\, W' .\,$ In particular,
the  incremental variables $\, \rho' \equiv \rho - \rho_{0} ,\, $ $
\,u' \equiv u - u_{0} ,\,$ $ s' \equiv s - s_{0} \,$  define the incremental vector
$\, Z' \,$ of nonconservative variables~: 

\noindent  (3.1.2) $\qquad \displaystyle 
Z' \,\,= \,\, \bigl( \, \rho' \,,\, u' \,,\, s' \, \bigr)^{\rm \displaystyle t} \,$

\noindent
and  the difference of pressure $\,\, p' \equiv p - p_{0}\,$ is given at the first
order as a function of the incremental thermodynamic  variables $\,  \rho'\,$ and $\, s'
\,$~:

\noindent  (3.1.3) $\qquad \displaystyle 
p' \,\,\equiv \,\, p \,-   p_{0} \,\,= \,\,  c_{0}^2 \, \rho'\,+\, {{ \partial
p}\over{\partial s}}(W_{0}) \, s'  \,.\,$ 

\noindent  
Starting from expression (1.2.15) of the Euler equations of gas dynamics, we get at
the same level of approximation~: 

\noindent  (3.1.4) $\qquad \displaystyle 
{{\partial Z'}\over{\partial t}} \,+\, B(W_{0}) \, {{\partial Z'}\over{\partial x}}
\,\,= \,\,0 \,.\,$ 

\smallskip \noindent   $\bullet \qquad \,\,\, $ 
As in section 1.2, we diagonalize matrix $\,B(W_{0}) \,$ whose expression has been
given at relation (1.2.14) and eigenvectors $\, \widetilde{R}_j(W_{0}) \,$
in (1.2.21). We can express the components $\,  \varphi_{j} \,$ of incremental
vector $\, Z' \,$ in the basis of vectors $\, \widetilde{R}_j(W_{0}) .\,$
These variables $\,  \varphi_{j} \,$ are called the {\bf characteristic variables}~: 

\noindent  (3.1.5) $\qquad \displaystyle 
Z' \,\,= \,\, \sum_{j = 1}^{3} \, \varphi_{j} \,\,  \widetilde{R}_j(W_{0})  \,\,
\equiv \,\, \varphi \,  {\scriptstyle \bullet} \,  \widetilde{R}(W_{0})  \,\,$

\noindent 
and we have from  (1.2.21) the following  expressions~: 

\setbox11=\hbox {$\displaystyle \,\,\varphi_{1} \,\, = \,\, {{1}\over{2
\, \rho_0 \, c_0^2 }} \, \bigl( \, p' \,-\, \rho_0\, c_0 \, u' \, \bigr) $ }
\setbox12=\hbox {$\displaystyle \,\,\varphi_{2} \,\, = \,\, -{{1}\over{c_0^2 }} \, s'
$ }
\setbox13=\hbox {$\displaystyle \,\,\varphi_{3} \,\, = \,\, {{1}\over{2 \,
 \rho_0 \,c_0^2 }} \, \bigl( \, p' \,+\,   \rho_0 \, c_0 \, u' \, \bigr) \,. \, $ }
\setbox40= \vbox {\halign{#&# \cr \box11 \cr  \box12  \cr \box13  \cr }}
\setbox41= \hbox{ $\vcenter {\box40} $}
\setbox44=\hbox{\noindent  (3.1.6) $\qquad \displaystyle \left\{ \box41 \right. $}

\noindent $ \box44 $

\smallskip \noindent   $\bullet \qquad \,\,\, $ 
The change of variables $\,\, \R^3 \, \ni \,  Z' \longmapsto \varphi \, \in \R^3 \, 
\,$ allows to decouple the system of  equations (3.1.4) into three {\bf uncoupled}
advection equations~: 

\noindent  (3.1.7) $\qquad \displaystyle 
{{\partial \varphi}\over{\partial t}} \,+\, \Lambda( W_{0})\, {{\partial
\varphi}\over{\partial x}} \,\,= \,\, 0 \,\,$ 

\noindent 
where $\,\, \Lambda( W_{0}) \, \equiv \, {\rm diag} \, (u_0 - c_0 \,,\, u_0 \,,\, u_0
+ c_0 ) .\, $ The system (3.1.7) is called the characteristic form of the linearized
Euler equations. The above study motivates the mathematical study of linearized
hyperbolic linear systems. 

% \vfill \eject
\smallskip  \smallskip 
\noindent {\smcaps 3.2 } $ \,  $ { \bf  Boundary  problem for
linear hyperbolic systems.} 

\noindent   $\bullet \qquad  \,\,\, $  
In a classical article, Kreiss [Kr70]  developed the notion of {\bf  well posed
problem}  for   initial boundary value problems  in the quarter of space $\,\, x
\leq \, L \,\,$ and  $\, t \geq 0 \,$  associated with linear hyperbolic systems of the
type (3.1.4) or (3.1.7). We suppose that the boundary $\, x=L \,$ is {\bf
non-characteristic}, {\it id est} 

\noindent  (3.2.1) $\qquad \displaystyle 
u_0 - c_0 \, \neq \, 0 \,,\quad  u_0  \, \neq \, 0 \,,\quad u_0 + c_0 \, \neq \, 0
\,, \,$ 

\noindent 
and we denote by $\, \Lambda_0^- \,$ (respectively $\,
\Lambda_0^+ $) the negative part (respectively the positive part) of the nonsingular 
matrix $\, \Lambda( W_{0}) .\,$ We decompose also the characteristic vector $\,
\varphi \,$ into the ingoing  components $\, \varphi^- \,$ and the outgoing components $\,
\varphi^+ \,$~: 

\setbox11=\hbox {$\displaystyle \!\! \varphi \,\,\,\equiv \,\, \varphi^- \,+\,
\varphi^+ \, $}
\setbox12=\hbox {$\displaystyle  \!\!  \varphi^-  = \,\bigl\{ \, (\varphi^-_{j})_{1
\leq  j \leq 3} \,\, ,\,\,  \varphi^-_{j}=\varphi_{j} \,\,\,  {\rm if}
\,  \lambda_j(W_0) < 0 , \,  \varphi^-_{j} =\,0  \,\, {\rm if}
\,  \lambda_j(W_0) > 0 \,  \bigr\} \,$  }
\setbox13=\hbox {$\displaystyle  \!\!  \varphi^+  = \,\bigl\{ \, (\varphi^+_{j})_{1
\leq  j \leq 3} \,\, ,\,\,  \varphi^+_{j} \,=\, \varphi_{j} \,\, \, {\rm if}
\,  \lambda_j(W_0) > 0 , \,  \varphi^+_{j} = 0  \,\,{\rm if}
\,  \lambda_j(W_0) < 0 \,  \bigr\} \,.\, $  }
\setbox40= \vbox {\halign{#&# \cr \box11 \cr  \box12  \cr \box13  \cr }}
\setbox41= \hbox{ $\vcenter {\box40} $}
\setbox44=\hbox{\noindent  (3.2.2) $ \displaystyle \left\{ \box41 \right.$}  

\smallskip \noindent $ \box44 $

%\vfill \eject $\,$ 
% \smallskip \vskip 4,5cm  \smallskip  \qquad  \qquad   
% \special{illustration  fig.3.1.epsf scaled  600}  \smallskip  \smallskip 
\bigskip 
\centerline { \epsfysize=5,0cm    \epsfbox  {fig.3.1.epsf} }
\smallskip  \smallskip

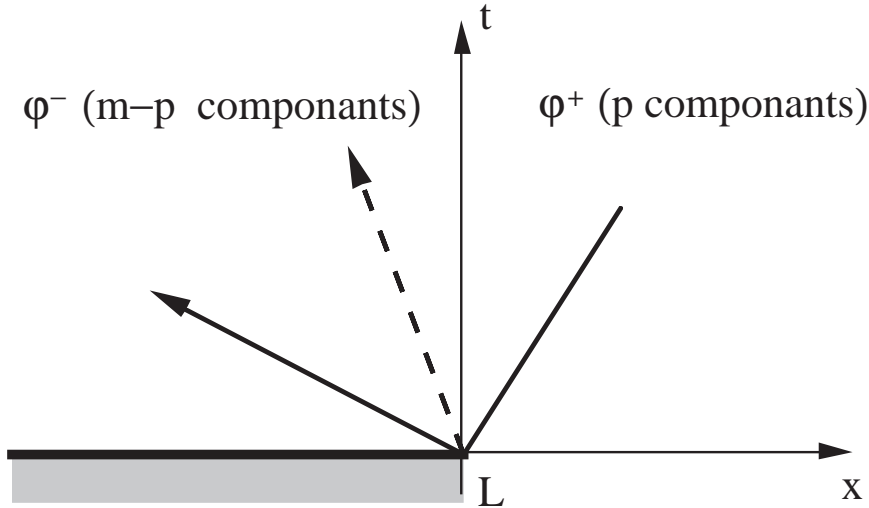
\noindent   {\bf Figure 3.1}	\quad   {\it Characteristic directions 
at the boundary  $\, {x=L} \,$  of the domain  $\, [0,\,L].\,$ }
\smallskip  \smallskip 

\noindent
Due to our choice to consider the  exterior of the domain  ``at the right'' of   the
domain $\,\, \{x\leq L \} \,\,$  of study, the vector $\, \varphi^- \,$ is associated
with negative eigenvalues of matrix $\, \Lambda(W_{0}) \,$ and vector $\, \varphi^+ \,$
corresponds to positive eigenvalues of the same matrix. With the above notations, Kreiss
has proved [Kr70] that in very general situations, the linear hyperbolic system
associated with the initial condition

\noindent  (3.2.3) $\qquad \displaystyle 
\varphi (x,\,0) \,\,= \,\, \varphi_{\,0}(x) \,, \quad \qquad x \leq  L \,$ 

\noindent 
and the boundary condition parameterized by the reflection operator $\Sigma \,$ that
associates to the outgoing characteristics variables $\,  \varphi^+ \,$ a linear
function $\,\, \Sigma \, {\scriptstyle \bullet} \, \varphi^+ \,\, $ which is an
incoming caracteristic variable~: 

\noindent  (3.2.4) $\qquad \displaystyle 
\varphi^-  \,\, = \,\, \Sigma \, {\scriptstyle \bullet} \, \varphi^+ \,\,+\,\, g \,,
\qquad x = L \,$ 

\noindent 
is {\bf well posed} in space $\, L^2 .\,$ 

\smallskip \noindent   $\bullet \qquad \,\,\, $ 
The boundary condition (3.2.4) can be interpreted in terms of {\bf characteristic
directions}~: the field $\,  \varphi^- \,$ along the boundary   is an {\bf affine
function} of the outgoing characteristics $\, \varphi^+  \,$ (see Figure 3.1).  The
previous result remains true   in more complicated situations. In particular, the
multidimensional case can be considered with the same arguments except that direction
$\, x \,$ has to be replaced by the normal direction (see {\it e.g.} Higdon [Hi86]). In the 
characteristic case where {\it e.g.} the reference velocity $\, u_0 \,$ is null, Majda and
Osher [MO75] have extended  Kreiss' result. We remark also that the boundary
condition  can naturally be written in terms of a {\bf boundary manifold}, even if it is
an affine manifold in the present case. We introduce the incremental input vector $\, 
Z^- \,$ and the associated output vector $\,  Z^+ \,$~:

\setbox11=\hbox {$\displaystyle \,\, Z'  \,\,\,\, \equiv \,  \, Z^- \,+\, Z^+ 
\, $}
\setbox12=\hbox {$\displaystyle \,\, Z^- \,\,= \,\, \sum_{\lambda_j(W_0) < 0} 
\varphi_{j} \, \widetilde{R}_{j}(W_{0}) \,\,=\,\, \sum_{j}  \varphi^-_{j} \,
\widetilde{R}_{j}(W_{0})\,$  }
\setbox13=\hbox {$\displaystyle \,\, Z^+ \,\,= \,\, \sum_{\lambda_j(W_0) > 0} 
\varphi_{j} \, \widetilde{R}_{j}(W_{0}) \,\,=\,\, \sum_{j}  \varphi^+_{j} \,
\widetilde{R}_{j}(W_{0}) \,,\, $  }
\setbox40= \vbox {\halign{#&# \cr \box11 \cr  \box12  \cr \box13  \cr }}
\setbox41= \hbox{ $\vcenter {\box40} $}
\setbox44=\hbox{\noindent  (3.2.5) $\qquad \displaystyle \left\{ \box41 \right. $}  

\noindent $ \box44 $

\noindent 
and define the manifold $\, {\cal M}_{right} \equiv  {\cal M}_{r} \,$ by the condition 

\noindent  (3.2.6) $\qquad \displaystyle 
{\cal M}_{r} \,\,= \,\, \bigl\{ \, Z' \, $ given by relations  (3.2.5), $ \,\, \,\, 
\varphi^-  \,=\, \Sigma \, {\scriptstyle \bullet} \, \varphi^+ \,\,+\,\, g \,\bigr\} \,
$ 

\noindent 
and re-write the boundary condition (3.2.4) under the equivalent form~: 

\noindent  (3.2.7) $\qquad \displaystyle 
Z' \, \in \,\, {\cal M}_{r} \, . \,$ 

\smallskip \noindent   $\bullet \qquad  \,\,\, $  
When we consider the boundary condition (3.2.7), we distinguish classically  between
four cases, following the sign of velocity $\, u_0 \,$ and the modulus $\,
\abs{u_0}\,$ compared with sound celerity $\, c_0 .\,$  If $\, u_0 < 0 \,$ at $\,
x=L ,\,$  the fluid enters inside the domain $\, ]-\infty,\,L] \,$ 
and if $\, u_0 > 0 \,$ the boundary $\, \{x=L\}\,$  is an output. If $\, \abs{u_0} \, <
\, c_0 ,\,$ the flow is subsonic whereas if $\, \abs{u_0} \, > \, c_0 ,\,$ it is
supersonic. 

(i) \quad {\bf Supersonic inflow} $\,\, (u_0 \,< \, -c_0) .\,\,$ The outgoing
component of $\, \varphi ,\,$ {\it i.e.} $\, \varphi^+ \,$ is reduced to zero, the
manifold $\,\, {\cal M}_{r} \,\,$ is reduced to the unique  point $\,\,  g\, 
{\scriptstyle \bullet} \,  \widetilde{R}(W_{0})  \,\,$  and $\, \,  {\cal M}_{r}
\,\,$ is of codimension $\, p=3.\,$ The boundary condition (3.2.7) is equivalent to
prescribe all the components of vector $\, Z' .\,$ 

(ii) \quad {\bf Subsonic inflow} $\,\, (-c_0 \,< \, u_0 \,< \, 0) .\,\,$ The number of
incoming characteristics is two and there is one outgoing characteristic
direction. The manifold $\, {\cal M}_{r} \,\,$ is of co-dimension $\, p=2 .\,$ The
linearized problem is well posed when any of  the following  pairs of {\bf  two} 
variables are given (see Oliger-Sundstrom [OS78] or Yee-Beam-Warming [YBW82])~: (density,
pressure), (velocity, pressure) or (enthalpy, entropy) ; we  extend this context  in
the next  section. 

(iii) \quad {\bf Subsonic outflow} $\,\,(0 \,< \, u_0 \,< \,c_0) .\,\,$ Only one
characteristic is going inside the domain of study and two are going outside. 
The manifold $\, {\cal M}_{r} \,\,$ is of co-dimension $\,1 \,$ and only {\bf one}
relation has to be given between the scalar data on the boundary~;  we  remark that the
classical  choice of imposing  the  pressure $\, p \, =\, p_0 \,$ can be written after
linearization~: 

\noindent  (3.2.8) $\qquad \displaystyle 
p' \,\,= \,\, 0 \,.\,$

\noindent 
This condition is equivalent to choose 

\noindent  (3.2.9) $\qquad \displaystyle 
\Sigma \,\,= \,\, \bigl( \, 0 \,,\, -1 \, \bigr) \,,\, \qquad g \,\,= \,\, 0 \,\,$ 

\noindent 
inside relations (3.2.5) and (3.2.6) because the boundary condition (3.2.8) can also
be written under the equivalent form~: 

\noindent  (3.2.10) $\qquad \displaystyle 
\varphi_{1} \,\,= \,\, -\varphi_{3} \,.\,$  

(iv) \quad {\bf Supersonic outflow} $\,\,( u_0 \,> \,c_0) .\,\,$ All the
characteristic directions are going in the direction opposite to the domain $\, \{ x
\leq L \} \,$ ({\it id est} $\, \varphi^- = 0 \,$) and condition (3.2.7) does not
carry any information~: {\bf  no numerical datum} has to be prescribed for a supersonic
outflow. 

\smallskip \noindent   $\bullet \qquad  \,\,\, $  
As a complement of the previous cases, we can add the important case where $\, u_0
\,=\,0 \,$ that corresponds physically to a {\bf  rigid boundary}. The linearized rigid
wall boundary condition takes the form~: 

\noindent  (3.2.11) $\qquad \displaystyle 
u' \,\,= \,\, 0 \,$ 

\noindent 
and the latter boundary condition  is equivalent to prescribe a {\bf  null} mass flux
through the boundary for the linearized equations. Note that the previous condition
(3.2.11) can also be written on the form (3.2.4) with the particular choice~: 

\noindent  (3.2.12) $\qquad \displaystyle
\Sigma \,\,= \,\, \bigl( \, 0 \,,\, 1 \, \bigr) \,,\, \qquad g \,\,= \,\, 0 \,\,$ 

\noindent
and conditions (3.2.11) (3.2.12)  correspond to the relation  

\noindent  (3.2.13) $\qquad \displaystyle 
\varphi_{1} \,\,= \,\, \varphi_{3} \,$  

\noindent
between  the characteristic variables.  It can be shown ({\it e.g.} Oliger-Sundstrom [OS78])
that the the initial boundary value  problem (3.1.4), (3.2.3) and (3.2.12)   is well
posed with the natural $\, L^2 \,$ condition.

\bigskip \noindent {\smcaps 3.3 } $ \,  $ {\bf  Weak Dirichlet nonlinear  boundary
condition. } 

\noindent   $\bullet \qquad  \,\,\, $  
For two particular  systems of  conservation laws of the type 

\noindent  (3.3.1) $\qquad \displaystyle
{{\partial W}\over{\partial t}} \,+\, {{\partial}\over{\partial x}} \, F(W) \,\,=
\,\, 0 \,,\, $ 

\smallskip \noindent
{\it i.e.} for   linear hyperbolic systems  and nonlinear scalar conservation
laws, following ideas developed by T. P. Liu ([Li77], [Li82]) for initial boundary value
problems and transonic flow in nozzles, we have remarked with  Le Floch in [DLF87] and 
[DLF88] that an efficient way to set a boundary condition of the type ``given state $\,
W_{right} \,$  outside the domain of study'', say $\, \{x \leq L \} \,$ to fix the ideas,
is to consider a {\bf weak form} of the Dirichlet boundary condition $\, W \,=\, W_r \,$
in the following sense~: 

\noindent  (3.3.2) $\qquad \displaystyle
W(L^-, t) \, \in \, {\cal B}(W_r) \,.\, $ 

\noindent
In relation (3.3.2), the state  $\, W(L^-, t) \,$ is the limit value of internal state
$\, W(x,\,t) \,$ as $\, x < L \,\,$ tends to the boundary, {\it i.e.} $\, x=L^- \,$ and
$\,  {\cal B}(W_r) \,$ is a set of {\bf admissible states at the boundary} associated
with  datum $\, W_0 .\,$ Even if other formulations are possible and are still in
development (see among others  Nishida-Smoller [NS77],  Bardos, Leroux and  N\'ed\'elec
[BLN79],  Audounet [Au84], Benabdallah [Be86], Benabdallah-Serre [BS87],   [DLF88],
Bourdel, Delorme and Mazet [BDM89],  Gisclon [Gi94], Serre [Se96]), we focus  here on
the use of the {\bf Riemann problem}  for taking into account a physical boundary
conditon, in particuler for gas dynamics. We define the set of admissible states at the
boundary as follows~: 

\setbox11=\hbox{  values at $\,\, {{x}\over{t}}=0^- \,\,$ of the entropic solution }
\setbox12=\hbox { of the Riemann problem $\, R(W,\, W_r) , \, W \in \Omega \,$  }
\setbox40= \vbox {\halign{#&# \cr \box11 \cr \box12 \cr}}
\setbox41= \hbox{ $\vcenter {\box40} $}
\setbox44=\hbox{\noindent  (3.3.3) $\qquad \displaystyle  {\cal B} \bigl( W_r \bigr) 
\,=\, \left\{ \box41 \right\} \, . $}  

\noindent $ \box44 $

\smallskip \noindent   $\bullet \qquad  \,\,\, $  
The advantage of the previous definition  is that the initial  boundary value problem
parameterized by {\bf  constant} boundary datum $\, W_r \,, $  {\bf  constant} 
initial condition $\, W_0 \,$   and set on space-time domain  $\, (x,\,t) \in \,
]-\infty,\, L[ \times [0,\, +\infty[ \,$ by the conditions 

\setbox11=\hbox{$\displaystyle   {{\partial W}\over{\partial t}} \,+\, 
{{\partial}\over{\partial x}} \, F(W) \,\,$  }
\setbox21=\hbox {$ = \,\, 0 \,$  }
\setbox31=\hbox {$ \qquad \qquad x < L \,,\qquad t > 0  \,$}
\setbox12=\hbox {$\displaystyle \,\, W(x,\, 0) \,\,$  }
\setbox22=\hbox {$ = \,\, W_0  \, $}
\setbox32=\hbox {$ \qquad \qquad x < L \,,\qquad t = 0  \,$}
\setbox13=\hbox {$\displaystyle \,\,W(L^- ,\,t) \,$  }
\setbox23=\hbox {$ \in \, {\cal B} \bigl( W_r \bigr) \, $}
\setbox33=\hbox {$ \qquad \qquad t > 0  \,  \,$}
\setbox40= \vbox {\halign{#&#&# \cr \box11 & \box21 & \box31 \cr \box12 & \box22  &
\box32 \cr \box13 & \box23  & \box33 \cr}}
\setbox41= \hbox{ $\vcenter {\box40} $}
\setbox44=\hbox{\noindent  (3.3.4) $\qquad \displaystyle \left\{ \box41 \right. $}  

\noindent $ \box44 $

\noindent
is {\bf well posed} in  conditions analogous to Theorem 1  for the Riemann problem in
Lax's theory (see Section 2.2). Moreover, the set $\, {\cal B} \bigl( W_r \bigr) \,$
is easy to evaluate explicitely. In [DLF88], we have calculated the boundary
set $\, {\cal B} \bigl( W_r \bigr)  \,$ in the particular case of Euler-Saint
Venant equations of isentropic gas dynamics. Even when the state  $\, W_r \,$
corresponds to a supersonic inflow in the linearized analysis, the admissible set 
$\, {\cal B} \bigl( W_r \bigr)  \,$ is reduced to  $\, \{ W_r \} \,$ only in a
vicinity of state  $\,   W_r  .\,$ For large differences between limit state $\,W(L^-
,\,t) \,$ and weakly imposed boundary state  $\,   W_r  \,$ these two states can be
linked together with  an entire family  of waves that compose the Riemann
problem, and state  $\,W(L^- ,\,t) \,$ can even correspond to a supersonic outflow~! 

\smallskip \noindent   $\bullet \qquad  \,\,\, $  
We propose here to extend the previous {\bf  weak} boundary Dirichlet condition (3.3.4) 
to a {\bf  manifold}  $\,\, {\cal M}_r .\,\,$ We first  introduce a new  set $\, \, \beta
\bigl( {\cal M}_r \bigr)  \,\, $ of admissible states  at the boundary~:

\setbox11=\hbox{  values at $\,\, {{x}\over{t}}=0^- \,\,$ of an entropic solution of
the }
\setbox12=\hbox { partial Riemann problem $\, P(W,\, {\cal M}_r) , \, W \in
\Omega \,$  }
\setbox40= \vbox {\halign{#&# \cr \box11 \cr \box12 \cr}}
\setbox41= \hbox{ $\vcenter {\box40} $}
\setbox44=\hbox{\noindent  (3.3.5) $\qquad \displaystyle   \beta \bigl( {\cal M}_r
\bigr)  \,=\, \left\{ \box41 \right\} \, . $}  

\noindent $ \box44 $

\noindent
We remark that this set of admissible states  is a natural extension of definition
(3.3.3) relative to admissible states associated with a single state $\, W_l .\,$  The
initial boundary value problem associated with datum $\,\,  {\cal M}_r\,\,$ is now 
formulated as follows~: 

\setbox11=\hbox{$\displaystyle   {{\partial W}\over{\partial t}} \,+\, 
{{\partial}\over{\partial x}} \, F(W) \,\,$  }
\setbox21=\hbox {$ = \,\, 0 \,$  }
\setbox31=\hbox {$ \qquad \qquad x < L \,,\qquad t > 0  \,$}
\setbox12=\hbox {$\displaystyle \,\, W(x,\, 0) \,\,$  }
\setbox22=\hbox {$ = \,\, W_0  \, $}
\setbox32=\hbox {$ \qquad \qquad x < L \,,\qquad t = 0  \,$}
\setbox13=\hbox {$\displaystyle \,\,W(L^- ,\,t) \,$  }
\setbox23=\hbox {$ \in \,  \beta \bigl( {\cal M}_r \bigr) \, $}
\setbox33=\hbox {$ \qquad \qquad t > 0  \, . \,$}
\setbox40= \vbox {\halign{#&#&# \cr \box11 & \box21 & \box31 \cr \box12 & \box22  &
\box32 \cr \box13 & \box23  & \box33 \cr}}
\setbox41= \hbox{ $\vcenter {\box40} $}
\setbox44=\hbox{\noindent  (3.3.6) $\qquad \displaystyle \left\{ \box41 \right. $}  

\smallskip \noindent $ \box44 $

\smallskip \noindent   $\bullet \qquad  \,\,\, $  
Theorem 2 shows that when constant initial datum $\, W_0 \,$ and boundary manifold $\,
{\cal M}_r \,$ are closed enough, the problem (3.3.6) is well posed in the family of
solutions composed by the $p$ ``first'' waves as described in relation (2.4.2). The
solution of the initial boundary value problem (3.3.6) is the self-similar solution $\,
\, (x,\,t) \longmapsto U({{x-L}\over{t}} \,;\, W_l ,\, {\cal M}_r) \,\,$ described at
relations (2.4.3) in Theorem 2. We remark that only the waves with negative celerities
contribute to problem (3.3.6) even if the partial Riemann problem contains  waves
with positive celerities (see {\it {\it e.g.}} figure 2.3 or~2.4).

\bigskip 
% \vfill \eject 
\noindent {\smcaps 3.4 } $\,${ \bf  Some  fluid boundary conditions. } 

\noindent   $\bullet \qquad  \,\,\, $  
We show in this section that for a family of classical  examples, it is possible to
introduce a partial Riemann problem in order to consider fluid boundary conditions
for the gas dynamics equations. We examine five particular cases~: given state at
infinity or supersonic inflow, subsonic inflow associated with a jet or a nozzle,
subsonic pressure outflow of given static pressure and supersonic outflow. Even if we
still adopt a classical denomination for these boundary conditions, the fact that the
fluid enters into the domain (``inflow'') or goes outside it (``outflow'') has no
influence for the classification of the boundary conditions. The construction of a link
between  physical  data and  mathematical model is done {\it via} a precise choice of a 
boundary manifold.

\smallskip  \noindent   $\bullet \qquad \,\,\, $  
{\bf Given state.} This case corresponds typically  to external aerodynamics problems.
For this kind of flow, a state  $\, W_r = W_{\infty} \,$ is known at a sufficient large
distance  between the object of study and the fluid boundary. In our present model
this particular boundary is located at $\, x=L .\,$  It is natural to set in a weak
sense the boundary condition thanks to a (classical) Riemann problem and we have in
this particular case a  manifold $\,\, {\cal M}^{\rm state}_{W_{\infty}} \,\,$
reduced to a single state~: 

\noindent  (3.4.1) $\qquad \displaystyle
{\cal M}^{\rm state}_{W_{\infty}} \,\,\equiv \,\, \{ \, W_{\infty} \, \} \,. \,$ 

\noindent
We remark that this particular case is mathematically equivalent to the classical 
so-called ``supersonic inflow'' boundary condition. The knowledge of the entire state
$\,  W_{\infty} \,$  is given and it interacts through the boundary with all the
waves of the  Riemann problem. 

%%%%%%%%%%%%%%%%%%%%%%%%%%%%%%%%%%%%%%%%%%%%%%%%%%%%%%%%%%%%%%%%%%%%%%%%%%%%%%%%%%%%%%
% \vfill \eject $\,$  
% \smallskip  \smallskip
\centerline { \epsfysize=3,5cm    \epsfbox  {fig.3.2.epsf} }
% \smallskip  \smallskip

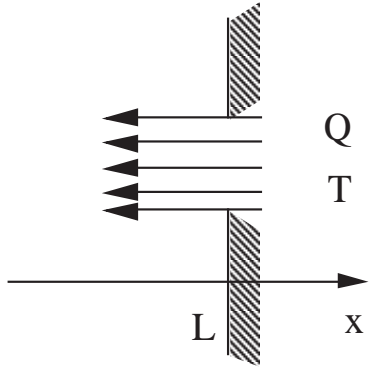
\noindent   {\bf Figure 3.2}	\quad   {\it Subsonic jet inflow at the boundary of the
domain $\, \{x \le L \}.\,$  It is described by a manifold  
$\, {\cal M}^{\rm jet}_{Q, T} \,$  parameterized by  mass flux  $\, Q < 0 \,$  and
temperature~$\, T. \,$ }
% \smallskip  \smallskip
%%%%%%%%%%%%%%%%%%%%%%%%%%%%%%%%%%%%%%%%%%%%%%%%%%%%%%%%%%%%%%%%%%%%%%%%%%%%%%%%%%%%%%

\smallskip  \noindent   $\bullet \qquad \,\,\, $  
{\bf 	Subsonic jet inflow.} In this case introduced in [CDV92] and  described on Figure
3.2, we suppose that the fluid has a given  flow rate $\, Q \,$ and a given temperature
$\, T .\,$ We note that in order to respect the positive sign of the external normal of
the boundary, the mass flux must be chosen negative  $\,\, ( Q<0 ) .\, \,$ We introduce
the boundary manifold $\,\, {\cal M}^{\rm jet}_{Q, T} \,\,$ equal to the set of states $\,
W \,$  that respect exactly the boundary condition~: 

\noindent  (3.4.2) $ \quad \displaystyle
{\cal M}^{\rm jet}_{Q, T} \,\equiv \, \Bigl\{ \, W=(\rho\,,\, q \,,\,
\epsilon)^{\displaystyle \rm t}  \in \Omega ,\,  \quad  q \,=\, Q \,, \quad
\epsilon \,=\, \rho \, C_v \,T \,+\, {1\over2} {{Q^2} \over{\rho}} \, \Bigr\} \,. \,$ 

\noindent
The resolution of the partial Riemann problem $\,\, P(W_l,\, {\cal M}^{\rm jet}_{Q, T}
) \,\,$ can be conducted in the plane of velocity and pressure as it is the case for
the classical Riemann problem. We first express the internal energy $\, e \,$ in terms
of temperature $\, T \,: \,$ 

\noindent  (3.4.3) $\qquad \displaystyle
e \,\, = \,\, C_v \, T \,$ 

\noindent 
and extract density $\, \rho \,$ from mass flow $\, Q \,$ and velocity $\, u \,: $ 

\noindent  (3.4.4) $\qquad \displaystyle
\rho \,\, = \,\, {{Q}\over{u}} \,.  \,$ 

\noindent 
Then the thermostatic law (1.1.2) for polytropic perfect gas  can be written in this
context~: 

\noindent  (3.4.5) $\qquad \displaystyle
p \, u \,\, = \,\, (\gamma-1) \, C_v \, Q \, T \,  \,$ 

\noindent 
and the boundary manifold $\,\, {\cal M}^{\rm jet}_{Q, T} \,\,$  is represented in
velocity-pressure plane with an hyperbola, as shown on Figure 3.3. The resolution of
the partial Riemann problem  in space-time is easy~: manifold  $\,\, {\cal M}^{\rm
jet}_{Q, T} \,\,$ is of codimension 1 and from a given state $\, W_l \,$ in the
vicinity of  $\,\, {\cal M}^{\rm jet}_{Q, T} , \,$ there exists a unique state $\,
W_r \in {\cal M}^{\rm jet}_{Q, T} , \,$ issued from $\, W_l \,$ through a 1-wave. In
the case proposed on Figure 3.4, this wave is a shock wave and the effect of the
interaction at the boundary is the incoming of a shock wave inside the domain $\,
\{x \leq L \} .\,$ 

% \vskip 5,0cm    \qquad  \qquad   
% \special{illustration  fig.3.3.epsf scaled  600}  \smallskip  \smallskip 
\bigskip 
\centerline { \epsfysize=5,0cm    \epsfbox  {fig.3.3.epsf} }
\smallskip  \smallskip

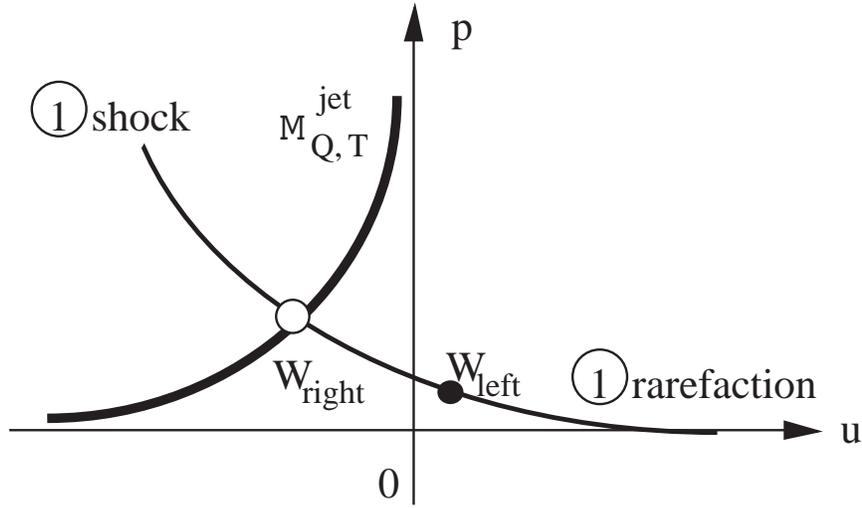
\noindent   {\bf Figure 3.3}	\quad   {\it  Resolution in velocity-pressure plane of the
subsonic jet inflow boundary condition associated with the manifold  $\, {\cal M}^{\rm
jet}_{Q, T} . \,$ }
\smallskip  \smallskip 

%\vfill \eject $\,$ 
% \bigskip \smallskip 
% \vskip 4,0cm  \smallskip  \qquad  \qquad   
% \special{illustration  fig.3.4.epsf scaled  600}  \smallskip  
\bigskip 
\centerline { \epsfysize=5,0cm    \epsfbox  {fig.3.4.epsf} }
\smallskip  \smallskip

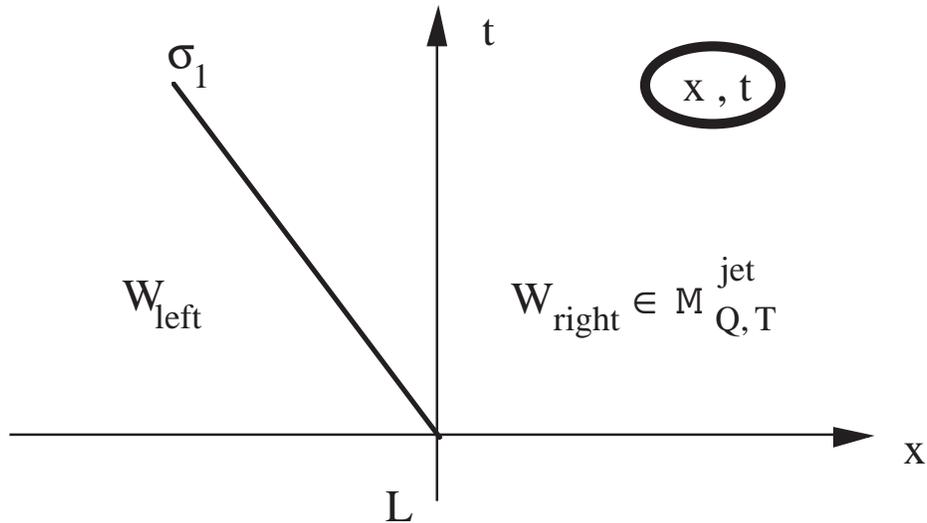
\noindent   {\bf Figure 3.4}	\quad   {\it   Resolution in space-time of the subsonic jet
inflow boundary condition associated with the manifold   $\, {\cal M}^{\rm jet}_{Q, T} .
\,$ }
\smallskip   \smallskip

\smallskip  \noindent   $\bullet \qquad \,\,\, $  
{\bf 	Subsonic nozzle inflow.} This case described on Figure 3.5. A given
tank such the chamber of combustion of an engine to fix the ideas,  contains
some stagnation gas at rest  at high temperature $\, T_s \,$  and high pressure $\,
p_s .\,$  This   gas is conducted to the inflow boundary through a compression wave in
a  convergent nozzle that maintains both total enthalpy $\, H \,$ and specific entropy
$\,\Sigma .\,$  Recall that in a general manner,  total enthalpy is defined according
to 

\noindent  (3.4.6) $\qquad \displaystyle
H \,\, \equiv \,\, {1\over2} u^2 \,+\, {{\gamma \,p}\over{(\gamma-1)\, \rho}} \,$ 

\noindent  
and specific entropy has been used at relation (1.1.12). When stagnation
temperature $\, T_s \,$ and stagnation pressure $\, p_s \,$ are fixed, the total enthalpy
$\, H \,$  and the specific entropy $\, \Sigma \,$ are given by 

\noindent  (3.4.7) $\qquad \displaystyle
H \,\,= \,\, C_p \, T_s \,,\qquad \Sigma \,\,= \,\, C_v \, {\rm Log} \biggl( \, 
\Bigl( {{T_s}\over{T_0}} \Bigr) ^{\gamma} \,   \Bigl( {{p_s}\over{p_0}} \Bigr)
^{1-\gamma} \, \biggr) \,$ 

\noindent
where $\, T_0 \,$ is the  temperature of reference state in relation
(1.1.12).  We define the associated boundary manifold   $\,\, {\cal M}^{\rm
nozzle}_{H,\Sigma} \,\,$ as the set of states such that total enthalpy and entropy
are equal respectively to $\, H \,$ and $\, \Sigma \,$ ; we set 

\setbox11=\hbox{  $  \displaystyle  \,\,  W=(\rho\,,\, q
\,,\, \epsilon)^{\displaystyle \rm t}  \in \Omega \, ,\quad   q \,= \,  \rho \,  u
\,,$  } 
\setbox12=\hbox{  $ \qquad  \epsilon \,=\, {1\over2} \rho \, u^2 \,+\, \rho\,e \,,
\quad  p \,=\, (\gamma-1) \, \rho \, e \, ,  $  }
\setbox13=\hbox {  $  \displaystyle   \qquad   {1\over2} u^2 \,+ \,{{\gamma
\,p}\over{(\gamma-1)\, \rho}} \, =\, H \,,\quad  {\rm Log} \Bigl( {{p \,
\rho_0^{\gamma}}\over{p_0 \, \rho^{\gamma}}} \Bigr) \,= \, {{\Sigma}\over{C_v}} \,$  }
\setbox40= \vbox {\halign{#&# \cr \box11 \cr \box12 \cr \box13 \cr }}
\setbox41= \hbox{ $\vcenter {\box40} $}
\setbox44=\hbox{\noindent  (3.4.8) $\quad \displaystyle  
{\cal M}^{\rm nozzle}_{H,\Sigma} \,\equiv \,\left\{ \box41 \right\} \, . $}  

\noindent $ \box44 $

\bigskip \smallskip 
% \vskip 2,9cm  \smallskip  \qquad  \qquad   
% \special{illustration  fig.3.5.epsf scaled  600}  \smallskip 
\bigskip 
\centerline { \epsfysize=3,9cm    \epsfbox  {fig.3.5.epsf} }
\smallskip  \smallskip

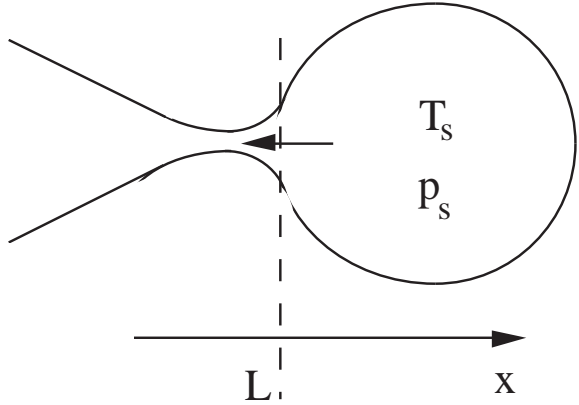
\noindent   {\bf Figure 3.5}	\quad   {\it  Subsonic nozzle inflow at the boundary of the
domain $\, {x 2 L}. \,$ It is described by a manifold $\, {\cal M}^{\rm
nozzle}_{H, \Sigma} \,$  parameterized by total enthalpy  $\, H \,$  and specific
entropy  $\, \Sigma \,$   and in an equivalent way by stagnation temperature  $\, T_s
\,$    and stagnation pressure  $\, p_s . \,$ }
\smallskip  \smallskip 

\noindent
The ``projection'' of manifold $\, \, {\cal M}^{\rm nozzle}_{H,\Sigma} \,\,$ in
velocity-pressure plane is obtained by elimination of density $\, \rho \,$ inside the
relations presented at definition (3.4.8). It comes 

\noindent  (3.4.9) $\qquad \displaystyle
{1\over2} \, u^2 \,+\, {{\gamma}\over{\gamma-1}}\, {{p_0^{1/\gamma}}\over{\rho_0}} \,
{\rm e}^{\Sigma / (\gamma \, C_v)} \, p^{{{\gamma-1}\over{\gamma}}}  \,\,= \,\, H \,$

\noindent
and in terms of stagnation parameters~: 

\noindent  (3.4.10) $\qquad \displaystyle
{{ u^2 }\over{2 \, C_p \, T_s}} \,\,+ \,\, \Bigl( {{p}\over{p_s}}
\Bigr)^{{{\gamma-1}\over{\gamma}}}  \,\, = \,\, 1 \,.\,  $ 

\noindent
The corresponding curve is represented on Figure 3.6  with the particular value of $\,
\gamma = {7\over5}  \,$ {\it i.e.} $\,\, {{\gamma-1}\over{\gamma}} = {2\over7} .\,$ 
We have also represented a graphical resolution of the
partial Riemann problem  $\,\, P(W_l,\,  {\cal M}^{\rm nozzle}_{H,\Sigma}) \,\,$ in
the same plane. 

%\vfill \eject $\,$
% \bigskip \smallskip 
% \vskip 3,5cm  \smallskip  \qquad  \qquad   
% \special{illustration  fig.3.6.epsf scaled  600}  \smallskip
\bigskip 
\centerline { \epsfysize=4,5cm    \epsfbox  {fig.3.6.epsf} }
\smallskip  \smallskip 

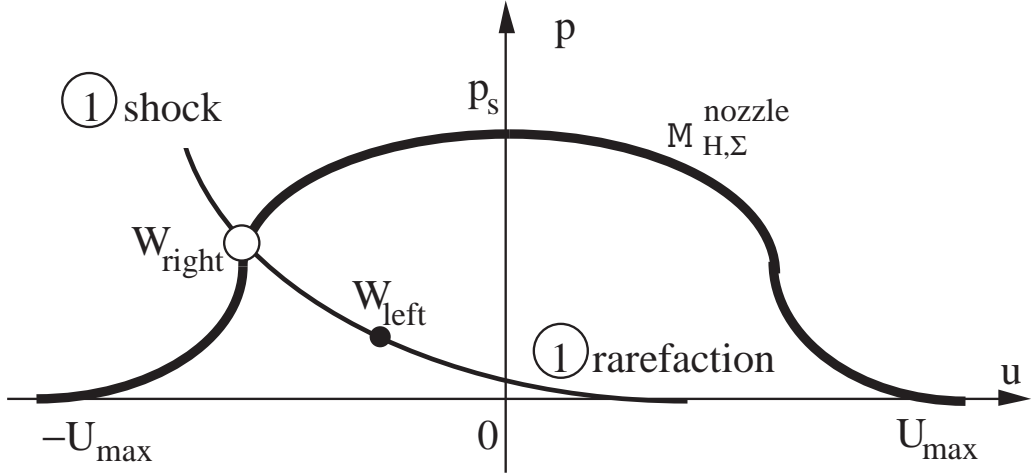
\noindent   {\bf Figure 3.6}	\quad   {\it  	Resolution in velocity-pressure plane of the
subsonic nozzle inflow boundary condition associated with the manifold $\, {\cal M}^{\rm
nozzle}_{H, \Sigma} .\,$   }
\smallskip  \smallskip 

\smallskip  \noindent   $\bullet \qquad \,\,\, $  
{\bf 	Subsonic pressure outflow.} At a subsonic exit, it is classical to prescribe
the static pressure. The linearized analysis shows that only one characteristic
enters inside the domain of study and in consequence only one parameter is a priori
necessary to close the boundary problem.   The  manifold associated to the pressure
datum $\, \Pi \,$ is denoted by  $\,\,  {\cal M}^{\rm pressure}_{\Pi} \,\,$ and is
particularily simple to define~: 

\setbox11=\hbox{  $  \displaystyle  
 W = \bigl( \rho\,,\, q = \rho
\, u \,,\, \epsilon =  \rho \, e \,+\, {1\over2} \rho \,u^2  \bigr)^{\displaystyle \rm
t}  \in \, \Omega ,\,\,  $  } 
\setbox12=\hbox{  $ \qquad \qquad \qquad    (\gamma-1) \, \rho \, e \,=\, \Pi  \,. \,$}
\setbox40= \vbox {\halign{#&# \cr \box11 \cr \box12 \cr  }}
\setbox41= \hbox{ $\vcenter {\box40} $}
\setbox44=\hbox{\noindent  (3.4.11) $\qquad \displaystyle  {\cal M}^{\rm pressure}_{\Pi} 
\equiv   \,\left\{ \box41 \right\} \, . $}  

\noindent $ \box44 $

\noindent  
The resolution of the partial Riemann problem $\,\, P(W_l,\, {\cal M}^{\rm
pressure}_{\Pi} )\,\,$ is presented on Figure 3.7. Note here that an inflow of the
fluid is absolutly compatible with the manifold $\,\,{\cal M}^{\rm pressure}_{\Pi} .\,$ 

% \vfill \eject $\,$
% \bigskip \smallskip 
% \vskip 4,0cm  \smallskip  \qquad  \qquad   
% \special{illustration  fig.3.7.epsf scaled  600}  \smallskip 
% \bigskip 
\centerline { \epsfysize=5,0cm    \epsfbox  {fig.3.7.epsf} }
\smallskip  \smallskip 

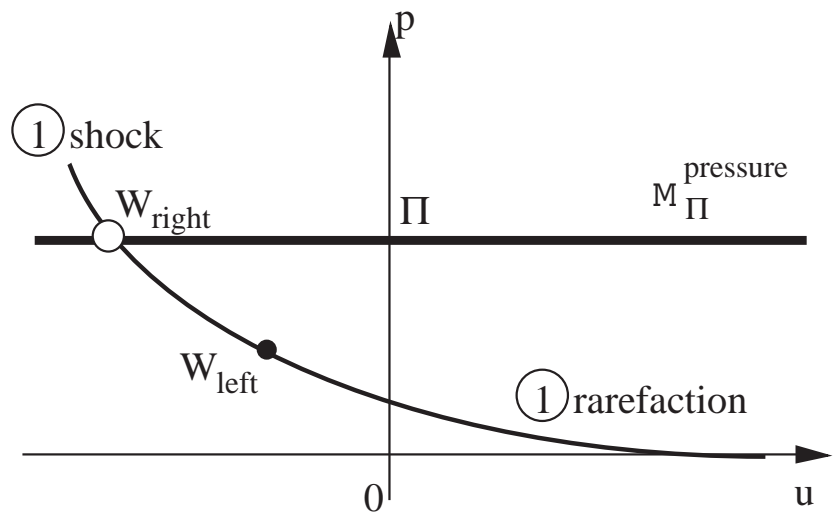
\noindent   {\bf Figure 3.7}	\quad   {\it   Subsonic pressure outflow at the boundary of
the domain $\, {x 2 L}. \,$ It is described by a manifold   $\, {\cal M}^{\rm
pressure}_{\Pi} .\,$   }
\smallskip  \smallskip 

\smallskip  \noindent   $\bullet \qquad \,\,\, $  
{\bf 	Supersonic outflow.} This case corresponds physically to natural boundary
condition at the exit of a De Laval nozzle with sonic neck. No numerical datum has to
be prescribed as we have seen previously with the linearized analysis. Nevertheless,
the fact that the limit state corresponds to a supersonic outflow ({\it i.e.}  $\, u
\geq c )\,$ can enter in conflict with the initial datum for example. In [Du87] and
[Du88], we have analyzed numerically this kind of problem. A natural way to prescibe
a ``supersonic outflow'' boundary condition is to introduce a boundary manifold 
$\,\, {\cal M}^{\rm super} \,\,$  that is now an half space as seen in section 2.5~:

\setbox11=\hbox{  $  \displaystyle  
 W = \bigl( \rho\,,\, q = \rho
\, u \,,\, \epsilon =  {{ \rho \, c^2}\over{\gamma \,(\gamma-1)}} \,+\, {1\over2} \rho
\,u^2  \bigr)^{\displaystyle \rm t}  \in \, \Omega ,\,\, $  } 
\setbox12=\hbox{  $ \qquad \qquad \qquad  \qquad  \qquad   u-c \geq \, 0   \,. \,$}
\setbox40= \vbox {\halign{#&# \cr \box11 \cr \box12 \cr  }}
\setbox41= \hbox{ $\vcenter {\box40} $}
\setbox44=\hbox{\noindent  (3.4.12) $\qquad \displaystyle  {\cal M}^{\rm super} \,
\equiv   \,\left\{ \box41 \right\} \, . $}  

\noindent $ \box44 $

\noindent 
and to take into account the boundary condition with a partial Riemann problem 
 $\,\, P(W_l,\, {\cal M}^{\rm super} )\,\,$ as seen previously at Proposition 2.

\bigskip \bigskip 
% \vfill \eject 
\noindent {\smcaps 3.5 } $ \,  $ {\bf  Rigid wall and moving solid boundary. } 

\noindent   $\bullet \qquad  \,\,\, $ 
The first important case in practice is the surface of some rigid body. In the
multidimensional case, if $\, {\bf n} \,$ is the external normal to the body and $\,
{\bf u} \,=\, (u,\,v) \, $ the bidimensional field of velocity,  this boundary
condition takes the form of an impermeability condition~: 

\noindent  (3.5.1) $\qquad \displaystyle
{\bf u} \, {\scriptstyle \bullet} \, {\bf n} \,\,= \,\, 0 \,. \,  \,$

\noindent
At one space dimension, the external normal is equal to $\, 1 \,$ and the condition
(3.5.1) reduces to a nonlinear version of condition (3.2.11), {\it id est } 

\noindent  (3.5.2) $\qquad \displaystyle
u \,\, = \,\,  0 \,. \,$

\noindent 
As in previous cases, this physically given boundary condition can be incompatible
with the state $\, W(L^-,\,t) \,$ present near the boundary. We propose here to write
in a weak way the boundary condition (3.5.2) with a partial Riemann problem
associated with boundary manifold  $\, \, {\cal M}_r \,= \, {\cal M}^{\rm velocity}_0
\,$~: 

\noindent  (3.5.3) $\qquad \displaystyle
{\cal M}^{\rm  velocity}_0 \,\,\equiv  \,  \, \Bigl\{  \, W = \bigl( \rho\,,\, q = \rho
\, u \,,\, \epsilon \bigr)^{\displaystyle \rm t}  \, \in \Omega ,\,\quad  u \, = \, 0 \,
\Bigr\} \,.  \,$

\smallskip \noindent  {\bf Proposition 3. $\quad$  Stationary state for a rigid wall. }

 \noindent 
Let $\,\, W_l \in \Omega \,\,$ be some given state,  $\,\, {\cal M}^{\rm 
velocity}_0 \,\,$ the manifold defined in (3.5.3) and $\,\, W(x,\,t) \,=\,
U({{x}\over{t}} \,;\, W_l,\,  {\cal M}^{\rm 
velocity}_0 )  \,\,$ the entropy solution of the partial Riemann problem $\,
P(W_l,\,  {\cal M}^{\rm  velocity}_0)  \,$ proposed at Theorem 2.  Then the velocity
$\, u^* \,$ of the  stationary state $\,\,  U(0 \,;\, W_l,\,  {\cal M}^{\rm 
velocity}_0 )  \,\,$ is null and   we have 

\noindent  (3.5.4) $\qquad \displaystyle
 U(0 \,;\, W_l,\,  {\cal M}^{\rm  velocity}_0 ) \,\,= \,\, \bigl( \rho^* \,,\, 0 \,,\,
{{p^*}\over{\gamma-1}}  \bigr)^{\rm \displaystyle t} \,.  \,$

% \vfill \eject $\,$
% \bigskip  \bigskip 
\smallskip 
% \vskip 2,5cm  \smallskip  \qquad  \qquad   
% \special{illustration  fig.3.8.epsf scaled  600}  \smallskip
\bigskip 
\centerline { \epsfysize=4.2cm    \epsfbox  {fig.3.8.epsf} }
\smallskip  \smallskip  

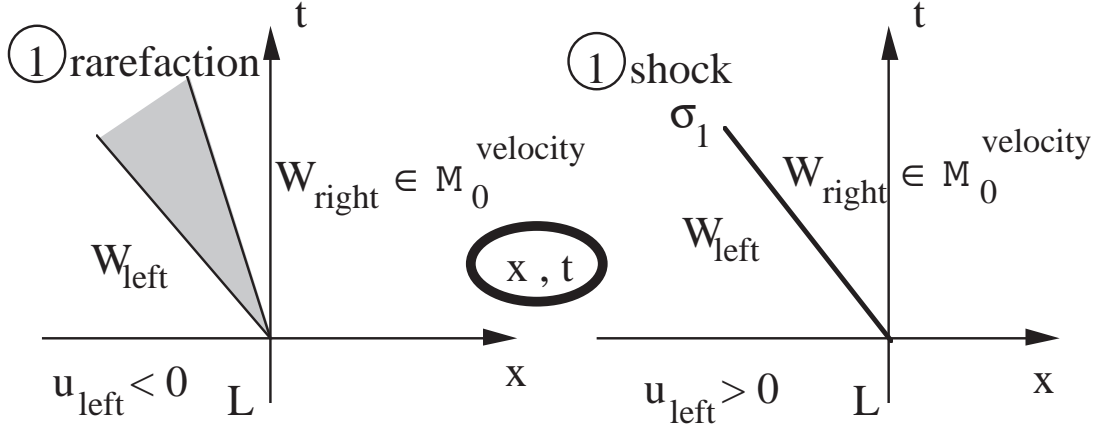
\noindent   {\bf Figure 3.8}	\quad   {\it  	Rigid wall boundary condition. The velocity
of the stationary state is equal to zero. }
\smallskip  \smallskip 

\smallskip \noindent   $\bullet \qquad \,\,\, $ 
The proof of Proposition 3 is clear if velocity $\, u_l \,$ of left state $\, W_l \,$
 is negative. In this case, the 1-wave that links state $\,\,  W_l \, \,$ and manifold
$\,\,  {\cal M}^{\rm  velocity}_0  \,\,$  is a rarefaction and the fan $\,\, ]\mu_1^-
\,,\, \mu_1^+ [ \,\,$  in space-time plane  is defined by the conditions 

\noindent  (3.5.5) $\qquad \displaystyle
u_l - c_l \,\, \equiv \,\, \mu_1^- \,\,  < \,\, \xi \,\,< \,\, \mu_1^+ \,\, \equiv \,\,
u_r - c_r \,=\, -c_r \,\, < \,\, 0 \,$ 

\noindent
as presented on Figure 3.8. Then we have 

\noindent  (3.5.6) $\qquad \displaystyle
U(0 \,;\, W_l,\,  {\cal M}^{\rm  velocity}_0 ) \,\,= \,\, W_r \,$ 

\noindent
and property (3.5.4) is true in this particular case.   If  $\, u_l > 0 ,\,$  the 1-wave
is a shock, we have from (1.6.14) and (1.4.16) and due to the classical
inequality $\,\, p^* \, \geq \, p_l \,$~: 

\noindent  (3.5.7) $\qquad \displaystyle
\sigma_1 \,\,\equiv \,\, u_l \,-\, \sqrt{ {{p^* \,+\, \mu^2 \,p_l}\over{(1\!-\!\mu^2) \,
\rho_l}} } \,\, \leq \,\,  u_l \,-\,  \sqrt{{{1\!-\!\mu^2}\over{\rho_l \,(p^* +\mu^2
\,p_l)}}} \, \bigl( p-p_l \bigr) \,\, \equiv \,\, u^* \,. \,$ 

\noindent
Property  (3.5.4) is a consequence of the general order (3.5.7)  for the different waves
: the celerity $\, \sigma_1 \,$ of the 1-wave is less or equal than the celerity $\, u^*
\,$ of the 2-wave. For the partial Riemann problem $\, P(W_l,\,  {\cal M}^{\rm 
velocity}_0)  ,\,$  we have $\, u^* = 0 \,$ and Proposition 3  is established.  $
\hfill \square \kern0.1mm $

\smallskip \noindent   $\bullet \qquad \,\,\, $ 
A consequence of Proposition 3 is the fact that  the only  wave of partial Riemann
problem $\,\, P(W_l,\,  {\cal M}^{\rm  velocity}_0)  \,\,$  is located in the
quarter of space $\, \{ x < L ,\,\, t > 0 \}\,  \,$ as depicted on Figure 3.8. We remark
also that in this quarter of space, the solution of the partial Riemann problem  $\,\,
P(W_l,\,  {\cal M}^{\rm  velocity}_0)  \,\,$ is {\bf identical} to the solution of the
Riemann problem $\, R( W_l,\, {\widetilde {W_l}})\,$  between left state $\, W_l \,$ and
its ``mirror'' $\,\, {\widetilde {W_l} } .\,$ Recall that, following {\it e.g.}  Roache [Ro72],
the mirror state is  defined by the conditions 

\noindent  (3.5.8) $\qquad \displaystyle
{\widetilde {\rho_l}} \,\,= \,\, \rho_l \,, \qquad 
{\widetilde {u_l}} \,\,= \,\, -u_l \,, \qquad  {\widetilde {p_l}} \,\,= \,\, p_l \,. \,$

\noindent
One advantage of the notion of partial Riemamm problem
is that we do not need to use this very  numerical notion in what follows.

% \vfill \eject $\,$ 
% \vskip 3,0cm  \smallskip  \qquad  \qquad   
% \special{illustration  fig.3.9.epsf scaled  600}  \smallskip 
 \bigskip 
\centerline { \epsfysize=5,5cm    \epsfbox  {fig.3.9.epsf} }
\smallskip  \smallskip  

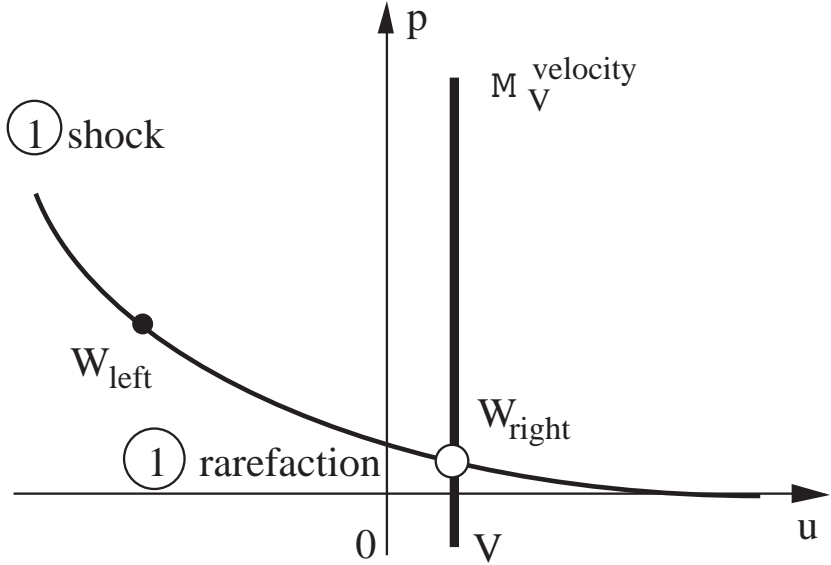
\noindent   {\bf Figure 3.9}	\quad   {\it  	Moving boundary condition at the boundary of
the domain  $\, \{x \le L \}. \,$  It is described by a manifold   $\, {\cal M}^{\rm
velocity}_{V} .\,$}
\smallskip  \smallskip 

% \bigskip \smallskip 
% \vskip 3,5cm  \smallskip  \qquad  \qquad   
% \special{illustration  fig.3.10.epsf scaled  600}  \smallskip 
\bigskip 
\centerline { \epsfysize=4,0cm    \epsfbox  {fig.3.10.epsf} }
\smallskip  \smallskip  

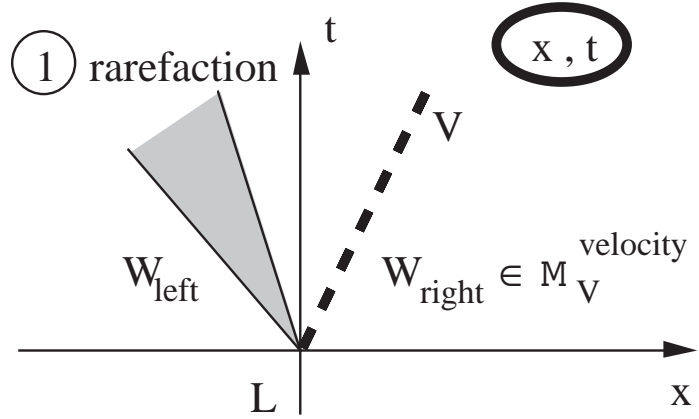
\noindent   {\bf Figure 3.10}	\quad   {\it   Solution in space-time of a moving boundary
condition at the boundary with imposed velocity  $\, V. \,$ }
\smallskip  \smallskip

\smallskip \noindent   $\bullet \qquad \,\,\, $ 
The next example is motivated by fluid-structure interaction problems. In this kind of
configuration, the boundary of the fluid domain is a {\bf  moving boundary} and the
displacements are sufficiently small compared to typical length scales of the problem
and  can be neglected. Nevertheless, the boundary is moving and the value $\, V \,$ of
velocity is supposed to be given. We consider the weak treatment of the boundary
condition $ \, u \, = \, V \,$ by the introduction of the manifold $\,\,  {\cal M}^{\rm 
velocity}_V \,\,$ defined as a generalization of (3.5.3)~:  

\noindent  (3.5.9) $\qquad \displaystyle
{\cal M}^{\rm  velocity}_V \,\,\equiv  \,  \, \Bigl\{  \, W = \bigl( \rho\,,\, q = \rho
\, u \,,\, \epsilon \bigr)^{\displaystyle \rm t}  \, \in \Omega ,\,\quad  u \, = \, V \,
\Bigr\} \,.  \,$

\noindent
The resolution of the partial Riemann problem $\,\, P(W_l,\, {\cal M}^{\rm 
velocity}_V) \,\,$ at the boundary  is straightforward (see [Du99]). The unknown
pressure $\, p^*(V) \,$ on the boundary is solution of the following equations introduced
in (1.6.6)~:  

\setbox21=\hbox {$\displaystyle  V \,- \, u_l \,+\,  \psi( p^*(V) \,;\, \rho_l \,,\,
p_l \,;\, \gamma ) \,\,=\,\, 0 \,\,,\qquad p^*(V) \,< \, p_l \,   $}
\setbox22=\hbox {$\displaystyle  V \,- \, u_l \,+\,  \varphi( p^*(V) \,;\, \rho_l \,,\,
p_l \,;\, \gamma ) \,\,=\,\, 0 \,\,,\qquad p^*(V) \,> \, p_l \, ,  \,  $}
\setbox30= \vbox {\halign{#&# \cr \box21 \cr \box22    \cr   }}
\setbox31= \hbox{ $\vcenter {\box30} $}
\setbox44=\hbox{\noindent  (3.5.10) $\displaystyle  \qquad \left\{ \box31 \right. $} 

\noindent $ \box44 $

\noindent
that can be explicited in this case~:

\setbox21=\hbox {$\displaystyle   V \,- \, u_l \,+\,  {{2\, c_l } \over {\gamma - 1}}
\,  \Bigl[ \, \Bigl( {{p^*(V)} \over {p_l}} \Bigr)^{{\gamma-1}\over{2\,\gamma}} \, - \, 1
\Bigr]  \,\,= \,\, 0  \quad \qquad \quad $if$ \,\,\,\,    V \,- \, u_l\, \geq \, 0 \, $}
\setbox22=\hbox {$\displaystyle  V \,- \, u_l \,+\,  {{\sqrt{2}\, \bigl( p^*(V) - p_l 
\bigr) } \over {\sqrt{\rho_l \bigl[ (\gamma + 1) \, p^*(V) + (\gamma - 1) \, p_l \bigr]}}} 
\,\, \,\,= \,\, 0  \quad $if$ \,\,\,\,    V \,- \, u_l\, \leq \, 0 \,. \, $}
\setbox30= \vbox {\halign{#&# \cr \box21 \cr \box22    \cr   }}
\setbox31= \hbox{ $\vcenter {\box30} $}
\setbox44=\hbox{\noindent  (3.5.11) $\displaystyle  \,\,\, \left\{ \box31 \right. $} 

\noindent $ \box44 $

\noindent
With this approach, it is very easy to precise the difference $\, \, p^*(V) -
p^*(0) \,\, $ of boundary pressures associated with moving boundary
with velocity $\, V \,$ and  a rigid wall.   At the first order, we find~:

\smallskip \noindent (3.5.12) $\qquad \displaystyle
 p^*(V) \,\,=\,\,   p^*(0)  \, - \,  \rho_l \,  c_l \, V \,+\, O(V^2) \,  .\,\,  $

\noindent
We have developed in  [Du99] a complete procedure to take into account numerically a
displacement of the boundary of small amplitude and small celerity  with the so-called
{\bf limiting flux for moving boundary}.

% \vfill \eject
%%%%%%%%%%%%%%%%%%%%%%%%%%%%%%%%%%%%%%%%%%%%%%%%%%%%%%%%%%%%%%%%%%%%%%%%%%%%%%%%%%%%%%%%%%%%%%%
\bigskip  \bigskip   
\noindent  {\smcaps 4)$ \quad $ Application to the finite volume method.}

\noindent  \noindent {\smcaps 4.1 } $ \, $ {\bf  Godunov finite volume method. }

\noindent   $\bullet \qquad \,\,\, $ 
The system of conservation laws

\noindent  (4.1.1) $\qquad \displaystyle
{{\partial }\over{\partial t}} W(x,\,t)\, + \, {{\partial}\over{\partial x}}\,
F(W(x,\,t)) \,\,=\,\,0\, $

\noindent
is now discretized with the finite volume method. We consider a finite one-dimensional
domain of study $\,\, ]0,\,L[ \, \,$ and we divide it with a mesh $\, {\cal T} \,$
initially  composed by a family $\,\,  {\cal S}\ib{\cal T} \,\,$ of $\, J \!+\!1 \,$
vertices (note that integer $\, J \,$ depends on the  mesh  $\, {\cal T} \,$)~:

\noindent  (4.1.2) $\qquad \displaystyle
{\cal S}\ib{\cal T} \,\,= \,\, \bigl\{ \, 0 = x_{0} \,< \, x_{1} \,< \cdots <
x_{j} \,<\, x_{j+1} \, <  \cdots < x_{J} = L \, \bigr\}\, .\, $

\noindent
We construct from vertices $\, \bigl(x_j)_{0 \leq j\leq J} \,$ the 
family $\,\, {\cal E}\ib{\cal T} \,\,$ of finite elements or control volume cells  $\,
K \,$~:

\noindent  (4.1.3) $\qquad \displaystyle
]0,\,L[ \,\,= \,\, \bigcup_{K \in {\cal E}\ib{\cal T}} {\overline K} \,,\quad  K =\,
]x_{j} \,,\, x_{j+1}[ \,\,\, \in  {\cal E}\ib{\cal T} \,,\quad  0 \leq j \leq J \! -
\! 1 \,,\, $

\noindent 
we remark that   two different cells $\, K_1\,$ and $\, K_2 \,$  have an intersection
with null measure~:

\noindent  (4.1.4) $\qquad \displaystyle
{\rm mes} \, (K_1 \cap K_2) \,\,= \,\, 0 \,, \qquad K_1 \not= K_2 \,\, \in  {\cal
E}\ib{\cal T} \,\,  $

\noindent 
and we introduce the two vertices $\, S_-(K) \,$ and $\, S_+(K) \,$ that define the
boundary $\, \partial K \,$ of element $\, K \,: $ 

\noindent  (4.1.5) $\qquad \displaystyle
K \,\,= \,\, ] S_-(K) \,,\, S_+(K) [ \,,\qquad \partial K \,\,= \,\,\bigl\{ \,  S_-(K)
\,,\, S_+(K) \, \bigr\} \,,\quad K \in  {\cal E}\ib{\cal T} \,. \,  $

\smallskip \noindent   $\bullet \qquad \,\,\, $ 
For each element $\, K \, $ of the mesh $\,(K \in {\cal E}\ib{\cal T} ) ,\, $ we
introduce the {\bf mean value} $\, W_K \,$ of the solution of conservation law
(4.1.1)~: 

\noindent  (4.1.6) $\qquad \displaystyle
W_K \,\,=\,\, {{1}\over{\mid K \mid}} \, \int_{K} \, W(x) \, {\rm d} x \,, \qquad  K
\in  {\cal E}\ib{\cal T} \,$ 

\noindent
and we consider the ordinary differential equation satisfied by the functions $ 
[0, +\infty[   $ $ \ni t \longmapsto W_K(t) \in \Omega .\,\,$ We integrate the
conservation law (4.1.1) in space in each  control volume $\, K . \,$ After integrating
by parts the divergent term $\, {{\partial F(W)}\over{\partial x}} ,\,$ it comes 

\noindent  (4.1.7) $\quad \displaystyle
\abs{K} \,  {{{\rm d}W_K}\over{{\rm d}t}} \,+\, F\Bigl( W \bigl(S^+(K),\,t \bigr)
\Bigr)  -  F\Bigl( W \bigl(S^-(K),\,t \bigr) \Bigr)  \, =\, 0  \,,\quad   K \in 
{\cal E}\ib{\cal T} \,. \,$ 

\noindent
The expression (4.1.7) shows that a discretization procedure can be achieved if we are
able to define the flux $\, F_S \,$ for each vertex $\, S \in {\cal S}\ib{\cal T} \,$
in terms of the data, {\it i.e. } the mean values $\, \smash{ \{W_K \,,\, K  \in  {\cal
E}\ib{\cal T} \}} \,$ and of the boundary conditions. This {\bf numerical modelling} 
characterizes  the so-called finite volume method  (see {\it e.g.} Godunov {\it
et al}   [GZIKP79],  Patankar [Pa80], Harten, Lax and Van Leer [HLV83] or  Faille, 
Gallou\"et and  Herbin  [FGH91] among others)  and we will precise some efficient choices
in practice in what follows. 

\smallskip \noindent   $\bullet \qquad \,\,\, $ 
We restrict ourselves for a time to vertices $\,\,  S \! \in \! {\cal S}\ib{\cal T} \,\,$
that are such that  {\bf two} finite elements $\, K_l(S) \,$ and $\, K_r(S) \,$ 
 possess the vertex $\, S \,$  in their boundary. It is the case when vertex $\, S \,$
is internal to the domain $\, ]0,\, L[ \,~:$ 

\noindent  (4.1.8) $\qquad \displaystyle
\partial K_l (S) \, \cap \, \partial K_r(S)  \,= \, \{ S \}  \,, \qquad S \in  {\cal
S}\ib{\cal T} \, $ not located on the boundary. 

\noindent 
In practice, $\, S \,$ is one of the vertices $\, x_{1},\cdots ,\, x_{J-1} \,$ and
$\,\, K_{l}(x_{j}) \,=\, ]x_{j-1}  \,,\, x_{j} [\,,$ $ \,\,  K_{r}(x_{j}) \,=\,
]x_{j} \,,\, x_{j+1} [ \,. \,\,$  With Godunov  [Go59] we propose a numerical
model for  the numerical flux $\, F_S \,$ at vertex $\, S \,$ with the help of the
Riemann problem. We denote by $\, U(\xi \equiv {{x}\over{t}} \,;\, W_l \,,\, W_r ) \,\,$
the entropy solution of the Riemann problem $\, R( W_l ,\, W_r ) . \,$ We set for
internal vertices 

\noindent  (4.1.9) $\qquad \displaystyle
F_{S} \,\,= \,\, F \Bigl( U \bigl( 0\,;\, W_{K_l(S)} \,,\, W_{K_r(S)} \bigr) \Bigr)
\,,\quad \, S \in \bigl\{ \, x_{1},\cdots ,\, x_{J-1} \, \bigr\} \,$ 

\noindent
and by this way, we have constructed  an ordinary differential equation 

\noindent  (4.1.10) $\qquad \displaystyle
 {{{\rm d}W_K}\over{{\rm d}t}} \,+\, {{1}\over{\abs{K}}} \, \Bigl( \, F_{S_+(K)} \,-
\,  F_{S_-(K)} \, \Bigr)  \,\,= \,\, 0 \,,\qquad   K \in  {\cal E}\ib{\cal T} \, \,$ 

\noindent
for control volumes $\, K \,$ such that their boundary $\, \partial K \,$  is not located
on the boundary. 

\smallskip \noindent  \noindent {\smcaps 4.2 } $ \, $ {\bf  Boundary fluxes.}

%%%%%%%%%%%%%%%%%%%%%%%%%%%%%%%%%%%%%%%%%%%%%%%%%%%%%%%%%%%%%%%%%%%%%%%%%%%%%%%%% 
\smallskip 
\centerline { \epsfysize=2,0cm    \epsfbox  {fig.4.1.epsf} }
\smallskip  \smallskip 

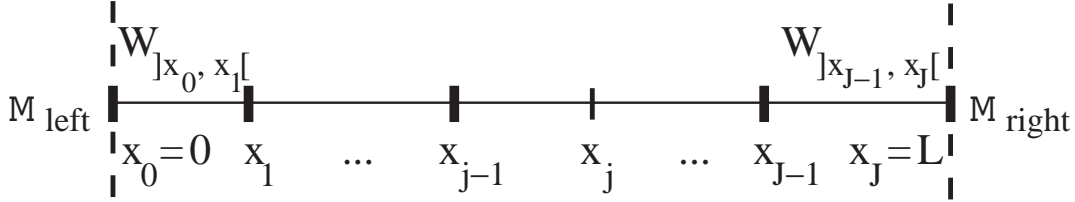
\noindent   {\bf Figure 4.1}	\quad   {\it  Flux boundary conditions for the finite
volume method. At the two boundaries of this unidimensional domain, the partial Riemann
problems have to be set to compute the flux at the boundary.} 
\smallskip \smallskip 
%%%%%%%%%%%%%%%%%%%%%%%%%%%%%%%%%%%%%%%%%%%%%%%%%%%%%%%%%%%%%%%%%%%%%%%%%%%%%%%%%

\noindent   $\bullet \qquad \,\,\, $ 
We wish now to give a sense to equation (4.1.10) for all control volumes $\, K \,$ of
the mesh $\, {\cal T}.\,$ We make the hypothesis that at the left interface $\, \{x=0 \}
\,$  and  at the right interface $\, \{x=L \} \,$ of the domain, a left boundary
manifold $\,\, {\cal M}_l \,\,$ and a right boundary manifold  $\,\, {\cal M}_r \,\,$
are determined by the  physical  data (see Figure 4.1).  Each boundary manifold is
supposed to be chosen inside the family proposed previously  at relations (3.4.1) for
given state $\, W_{\infty} $, (3.4.2) for  subsonic jet inflow $\, {\cal M}^{\rm
jet}_{Q, T} \,$  (with mass flux parameter $\, Q \,$ and temperature $\, T $), (3.4.8)
for  subsonic nozzle inflow $\,  {\cal M}^{\rm nozzle}_{H,\Sigma} \,$ (with numerical
datum defined with  total enthalpy $\, H \,$ and  specific entropy $\, \Sigma $),
(3.4.11) for  subsonic pressure outflow $\, {\cal M}^{\rm pressure}_{\Pi}  \,$ (with
given static pressure $\, \Pi $), (3.4.12) for  supersonic outflow $\, {\cal M}^{\rm
super} \,$ or (3.5.9) for a  moving boundary condition $\,  {\cal M}^{\rm  velocity}_V
\,$ (with celerity~$V $)~:  

\noindent  (4.2.1) $\quad \displaystyle
 {\cal M}_l \,,\,  {\cal M}_r \, \in \, \bigl\{   \{ W_{\infty} \} \,,\, {\cal
M}^{\rm jet}_{Q, T} \,,\,  {\cal M}^{\rm nozzle}_{H,\Sigma} \,,\,  {\cal M}^{\rm
pressure}_{\Pi} \,,\,  {\cal M}^{\rm super} \,,\,  {\cal M}^{\rm  velocity}_V    \bigr\}
\,. \,$

\smallskip \noindent   $\bullet \qquad \,\,\, $ 
We denote by  $\, U(\xi \equiv {{x}\over{t}} \,;\, {\cal M}_l \,,\, W_r ) \,\,$ and 
$\, U(\xi \equiv {{x}\over{t}} \,;\, W_l \,,\,{\cal M}_r ) \,\,$ the selfsimilar
 entropy solution of the partial Riemann problems $\, P({\cal M}_l \,,\, W_r ) \,$ and
$\, P( W_l \,,\,{\cal M}_r ) \,$  proposed at Theorem 2 and relation (2.4.14). The
boundary flux $\, F_{x_0} \,$ at the left interface is evaluated with the resolution of
the partial Riemann problem $\,\, P( {\cal M}_l \,,\,$ $ W_{]x_0 ,\, x_1 [} ) \,\,$ 
between the left boundary datum $\, {\cal M}_l \,$ and the state value $\, W_r \,$ in
the first cell $\, ]x_{0} ,\, x_{1} [ \,$: 

\noindent  (4.2.2) $\qquad \displaystyle
F_{x_{\scriptscriptstyle 0}} \,\,= \,\, F \Bigl( U \bigl( 0\,;\, {\cal M}_l \,,\,
W_{]x_{\scriptscriptstyle 0} ,\, x_{\scriptscriptstyle 1} [} \bigr) \Bigr) \,.\, $

\noindent
In a similar way, the boundary flux $\, F_{x_{\scriptscriptstyle J}}  \,$ at the
right interface $\, \{x=L \} \,$  is obtained thanks to  the resolution of the partial
Riemann problem $\,\,P(W_{]x_{\scriptscriptstyle  J \! - \! 1} ,\,
x_{\scriptscriptstyle J} [} \,,\, {\cal M}_r) \,\,$  between the left
state $\, W_l \,$ in the last cell  $\, ]x_{J-1} ,\, x_{J} [ \,$  and the right
boundary datum~$ \, {\cal M}_r \,$~: 

\noindent  (4.2.3) $\qquad \displaystyle
F_{x_{\scriptscriptstyle J}} \,\,= \,\, F \Bigl( U \bigl( 0\,;\,
W_{]x_{\scriptscriptstyle  J \! - \! 1} ,\, x_{\scriptscriptstyle J} [} \,,\,
{\cal M}_r  \bigr) \Bigr) \,.\, $

\noindent 
With this choice of boundary fluxes $\,\,  F_{x_{\scriptscriptstyle 0}} \,\,$ and
$\,\,  F_{x_{\scriptscriptstyle J}} ,\,$ the ordinary differential equation (4.1.10)
takes  a mathematical sense for all the control volumes $\, K \in {\cal E}\ib{\cal T}
.\,$ 

\smallskip \noindent  \noindent {\smcaps 4.3 } $ \, $ {\bf  Strong nonlinearity at the
boundary. }

\noindent   $\bullet \qquad \,\,\, $ 
We focus here on the fact that the introduction of a partial Riemann problem allows the
treatment of strongly nonlinear effects at the boundary.  Consider to fix the ideas the
boundary $\, \{x=L \} \,$ associated with computational domain $\, ]0,\,L [ . \,$ We
suppose that the physical conditions at this boundary are taken into account with the
help of some manifold $\, {\cal M}_r .\,$ In order to consider  weakly the boundary
condition, we have introduced in (3.3.5)  the set  $\, \beta ( {\cal M}_r )\,$ of
admissible values at the boundary~: 

\setbox11=\hbox{  values at $\,\, {{x}\over{t}}=0^- \,\,$ of the entropic solution  of
the}
\setbox12=\hbox { partial Riemann problem $\, P(W,\,  {\cal M}_r) , \, W \in \Omega \,$ }
\setbox40= \vbox {\halign{#&# \cr \box11 \cr \box12 \cr}}
\setbox41= \hbox{ $\vcenter {\box40} $}
\setbox44=\hbox{\noindent  (4.3.1) $\qquad \displaystyle  \beta \bigl( {\cal M}_r \bigr) 
\,=\, \left\{ \box41 \right\} \, . $}  

\noindent $ \box44 $

\noindent
Then the boundary condition at $\, x=L \,$ is considered in the continuous case in
(3.3.4) and we have by definition~: 

\noindent  (4.3.2) $\qquad \displaystyle
W(L^- ,\,t) \,\in \, \beta \bigl( {\cal M}_r \bigr) \,. \,$

\noindent 
This definition is compatible with the proposed implementation with the finite volume
method~: the boundary flux evaluated in (4.2.3) is the flux of the particular state at
$\, \xi=0 \,$ in the selfsimilar solution $\,\,U(\xi \,;\, W_l  ,\, {\cal M}_r ) \,\,$
of the partial Riemann problem $\, P( W_l  ,\, {\cal M}_r ) . \,$

\smallskip \noindent   $\bullet \qquad \,\,\, $ 
In lot of cases,  the boundary condition acts as in the linear regime, {\it i.e.} the
left state $\, W_l \,$ in the partial Riemann problem $\,\,  P( W_l  ,\, {\cal M}_r ) \,
\,$ lies in the vicinity of the manifold $\, {\cal M}_r \,$ and the stationary value $\,
\, U(\xi=0 \,;\,  W_l  ,\, {\cal M}_r ) \,$ is located  {\bf  inside} the boundary
manifold $\,  {\cal M}_r  \,$~:

\noindent  (4.3.3) $\qquad \displaystyle
U(\xi=0 \,;\,  W_l  ,\, {\cal M}_r ) \,\, \in  {\cal M}_r  \,\,~:  \quad $ weak
nonlinearity at the boundary. 

\noindent
This case occurs typically for a manifold of codimension 1. Consider to fix the ideas a
manifold $\,  {\cal M}^{\rm pressure}_{\Pi} \,$ where the static pressure is fixed by
the physical outflow conditions. The linearized approach in the vicinity of the
(unknown) state $ \,\, W^* \,=\, U(\xi=0 \,;\,  W_l  ,\, {\cal M}_r ) \,\,$  is in
general correct (see Figure 4.2) and the classical approach (see {\it e.g.} Viviand and
Veuillot [VV78], Chakravarthy  [Ch83],  Osher and Chakravarthy [OC83])  consists in
determining a state $\, W^*  \,$ at the boundary that   satifies strongly the boundary
condition,  {\it i.e.}   satisfies the condition  $\,\,  W^*  \! \in \! {\cal M}_r \,\, $
and satisfies also  the linearized gas dynamics equations  written under a characteristic
form along the two characteristics going outside  the domain of study.

% \bigskip \smallskip 
% \vskip 3,0cm  \smallskip  \qquad  \qquad   
% \special{illustration  fig.4.2.epsf scaled  600}  \smallskip 
\bigskip 
\centerline { \epsfysize=4,0cm    \epsfbox  {fig.4.2.epsf} }
\smallskip  \smallskip 

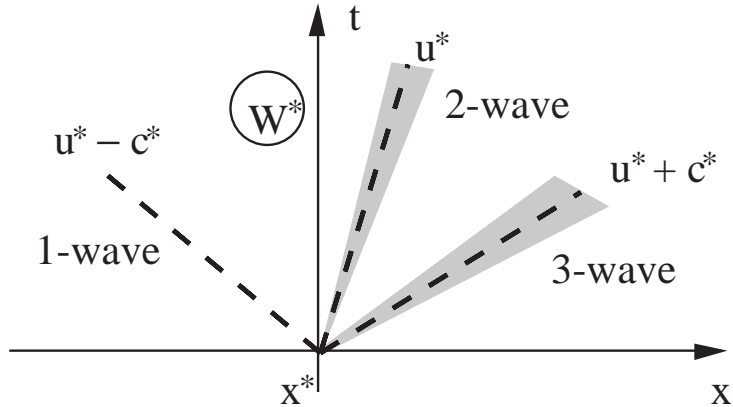
\noindent   {\bf Figure 4.2}	\quad   {\it  Weakly nonlinear boundary condition. The
state  W*  at  the boundary  x=x*  belongs to the boundary manifold. } 
\smallskip \smallskip 

\smallskip \noindent   $\bullet \qquad \,\,\, $ 
This approach is in defect when the limit state is far from the boundary manifold $\, 
{\cal M}_r .\,$ By example, a supersonic boundary condition is prescribed but the
limiting state $\, W(L^-,t) \,$ corresponds to a subsonic outflow. On the contrary, a
subsonic outflow with pressure $ \, \Pi \,$ is supposed to be given but the limiting
state $\, W(L^-,t) \,$ is associated with a supersonic exit ! The flexibility of the
partial Riemann problem allows the treatment with strong nonlinear waves. When
condition (4.3.3) is in defect, the boundary condition (4.2.2) or (4.2.3)  takes into
account the physical data that constructs the interaction at the boundary whereas
these data are not directly used  for the final computation of the boundary flux  (4.2.2)
or (4.2.3).

\bigskip  \vfill \eject 
\noindent {\smcaps 4.4 } $ \, $ {\bf    Extension to  second order
accuracy and to two space dimensions. }

\noindent   $\bullet \qquad \,\,\, $ 
We detail in this section a generalization  for unstructured meshes of the  
``Multidimensional Upwindcentered Scheme for Conservation Laws''  proposed by   Van
Leer [VL79]. At one space dimension on a uniform mesh, it is classical to consider a
scalar field $\,\, z \,\, $ among the primitive variables, {\it i.e.} 

\noindent  (4.4.1) $\qquad \displaystyle
z \,\, \in \,\, \{ \,  \rho \,,\, u \,,\, v \,,\, p \, \} \qquad \, $  
(primitive variables)

\noindent 
and instead of computing the interface  flux with relation (4.1.9), to first construct
two interface states $\,\, W_S^- \,\,$ and $\,\, W_S^+ \,\,$ on each side of the
interface $\, S .\,$ Then the flux is evaluated by the decomposition of the
discontinuity~: 

\noindent  (4.4.2) $\qquad \displaystyle
F_{S} \,\,= \,\, F \Bigl( U \bigl( 0\,;\,  W_S^- \,,\, W_S^+  \bigr) \Bigr)
\,,\quad \, S \in \bigl\{ \, x_{1},\cdots ,\, x_{J-1} \, \bigr\} \,. \,$ 

\noindent 
This nonlinear interpolation is done with a so-called ``slope limiter''  $\,\, \varphi(
{\scriptstyle \bullet}) \,\,$  that operates on each variable proposed in (4.4.1) and
we have typically when a left-right invariance is assumed [Du91]~:

\noindent  (4.4.3) $\quad \displaystyle 
z_{S}^- \,\,\,=\,\,\, z_{j-1/2}  \,\,+ \,\,{1\over2} \, \varphi \biggl( 
 {{z_{j-1/2}-z_{j-3/2}}\over {z_{j+1/2}-z_{j-1/2}}} \biggr)
\, \bigl( z_{j+1/2}-z_{j-1/2}  \bigr) \,, \quad \, S \,=\, x_{j}\,$ 

\noindent  (4.4.4) $\quad \displaystyle 
z_{S}^+ \,\,\,=\,\,\, z_{j+1/2}  \,\,- \,\,{1\over2} \, \varphi \biggl( 
 {{z_{j+3/2}-z_{j+1/2}} \over {z_{j+1/2}-z_{j-1/2}}} \biggr)
\, \bigl( z_{j+1/2}-z_{j-1/2}  \bigr) \,, \quad \, S \,=\, x_{j} \,. \,$

% \bigskip \smallskip 
% \vskip 4,5cm  \smallskip  \qquad  \qquad   
% \special{illustration  fig.4.3.epsf scaled  600}  \smallskip 
\smallskip 
\centerline { \epsfysize=5,0cm    \epsfbox  {fig.4.3.epsf} }
\smallskip  \smallskip 

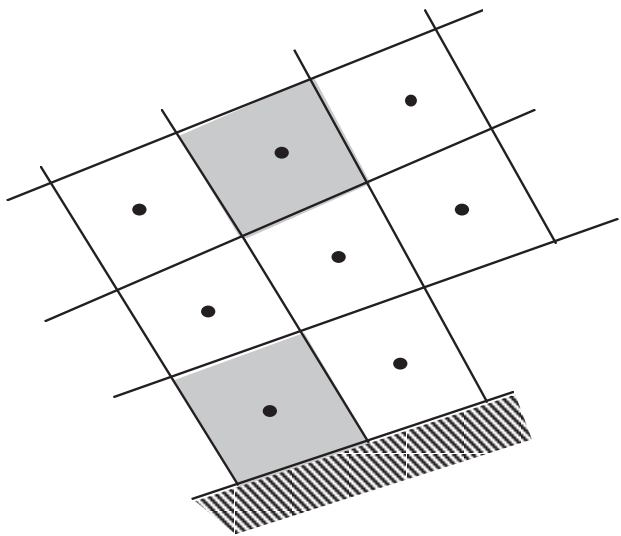
\noindent   {\bf Figure 4.3}	\quad   {\it   Structured cartesian mesh. The control
volumes are exactly the elements of mesh $\,  {\cal T} .\, $} 
\smallskip \smallskip 

% \bigskip \smallskip 
% \vskip 3,5cm  \smallskip  \qquad  \qquad   
% \special{illustration  fig.4.4.epsf scaled  600}  \smallskip
\bigskip 
\centerline { \epsfysize=4,0cm    \epsfbox  {fig.4.4.epsf} }
\smallskip  \smallskip  

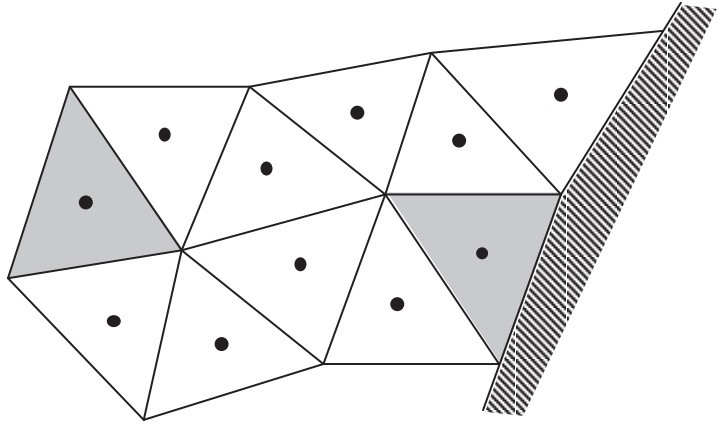
\noindent   {\bf Figure 4.4}	\quad   {\it  	Unstructured mesh composed by triangular
elements.  The control volumes are exactly the elements of mesh  $\,  {\cal T} .\, $} 
\smallskip \smallskip 

\smallskip \noindent   $\bullet \qquad \,\,\, $ 
We focus now on two problems for  the extension to second order   accuracy of Godunov 
finite volume method. One is  associated with the use of unstructured meshes and the
other with the  treatment of boundary conditions. As in the one-dimensional case, the
domain of study is decomposed into finite elements (or control volumes)  $\,\,  K \in
{\cal E}\ib{\cal T} \,\,$  than can be structured in a cartesian way (Figure 4.3) or
with a cellular complex as in Figure 4.4. In both cases,  the intersection of two finite
elements define an interface $\, \, f \! \in  \! {\cal F}\ib{\cal T} .\,\,$ We denote
by   $\, {\bf n}_{f} \,$ the normal at the interface~$\, f \, $ that separates a left
control volume $\, K_l(f) \,$ and a right control volume $ \, K_r(f) .\,$ The ordinary
differential equation (4.1.7) is replaced by a multidimensional version~: 

\noindent  (4.4.5) $\qquad \displaystyle 
\abs{K} \,\, {{{\rm d}W_{K}}\over{{\rm d}t}} \,+\, \sum_{f \subset \partial K}
\,\abs{f} \, \Phi\bigl(  W_{K} \,,\,  {\bf n}_{f} \,,\, W_{ K_r(f)} \bigr)  \,\,=
\,\, 0 \,,\qquad   K \in  {\cal E}\ib{\cal T} \,. \,$ 

\noindent 
For internal interfaces, the function  $\, \Phi\bigl(   {\scriptstyle \bullet} \,,\, 
{\bf n}_{f} \,,\,  {\scriptstyle \bullet}) \,\,$ is equal to the flux of the 
solution at  $\,\, {{x}\over{t}}=0 \,\,$ of the Riemann problem  between states $\, 
W_{ K_l(f)} \, $ and $\,  W_{ K_r(f)} \, $ in the one-dimensional direction along
normal $\,  {\bf n}_{f} \,$  in order to take into account the invariance by rotation
of the equations of gas dynamics (see  [GR96]). We suppose also that for each
interface $\, f \,$ of the boundary, a boundary  manifold $\,\, {\cal M}_f \,$ of
codimension $\, p(f) \,$  is given and   the  normal direction $\,  {\bf n}_{f} \,$ is
by convention external to the domain of study. In consequence,  the control cell $\,K=
K_l(f) \,$ has the face $\, f \, $ in its boundary $\,\, (f \subset \partial ( K_l(f)))
\,\,$ and in relation (4.4.5),  the state  $\,\, W_{ K_r(f)} \,\,$ belongs to $\,\,
{\cal M}_f \,\,( W_{ K_r(f)} \! \in \! {\cal M}_f) \, \,$ and is equal to the state 
$\,\, W^{p(f)} \,\,$ introduced in relation (2.4.2) when solving   the partial Riemann
problem $\,\, P( W_{K_l(f)} \,,\,  {\bf n}_{f} \,,\,  {\cal M}_f ) \,\,$  between
state $\, W_{K_l(f)} \,$ and  manifold~$\,\, {\cal M}_f . \,$

% \bigskip \smallskip 
% \vskip 4,5cm  \smallskip  \qquad  \qquad   
% \special{illustration  fig.4.5.epsf scaled  600}  \smallskip 
\bigskip 
\centerline { \epsfysize=5,0cm    \epsfbox  {fig.4.5.epsf} }
\smallskip  \smallskip 

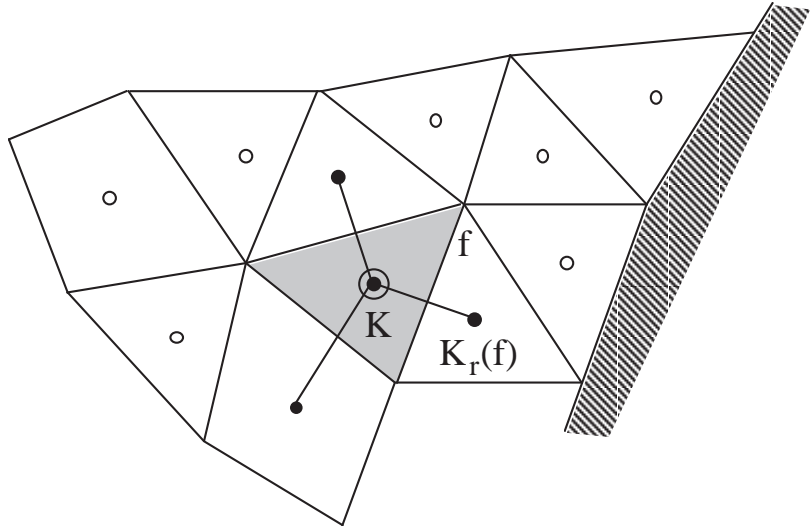
\noindent   {\bf Figure 4.5}	\quad   {\it   	Cellular complex mesh with triangles and
quadrangles. Three neighbouring cells are necessary to determine the gradient in
triangle $\, K \,$  and to limit eventually its variation.} 
\smallskip \smallskip 

\smallskip \noindent   $\bullet \qquad \,\,\, $ 
We consider now a finite element $\, K \,$ internal to the domain. The extension to
second order accuracy of the finite volume scheme consists in replacing the arguments 
$\,  W_{ K_l(f)} \, $ and $\,  W_{ K_r(f)} \, $ in relation (4.4.5) by
nonlinear extrapolations  $\,\,  W_{K_l(f),\,f} \,\,$  and $\,\,  W_{K_r(f),\,f}
\,\,$  on each side of the boundary of  state data and  evaluated as described in what 
follows. We first introduce the set $\, {\cal N}(K) \,$ of neighbouring cells of given
finite element $\, K \in {\cal E}\ib{\cal T} ,\,$ as illustrated on Figure 4.5~: 

\noindent  (4.4.6) $\qquad \displaystyle 
{\cal N}(K) \,\,= \,\, \bigl\{ \, L \in {\cal E}\ib{\cal T} , \quad \exists \, f \in
{\cal F}\ib{\cal T} , \quad  f \, \subset \,\, \partial K \cap \partial L \, \bigr\}
\,.\,$ 

\noindent 
For $\, \, L \! \in \! {\cal N}(K) ,\,$ we suppose by convention that the normal  $\, 
{\bf n}_{f} \,$ to the face $\, f  \subset  \partial K \cap \partial L \,$ is external
to the  element $\, K \,$ {\it id est} $ \,\, K_r(f) = K ,\, K_l(f) = L .\,$ We introduce
also the point $\,\, y_{K,\,f}   \,\,$ on the interface $\,\, f \subset \partial K
\,\, $  that links the barycenters $\, x_{K} \,$ and   $\, x_{K_r(f)} \,$~:

\setbox11=\hbox{   $\,\,  y_{K,\,f}   \,\, \equiv \,\,  (1-\theta_{K,\,f} ) \,
x_{K} \,\,+\,\, \theta_{K,\,f}  \,  x_{K_r(f)}  \,\,,\qquad    y_{K,\,f} \in f
\,, \, $ }
\setbox12=\hbox {  $\qquad \qquad  f  \subset  \partial K \, , \qquad \, K \,$ finite 
 element internal to mesh $\,\, {\cal T} \,. \,$  }
\setbox40= \vbox {\halign{#&# \cr \box11 \cr \box12 \cr}}
\setbox41= \hbox{ $\vcenter {\box40} $}
\setbox44=\hbox{\noindent  (4.4.7) $\qquad \left\{ \box41 \right. \,  $}  

\noindent $ \box44 $

\smallskip \noindent
Then, following Pollet [Po88], for $\, z \,$ equal to one scalar variable of the
 family~: 

\noindent  (4.4.8) $\qquad \,\, \displaystyle
z \,\, \in \,\, \{ \,  \rho \,,\, \rho \,u \,,\, \rho \,v \,,\, p \, \} \qquad \, $  

\noindent
we evaluate a mean value $\, \overline{z_{K,\,f}} \,\, $ on the interface $\,f \,: $ 

\noindent  (4.4.9) $\qquad \,\,  \displaystyle
\overline{z_{K,\,f}} \,\,= \,\,  (1-\theta_{K,\,f} ) \, z_{K} \,\,+\,\,
\theta_{K,\,f}  \,  z_{K_r(f)} \,$

\noindent
and the gradient $\,\, \nabla z(K) \,\,$ of field $\,\, z({\scriptstyle \bullet}) \,
\,$  in volume $\, K \,$ with a Green formula~:

\noindent  (4.4.10) $\quad \displaystyle
\nabla z(K) \,=\, {{1}\over{\mid K \mid}} \, \int_{\partial K} \, \overline{z} \, 
\, {\bf n} \, {\rm d}\gamma \,=\,  {{1}\over{\mid K \mid}} \, \sum_{f \subset
\partial K} \, \abs{f} \,  \overline{z_{K,\,f}}  \, \,\,  {\bf n}_{f} \,,\,\, 
  K \in  {\cal E}\ib{\cal T} \,. \,$ 

\smallskip  \vfill \eject  \noindent   $\bullet \qquad \,\,\, $ 
An ideal extrapolation of field $\,\,  z({\scriptstyle \bullet}) \, \,$ at the
interface $\, f \,$ would be~:

\noindent  (4.4.11) $\qquad \displaystyle
z_{K,\,f} \,\,= \,\, z_{K} \,+\,  \nabla z(K) \, {\scriptstyle \bullet} \, \bigl( 
 y_{K,\,f}  - x_{K} \bigr) \,$
 
\noindent
but the corresponding scheme is unstable as seen by Van Leer [VL77]. When the variation
$\,\,  \nabla z(K) \, {\scriptstyle \bullet} \, \bigl(   y_{K,\,f}  - x_{K} \bigr)
\,\,$  is too important, it has to be ``limited'' as first suggested by Van Leer [VL77].
For doing this in a very general way, we introduce the minimum $\,  m_{K}(z) \,$ 
and the  maximum $\,  M_{K}(z) \,$  of field  $\,\,  z({\scriptstyle \bullet}) \, \,$
in the neighbouring cells~: 

\noindent  (4.4.12) $\qquad \displaystyle
 m_{K}(z) \,\,\,=\,\, {\rm min} \, \bigl\{ \, z_{L} \,, \quad L \in {\cal N}(K) 
 \, \bigr\} \, $

\noindent  (4.4.13) $\qquad \displaystyle
M_{K}(z) \,\,=\,\, {\rm max} \, \bigl\{ \, z_{L} \,, \quad L \in {\cal N}(K) 
 \, \bigr\} \, . \, $

\noindent
If the value $\,  z_{K} \,$ is extremum among the neighbouring ones, {\it i.e.} if 
$\,\, z_{K} \, \leq \,  m_{K}(z) \,,\,$ or $\,\, z_{K}\, \geq \,  M_{K}(z) ,\,\,$
we impose that the interpolated value $\,\, z_{K,\,f} \,\,$ is equal to the cell value
$\,\,  z_{K} \,$~: 

\noindent  (4.4.14) $\quad \displaystyle
z_{K,\,f} \,\,= \,\,  z_{K} \qquad {\rm if} \quad  z_{K} \, \leq \,  m_{K}(z)
\quad {\rm or} \quad   z_{K}\, \geq \,  M_{K}(z) \,\,,\quad f \subset \partial K \,.
\,$ 

\noindent 
When on the contrary $\,\,  z_{K} \,\,$ lies inside the interval $\,\, [ m_{K}(z)
,\,  M_{K}(z)] ,\,$ we impose that the variation $\,\, z_{K,\,f} - z_{K} \,\,$ is
limited by some coefficient $\, k \,( 0 \leq k \leq 1) \,\,$ of the variations $\,\, 
z_{K} - m_{K}(z) \,\,$ and $\,\, M_{K}(z) - z_{K} . \,\,$ We introduce a
nonlinear extrapolation of the field $\,\,  z({\scriptstyle \bullet})\, \,$ between
center $\, x_{K} \,$ and boundary face $\,  y_{K,\,f} \, (f \subset \partial K)
$~: 

\noindent  (4.4.15) $\qquad \displaystyle
z_{K,\,f} \,\,= \,\,  z_{K} \,+\, \alpha_{K}(z) \, \nabla z(K) \,{ \scriptstyle
\bullet} \, \bigl( y_{K,\,f}  - x_{K}  \bigr) \,\,,\qquad f \subset \partial K \,$

\noindent 
with a limiting coefficient $\,\,  \alpha_{K}(z) \,\,$ satisfying the following
conditions~:

\setbox10=\hbox{   $\!\!0 \, \leq  \, \alpha_{K}(z)  \,\, \leq \,\, 1 \,\,, \quad 
z({\scriptstyle \bullet}) \, \,$ scalar field defined in (4.4.8),$\quad K \in {\cal
E}\ib{\cal T} \,$ }
\setbox11=\hbox{   $\!\!  	k \, (z_{K}- m_{K}(z)) \, \leq \,  \alpha_{K}(z)
\,\nabla z(K) \,{\scriptstyle \bullet} \,  \bigl( y_{K,\,f}  - x_{K}  \bigr) \, \leq
\,   	k \, (M_{K}(z) - z_{K}) \,$ }
\setbox12=\hbox { $ \qquad \qquad \qquad \qquad \qquad \qquad \qquad \qquad \qquad
\forall \, f \subset  \partial K \,, \quad K \in {\cal E}\ib{\cal T} \,. \,$ }
\setbox40= \vbox {\halign{#&# \cr \box10 \cr \box11 \cr \box12 \cr}}
\setbox41= \hbox{ $\vcenter {\box40} $}
\setbox44=\hbox{\noindent  (4.4.16) $\,\, \displaystyle  \left\{ \box41 \right. \, $}  

\noindent $ \box44 $

%%%%%%%%%%%%%%%%%%%%%%%%%%%%%%%%%%%%%%%%%%%%%%%%%%%%%%%%%%%%%%%%%%%%%%%%%%%%%%%%%%%%%
% \vskip 3,5cm 
% \bigskip 
\centerline { \epsfysize=4,5cm    \epsfbox  {fig.4.6.epsf} }
\smallskip  \smallskip 

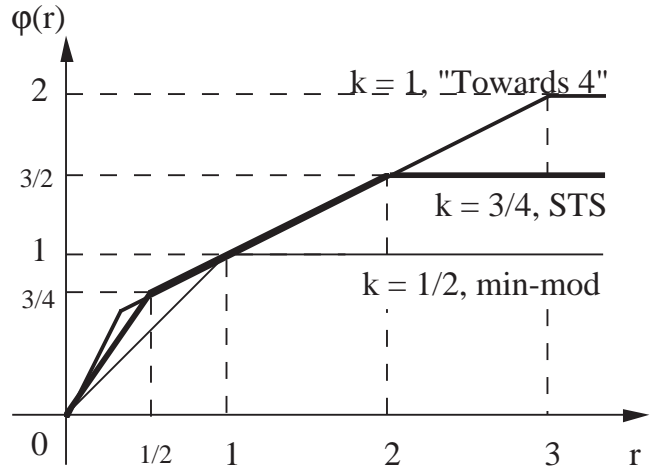
\noindent   {\bf Figure 4.6}	\quad   {\it     Examples of limiter functions that can be
easily extended  to unstructured meshes. } 
\smallskip \smallskip 
%%%%%%%%%%%%%%%%%%%%%%%%%%%%%%%%%%%%%%%%%%%%%%%%%%%%%%%%%%%%%%%%%%%%%%%%%%%%%%%%%%%%%

\noindent
Then $\, \alpha_{K}(z) \,\,$ is chosen as big as possible and inferior or equal to 1 in
order to satisfy the constraints (4.4.16) as displayed on Figure 4.6~: 

\noindent  (4.4.17) $\qquad \displaystyle
\alpha_{K}(z)\,\,\,=\,\,\, {\rm min} \, \biggl[ \, 1 \,,\, k \,\, {{ {\rm
min} \bigl( \, M_{K}(z) - z_{K} \,,\, z_{K}  - m_{K}(z) \, \bigr)  }\over{ {\rm
max} \, \bigl\{ \, \mid \nabla z(K) \, {\scriptstyle \bullet} \, (y_{K,\,f}
 - x_{K} ) \mid \,,\,f \subset \partial K \, \bigr\}  }} \, \biggr] \,. \,$v

\smallskip \noindent   $\bullet \qquad \,\,\, $ 
In the one dimensional case with a regular mesh, it is an exercice to re-write the
extrapolation (4.4.15) under the usual form (4.4.3)  in the context of finite
differences. In this particular case, some limiter functions $\,\, r \longmapsto
\varphi_{k}(r) \,\,$ associated with particular  parameters $\,\, k \,\,$ are shown on
Figure 4.6. For $\,\, k=1 ,\,\,$ we recover the initial limiter proposed by Van Leer in
the fourth paper of the family   ``Towards the ultimate finite difference scheme...'' 
[VL77] ; for this reason, we have named it the ``Towards~4'' limiter (see Figure 4.6).
When $\,\, k={1\over2}  \,\,$ we obtain the ``min-mod'' limiter proposed by Harten
[Ha83]. The intermediate value $\,\, k={3\over4}  \,\,$ is a good compromise between the
``nearly  unstable'' choice $\,\, k=1 \,\,$ and the ``too compressive''  min-mod choice.
We have named it  STS and it has been chosen for our Euler computations in [DM92]. 

%\vfill \eject $\,$
% \bigskip \smallskip 
% \vskip 4,0cm  \smallskip  \qquad  \qquad   
% \special{illustration  fig.4.7.epsf scaled  600}  \smallskip 
\bigskip 
\centerline { \epsfysize=4,5cm    \epsfbox  {fig.4.7.epsf} }
\smallskip  \smallskip 

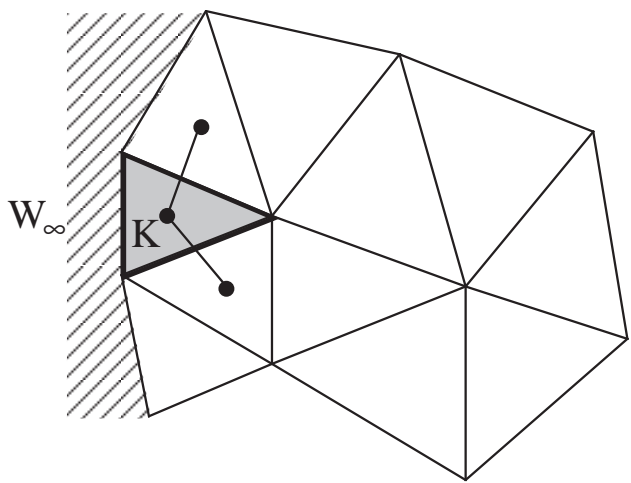
\noindent   {\bf Figure 4.7}	\quad   {\it  	Slope limitation at a fluid boundary.} 
\smallskip \smallskip 

\smallskip \noindent   $\bullet \qquad \,\,\, $ 
We explain now the way the preceeding scheme is adapted near the boundary. We first 
consider  a {\bf fluid boundary}.  When $\, K \,$ is a finite element with some
face $\, g \subset \partial K \,$  lying on the boundary, we still define the set $\,\,
{\cal N}(K) \,\,$  of   neighbouring cells by the relation (4.4.6) as shown on Figure
4.7. The number of neighbouring cells is just less important in this case. Then points
$\,\, y_{K,\,f}\,\,$ are introduced by relation (4.4.7)  if face $\, f \,$ does
not lie on the boundary and by taking the barycenter of face $\,\, g \,\,$ if it is 
lying on the boundary.  The only difference is the way the values $\,\, 
\overline{z_{K,\,g}} \,\,$ are extrapolated for the face $\,\, g \,\,$ that is on the
boundary ; we set 

\noindent  (4.4.18) $\quad \displaystyle
\overline{z_{K,\,g}} \,=\,   \, z_{K} \,\,,\quad g \subset \partial K 
\,\,,\quad g $ face  lying on the boundary of the domain. 

\noindent 
When values $\,\, \overline{z_{K,\,f}} \,\,$ are determined for all the faces $\, \, 
f \subset \partial K ,\,\,$ the gradient $\, \, \nabla z(K) ,\,\,$  the minimal 
 $\,\,  m_{K}(z) \,\, $ and maximal $\,\,  M_{K}(z) \,\, $ values among the
neighbouring cells  are still determined with the relations (4.4.10), (4.4.12) and
(4.4.13) respectively. The constraints (4.4.16) remain unchanged except that no
limitation process is due to the faces lying on the boundary. In a precise way, we set~: 

\noindent  (4.4.19) $  \displaystyle
\alpha_{K}(z) \!=\!{\rm min}  \biggl[  1 ,  { k \quad { {\rm
min} \bigl( \, M_{K}(z) - z_{K} \,,\, z_{K}  - m_{K}(z) \, \bigr)  }\over{ {\rm
max} \, \bigl\{  \abs{ \nabla z(K) \, {\scriptstyle \bullet} \, (y_{K,\,f} - x_{K} )
}, f \subset \partial K,   K_r(f) \!\in\! {\cal N}(K) \bigr\}  }}  \biggr] \,. \,$

\noindent
Then the interpolated values $\,\, z_{K,\,f} \,\,$  for all the faces $\,\, f \subset
\partial K \,\,$ are again predicted with the help of relation (4.4.15). 

\smallskip \noindent   $\bullet \qquad \,\,\, $ 
For a {\bf rigid wall}, the limitation process is a little modified, as presented at
Figure 4.8. First we introduce the limit face $\,g \,$ inside the set of neighbouring
cells~: 

%  \noindent  (4.4.20) $ \displaystyle 
%  {\cal N}(K) = \bigl\{  L \!\in\! {\cal E}\ib{\cal T} , \, \exists f \!\subset \! \partial
%  K \cap \partial L \, \bigr\} \cup   \bigl\{ g \!\in\! {\cal F}\ib{\cal T},\, g \subset 
%  \partial K ,\, g $ on the boundary $ \bigr\}  .\,$ 

\setbox11=\hbox{  $ \displaystyle  {\cal N}(K) = \big\{  L \!\in\! {\cal E}\ib{\cal T} , \, 
\exists f \!\subset \! \partial K \cap \partial L \, \big\}  \,  \cup  $} 
\setbox12=\hbox{ $ \displaystyle \qquad \qquad 
   \cup \,  \bigl\{ g \!\in\! {\cal F}\ib{\cal T},\, g \subset 
\partial K ,\, g $ on the boundary$ \bigr\}  .\, $}  
\setbox40= \vbox {\halign{#&# \cr  \box11 \cr \box12 \cr }}
\setbox41= \hbox{ $\vcenter {\box40} $}
\setbox44=\hbox{\noindent  (4.4.20) $\qquad  \displaystyle  \left\{ \box41 \right.\,$} 
\noindent $ \box44 $ 

%\vfill \eject $\,$
% \bigskip \smallskip 
% \vskip 5,5cm  \smallskip  \qquad  \qquad   
% \special{illustration  fig.4.8.epsf scaled  600}  \smallskip 
\bigskip 
\centerline { \epsfysize=6,0cm    \epsfbox  {fig.4.8.epsf} }
\smallskip  \smallskip 

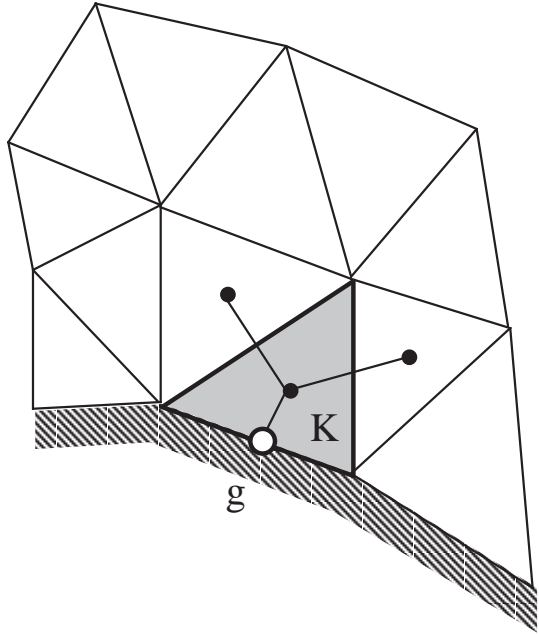
\noindent   {\bf Figure 4.8}	\quad   {\it 	Slope limitation at a solid boundary.}
\smallskip \smallskip 

\noindent
For the face(s) $\,\,  g \subset \partial K \,\,$ lying on the solid boundary, we
determine preliminary values $\,\,  \overline{z_{K,\,g}} \,\,$ by taking in
consideration    at this level the nonpenetrability boundary condition (3.5.1). 
We introduce the two components $\,\, n^x_{g} \,$ and $\,\, n^y_{g} \,$ of the
normal $\,\, {\bf n}_{g} \,\,$ at the boundary and we  set, in coherence with
variables (4.4.8)~: 

\setbox11=\hbox{   $\!\!  \overline{\rho_{K,\,g}} \,\,= \,\, \rho_{K}\, $} 
\setbox12=\hbox{   $\!\!  \overline{\rho_{K,\,g}} \,\,\,  \overline{u_{K,\,g}}
\,\,= \,\, \rho_{K}\, \bigl( u_{K} - ({\bf u}_{K} \, {\scriptstyle \bullet} \,
{\bf n}_{g}) \,  n^x_{g} \,\bigr) \,  $} 
\setbox13=\hbox{   $\!\!  \overline{\rho_{K,\,g}} \,\,\,  \overline{v_{K,\,g}}
\,\,= \,\, \rho_{K}\, \bigl( v_{K} \,- ({\bf u}_{K} \, {\scriptstyle \bullet} \,
{\bf n}_{g}) \,  n^y_{g} \,\bigr) \,  $} 
\setbox14=\hbox{   $\!\!  \overline{p_{K,\,g}} \,\,= \,\, p_{K}\,.\,  $} 
\setbox40= \vbox {\halign{#&# \cr  \box11 \cr \box12 \cr \box13 \cr \box14 \cr}}
\setbox41= \hbox{ $\vcenter {\box40} $}
\setbox44=\hbox{\noindent  (4.4.21) $\qquad  \displaystyle  \left\{ \box41 \right.\,$}  

\noindent $ \box44 $

\noindent
We consider  also these values for the limitation algorithm. We define ``external
values'' $\,\, z_{L} \,\,$ for $\,\,  L=g \,\,$ and face $\,\,g \,\, $  lying on the
boundary as equal to the ones defined in relation (4.4.21)~: 

\noindent  (4.4.22) $\quad \displaystyle 
z_{g} \,\equiv \, \overline{z_{K,\,g}} \,,\quad    z({\scriptstyle \bullet})
\, \,$  field defined in (4.4.21),$\,\,\,\,   g \subset  \partial K  $ on the
boundary.

\noindent
Then the extrapolation algorithm that conducts to relation (4.4.15) for extrapolated
values $\,\,  z_{K,\,f} \,\,$ is used as  in the internal case. 

\smallskip \noindent   $\bullet \qquad \,\,\, $ 
When all values  $\,\,  z_{K,\,f} \,\,$ are known for all control volumes $\,\, K \in
{\cal E}\ib{\cal T} ,\,\,$ all faces $\,\, f \!\subset \! K \,\,$ and all fields $\,\,
z({\scriptstyle \bullet}) \,\,$ defined at relation (4.4.8), extrapolated states  
$\,\,  W_{K,\,f} \,\,$ are naturally defined by going back to the conservative
variables. Then we introduce  these states as arguments of the flux function $\,\,
\Phi\bigl(   {\scriptstyle \bullet} \,,\,  {\bf n}_{f} \,,\,  {\scriptstyle \bullet})
\,\,$ and obtain by this way a new system of ordinary differential equations~: 

\noindent  (4.4.23) $\,\,\, \displaystyle 
\abs{K} \,\, {{{\rm d}W_{K}}\over{{\rm d}t}} \,+\, \sum_{f \subset \partial K}
\,\abs{f} \, \Phi\bigl(    W_{K,\,f}  \,,\,  {\bf n}_{f} \,,\,   W_{K_r(f),\,f} 
\bigr)  \,\,= \,\, 0 \,,\qquad   K \in  {\cal E}\ib{\cal T} \,. \,$ 

\noindent
The numerical integration of such kind of system is just a question of Runge-Kutta
scheme as presented in [CDV92]. We have used with success in [DM92] the Heun scheme of
second order accuracy for discrete integration of (4.4.23) between time steps $\,\,  n\,
\Delta t \,\,$ and $\,\, (n+1) \, \Delta t \,\,$~: 

\noindent  (4.4.24) $  \displaystyle 
{{\abs{K}}\over{\Delta t}} \,  \Bigl( \, \widetilde{ W_{K}} -  W_{K}^{n} \, \Bigr)
 + \! \sum_{f \subset  \partial K}   \abs{f} \, \Phi \Bigl(    W^{n}_{K,\,f}  \,,\, 
{\bf n}_{f} \,,\,   W^{n}_{K_r(f),\,f}  \Bigr)  \,=\, 0 \,,\,\,    K \in  {\cal
E}\ib{\cal T}   \,$ 

\noindent  (4.4.25) $  \displaystyle 
{{\abs{K}}\over{\Delta t}} \,  \Bigl( \, \widetilde{ \widetilde{ W_{K}}}  -  
\widetilde{ W_{K}} \, \Bigr)  + \!  \sum_{f \subset  \partial K}  \abs{f} \,
\Phi \Bigl(   \widetilde{ W} _{K,\,f}  \,,\,  {\bf n}_{f} \,,\,    \widetilde{ W}
_{K_r(f),\,f}  \Bigr)  \,=\, 0 \,,\,\,     K \in  {\cal E}\ib{\cal T} \, \,$ 

\noindent  (4.4.26) $\,\,\, \displaystyle 
W_{K}^{n+1} \, = \, {1\over2} \, \Bigl( \, \widetilde{ \widetilde{ W_{K}}} \, +
\,  W_{K}^{n}  \, \Bigr) \,,\quad   K \in  {\cal E}\ib{\cal T} \,. \,$ 

\smallskip \noindent   $\bullet \qquad \,\,\, $ 
{\it Acknowledgments}. The author thanks Edwige Godlewski and St\'ephanie Mengu\'e for
helpfull comments on the preliminary edition of this work. 

%%%%%%%%%%%%%%%%%%%%%%%%%%%%%%%%%%%%%%%%%%%%%%%%%%%%%%%%%%%%%%%%%%%%%%%%%%%%%%%%%%%%%%%%%%
 \bigskip \bigskip  
\noindent  {\smcaps 5) $ \quad $  References.}

\smallskip \hangindent=7mm \hangafter=1 \noindent  
 [Au84] \quad   J. Audounet. Solutions discontinues
param\'etriques des syst\`emes de lois de conservation et des probl\`emes aux limites
associ\'es, {\it S\'eminaire}, Universit\'e Toulouse 3, 1983-84.

\smallskip \hangindent=7mm \hangafter=1 \noindent  
[BLN79]  \quad    C. Bardos, A.Y. Leroux, J.C.
N\'ed\'elec. First Order Quasilinear Equations with Boundary Conditions, {\it Comm.
Part. Diff. Eqns.}, vol.~4, p.~1017-1034, 1979.

\smallskip \hangindent=7mm \hangafter=1 \noindent  
 [Be86] \quad   A. Benabdallah. Le ``p-syst\`eme'' sur un
intervalle,  {\it Comptes Rendus de l'Acad\'emie des
Sciences,} Paris, tome 303, S\'erie 1, p.~123-126, 1986. 

\smallskip \hangindent=7mm \hangafter=1 \noindent  
  [BS87]  \quad  A. Benabdallah, D. Serre.  Probl\`emes aux 
limites pour les syst\`emes hyperboliques non-lin\'eaires de deux \'equations \`a une
dimension d'espace, {\it Comptes Rendus de l'Acad. des Sciences,} Paris, tome 303,
S\'erie 1, p.~677-680, 1987. 

\smallskip \hangindent=7mm \hangafter=1 \noindent  
 [BDM89]  \quad    F. Bourdel, P. Delorme, P. Mazet. 
Convexity in Hyperbolic Problems. Application to a Discontinuous Galerkin Method for
the Resolution of the Polydimensional Euler Equations, Second International Conference
on Hyperbolic Problems (Ballmann-Jeltsch Editors), {\it Notes on Numerical Fluid
Dynamics,} vol.~24, p.~ 31-42, Vieweg, 1989. 

\smallskip \hangindent=7mm \hangafter=1 \noindent  
 [Ca85]  \quad   H.B. Callen. { \it Thermodynamics and an
introduction to thermostatics, second edition}, John Wiley \& Sons, New York, 1985.

\smallskip \hangindent=7mm \hangafter=1 \noindent  
 [Ch83]  \quad   H.B. Callen.   S. Chakravarthy. Euler solutions, Implicit 
schemes and Boundary Conditions, {\it AIAA Journal}, vol.~21,  n$^{\rm o}$ 5, 
p.~699-706, 1983. 

\smallskip \hangindent=7mm \hangafter=1 \noindent  
 [CDV92]  \quad   D. Chargy, F. Dubois, J.P. Vila.
M\'ethodes num\'eriques pour le calcul  d'\'ecoulements compressibles, applications
industrielles, {\it cours de l'Insti\-tut pour la Promotion des Sciences de
l'Ing\'enieur}, Paris,  septembre 1992.

\smallskip \hangindent=7mm \hangafter=1 \noindent  
 [CF48] \quad    R. Courant, K.O. Friedrichs. {\it 
Supersonic Flow and Shock Waves,}  Interscience, New-York, 1948.

\smallskip \hangindent=7mm \hangafter=1 \noindent  
 [Du87] \quad    F. Dubois.   Boundary Conditions and the
Osher Scheme for the Euler Equations of Gas Dynamics, {\it Internal Report n$^{0}$
170 of Centre de Math\'e\-matiques Appliqu\'ees de l'Ecole Polytechnique,} Palaiseau,
September 1987.

\smallskip \hangindent=7mm \hangafter=1 \noindent  
  [Du88]  \quad   F. Dubois.   Conditions aux limites fortement
non lin\'eaires pour les \'equations d'Euler, in {\it Cours CEA-EDF-INRIA } 
``M\'ethodes de diff\'eren\-ces finies et \'equations hyperboliques'', (P.L. Lions
Editor),  november 1988.

\smallskip \hangindent=7mm \hangafter=1 \noindent  
 [Du90] \quad    F. Dubois. Concavity of thermostatic entropy
and convexity of Lax's mathematical entropy,  {\it Aerospace Research},  n$^{0}$
1990-3, p.~77-80, may 1990.

\smallskip \hangindent=7mm \hangafter=1 \noindent  
 [Du91] \quad     F. Dubois. Nonlinear Interpolation and Total
Variation Diminishing Schemes, {\it  Third International Conference on Hyperbolic
Problems},  (Eng\-quist-Gustafsson Editors), Chartwell-Bratt, p.~351-359, 1991.

\smallskip \hangindent=7mm \hangafter=1 \noindent  
 [Du98] \quad    F. Dubois.   Tout ce que vous avez toujours
voulu savoir sur les conditions aux limites sans jamais oser le demander, {\it
Seminar,  Centre de Math\'e\-matiques Appliqu\'ees, Ecole Polytechnique}, Palaiseau,
juin 1998.

\smallskip \hangindent=7mm \hangafter=1 \noindent  
 [Du99] \quad   F. Dubois.  Flux limite de paroi mobile,
{\it Internal Report}  n$^{0}$ 328/99, CNAM-IAT, october 1999.

\smallskip \hangindent=7mm \hangafter=1 \noindent  
  [DLL91] \quad    F. Dubois, B. Larrouturou, A. Lerat.   
R\'esolution des \'equations d'Euler multidimensionnelles, {\it Cours au DEA de
l'Universit\'e Paris 6}, juin 1991.

\smallskip \hangindent=7mm \hangafter=1 \noindent  
 [DLF87]  \quad     F. Dubois, P. Le Floch.   Condition \`a
la limite pour un syst\`eme de lois de conservation,  {\it Comptes Rendus de
l'Acad\'emie  des Sciences,} Paris, tome 304, S\'erie 1, p.~75-78, 1987. 

\smallskip \hangindent=7mm \hangafter=1 \noindent  
 [DLF88] \quad     F. Dubois, P. Le Floch.   Boundary
Conditions for Nonlinear Hyperbolic Systems of Conservation Laws, {\it Journal of
Differential Equations,} vol.~71, n$^{\rm o}$1, p.~93-122, 1988. 

\smallskip \hangindent=7mm \hangafter=1 \noindent  
 [DLF89] \quad       F. Dubois, P. Le Floch.    Boundary
Conditions for Nonlinear Hyperbolic Systems of Conservation laws, Second International
Conference on Hyperbolic Problems (Ballmann-Jeltsch Editors), {\it Notes on
Numerical Fluid Dynamics,} vol.~24, p.~96-104, Vieweg, 1989.

\smallskip \hangindent=7mm \hangafter=1 \noindent  
 [DM92]   \quad        F. Dubois, O. Michaux.    Solution of 
the Euler Equations Around a Double Ellipso\"{\i}dal  Shape Using Unstructured Meshes
and Including Real Gas Effects, {\it Workshop on Hypersonic Flows for Reentry
Problems}, (D\'esid\'eri-Glowinski-P\'eriaux Editors), Springer Verlag, vol.~II,
p.~358-373, 1992.

\smallskip \hangindent=7mm \hangafter=1 \noindent  
  [FGH91]  \quad       I. Faille, T. Gallou\"et, R. Herbin. 
Les Math\'ematiciens d\'ecouvrent les Volumes Finis, {\it  Matapli,  n$^o$  23, 
p.~37-48, } octobre 1991. 

\smallskip \hangindent=7mm \hangafter=1 \noindent  
  [FL71] \quad     K.O. Friedrichs, P.D. Lax. Systems of
conservation equations with a convex extension, {\it Proc. Nat. Acad. Sci.
U.S.A.}, vol$.\,$68, n$^{\rm o}$8,  p.~1686-1688,  1971.

\smallskip \hangindent=7mm \hangafter=1 \noindent  
 [GB53] \quad     P. Germain, R. Bader. Unicit\'e des
\'ecoulements avec chocs dans la m\'ecanique de Burgers, {\it Note technique ONERA
OA}  n$^{0}$ 1/1711-1, mai 1953. 

\smallskip \hangindent=7mm \hangafter=1 \noindent  
 [Gi94] \quad     M. Gisclon. {\it Conditions aux limites et
approximation parabolique}, PhD thesis, Ecole Normale Sup\'erieure de  Lyon, 
1994.

\smallskip \hangindent=7mm \hangafter=1 \noindent  
 [GR96] \quad    E. Godlewski, P.-A. Raviart. {\it
Numerical Approximation of Hyperbolic Systems of Conservation Laws}, Applied
Mathematical Sciences, vol$.\,$118, Springer, New York,  1996. 

\smallskip \hangindent=7mm \hangafter=1 \noindent  
  [Go59] \quad   S.K. Godunov. A Difference Method for the
Numerical Computation of Discontinuous Solutions of the Equations of Fluid Dynamics,
{\it Math. Sbornik}, vol.~47, p.~271-290, 1959. 

\smallskip \hangindent=7mm \hangafter=1 \noindent  
  [Go61]   \quad   S.K. Godunov. An intersting class of
quasilinear Systems, {\it Doc. Akad. Nauk. SSSR}, vol.~139, p.~521-523 and {\it Soviet
Math.}, vol.~2, p.~947-949, 1961. 

\smallskip \hangindent=7mm \hangafter=1 \noindent  
  [GZIKP79]  \quad     S.K. Godunov, A. Zabrodine, M.
Ivanov, A. Kraiko, G. Prokopov.  {\it  R\'esolution num\'erique des probl\`emes
multidimensionnels de la dynamique des gaz,} Editions de Moscou, 1979.  

\smallskip \hangindent=7mm \hangafter=1 \noindent  
  [Ha83]  \quad     A. Harten.   High Resolution Schemes for
Hyperbolic Conservation Laws, {\it Journal of Computational Physics,} vol.~49, p.
357-393, 1983. 

\smallskip \hangindent=7mm \hangafter=1 \noindent  
 [HLV83] \quad   A. Harten, P.D. Lax, B. Van Leer.  
On Upstream Differencing and Godunov-type Schemes for Hyperbolic Conservation Laws,
{\it  SIAM Review,} vol.~25, n$^{0}$~1, p.~35-61,  January 1983. 

\smallskip \hangindent=7mm \hangafter=1 \noindent  
 [Hi86] \quad    R.L. Higdon. Initial-Boundary Value Problems
for Linear Hyperbolic Systems, {\it SIAM Review}, vol.~28,  n$^{0}$ 2, p.~177-217,
1986. 

\smallskip \hangindent=7mm \hangafter=1 \noindent  
 [Kr70] \quad     H.O. Kreiss. Initial Boundary Value Problems
for Hyperbolic Systems, {\it Comm. Pure Applied Math.}, vol.~23, p.~277-298, 1970. 

\smallskip \hangindent=7mm \hangafter=1 \noindent  
 [Kv70]  \quad    S. Kruzkov. First Order Quasi-Linear Systems
in Several Independent Variables, {\it Math. Sbornik}, vol.~123, p.~228-255, et {\it
Math. USSR Sb.}, vol.~10,  n$^{0}$~2, p.~217-243, 1970. 

\smallskip \hangindent=7mm \hangafter=1 \noindent  
 [LL54]  \quad   L. Landau, E. Lifchitz.  {\it Fluid Mechanics},
1954, Editions de Moscou, 1967.

\smallskip \hangindent=7mm \hangafter=1 \noindent  
  [Lax73]  \quad     P.D. Lax. Hyperbolic Systems of Conservation
Laws and the Mathematical Theory of Shock Waves, {\it Conf. Board in Mathematical
Sciences,} vol.~11, SIAM, Philadelphia, 1973. 

\smallskip \hangindent=7mm \hangafter=1 \noindent  
  [Li77] \quad    T.P. Liu. Initial-Boundary Value Problems
for Gas Dynamics, {\it Archive for Rational Mechanics and Analysis}, vol.~64, p.
137-168, 1977. 

\smallskip \hangindent=7mm \hangafter=1 \noindent  
 [Li82] \quad   T.P. Liu. Transonic Gas Flow in a Duct of
Varying Area,   {\it Archive for Rational Mechanics and Analysis}, vol.~80,  n$^{0}$
1, p.~1-18, 1982.

\smallskip \hangindent=7mm \hangafter=1 \noindent  
 [MO75] \quad     A. Majda, S. Osher. Initial-Boundary Value
Problems for Hyperbolic Equations with Uniformly Characteristic Boundaries, 
{\it Comm. Pure Applied Math.}, vol.~28, p.~607-675, 1975.  

\smallskip \hangindent=7mm \hangafter=1 \noindent  
  [MBGB87] \quad     P. Mazet, F. Bourdel, R. Greborio, 
J. Bor\'ee. Application de la m\'ethode variationnelle d'entropie \`a la
r\'esolution des \'equations d'Euler, {\it ONERA, Centre d'Etudes et de Recherches de
Toulouse, Internal Report}, 1987. 

\smallskip \hangindent=7mm \hangafter=1 \noindent  
 [NS77]  \quad       T. Nishida, J. Smoller. Mixed Problems for
Nonlinear Conservation Laws, {\it Journal of Differential Equations,} vol.~23, p.
244-269, 1977. 

\smallskip \hangindent=7mm \hangafter=1 \noindent  
 [Ol57] \quad       O. Oleinik. Discontinuous Solutions of
Nonlinear Differential Equations, {\it Usp. Mat. Nauk. (N.S.)}, vol.~12, p.~3-73, and
{\it Amer. Math. Transl.}, Ser. 2, vol.~26, p.~95-172, 1957.

\smallskip \hangindent=7mm \hangafter=1 \noindent  
  [OC83] \quad      S. Osher, S. Chakravarthy. Upwind Schemes and
Boundary Conditions with Applications to Euler Equations in General Geometries, {\it
Journal of  Computational Physics}, vol.~50, p.~447-481, 1983.

\smallskip \hangindent=7mm \hangafter=1 \noindent  
  [Pa80] \quad     S.V. Patankar. {\it Numerical Heat Transfer
and Fluid Flow}, Hemisphere publishing, 1980.

\smallskip \hangindent=7mm \hangafter=1 \noindent  
 [Po88]  \quad     M. Pollet.  M\'ethodes de calcul relatives aux
interfaces missiles-propul\-seurs, {\it Internal report}, 
Aerospatiale Les Mureaux, 1988.

\smallskip \hangindent=7mm \hangafter=1 \noindent  
  [Ro72]  \quad     P.J. Roache. {\it Computational Fluid
Dynamics}, Hermosa Publishers, Albukerque, 1972. 

\smallskip \hangindent=7mm \hangafter=1 \noindent  
  [Se96]  \quad   D. Serre. {\it Syst\`emes de lois de
conservation.} vol.~1.  {\it Hyperbolicit\'e, entropies, ondes de choc.} vol.~2.  {\it
structures g\'eom\'etriques, oscillations et probl\`emes mixtes}, Diderot, Paris,
1996. 

\smallskip \hangindent=7mm \hangafter=1 \noindent  
  [Sm83]   \quad   J. Smoller. {\it Shock Waves and
Reaction-Diffusion Equations}, Springer Verlag, Berlin, 1983. 

\smallskip \hangindent=7mm \hangafter=1 \noindent  
   [VL77]   \quad   B. Van Leer. Towards the Ultimate 
Conservative Difference Scheme IV. A New Approach to Numerical Convection, {\it
Journal of Computational Physics,} vol.~23, p.~276-299, 1977.

\smallskip \hangindent=7mm \hangafter=1 \noindent  
   [VL79]   \quad    B. Van Leer.   Towards the Ultimate
Conservative Difference Scheme V. A Second Order Sequel to Godunov's Method,
{\it Journal of Computational Physics,} vol.~32, n$^{\rm o}$ 1, p.~101-136, 1979. 

\smallskip \hangindent=7mm \hangafter=1 \noindent  
  [VV78]  \quad   H. Viviand, J.P.~Veuillot. M\'ethodes
pseudo-instationnaires pour le calcul d'\'ecoulements transsoniques, {\it Publication
ONERA}  n$^{\rm o}$  1978-4, 1978.  

\smallskip \hangindent=7mm \hangafter=1 \noindent  
   [YBW82]   \quad    H. Yee, R. Beam, R. Warming.
Boundary Approximations for Impli\-cit Boundary Approximations for Implicit
Schemes for One-Dimensional Inviscid Equations of Gas Dynamics, {\it AIAA Journal},
vol.~20, p.~1203-1211, 1982.

\bye